\documentclass[11pt]{article}

\usepackage[top=1in, bottom=1in, left=1in, right=1in]{geometry}

\usepackage[linktocpage,colorlinks,linkcolor=blue,anchorcolor=blue,citecolor=blue,urlcolor=blue,pagebackref]{hyperref}
\usepackage{euscript,mdframed}
\usepackage{microtype,todonotes,relsize}
\usepackage{amsmath,amsthm}
\usepackage{algorithm}
\usepackage{algpseudocode}

\usepackage{accents}
\usepackage{comment,url,graphicx,relsize}
\usepackage{amssymb,amsfonts,amsmath,amsthm,amscd,dsfont,mathrsfs,mathtools,nicefrac, bm}
\usepackage{float,psfrag,epsfig,color,xcolor,url}
\usepackage{epstopdf,bbm,mathtools,enumitem}
\usepackage{subcaption}
\usepackage{tablefootnote}
\graphicspath{{Plots/}}
\usepackage{diagbox}
\usepackage{tikz}
\usetikzlibrary{calc,patterns,angles,quotes}
\usepackage{makecell}
\usepackage{multirow}
\usepackage{booktabs}

\ifdefined\nonewproofenvironments\else
\ifdefined\ispres\else

\newenvironment{proof-sketch}{\noindent\textbf{Proof Sketch}
  \hspace*{0.em}}{\qed\bigskip\\}
\newenvironment{proof-idea}{\noindent\textbf{Proof Idea}
  \hspace*{0.em}}{\qed\bigskip\\}
\newenvironment{proof-of-lemma}[1][{}]{\noindent\textbf{Proof of Lemma {#1}.}
  \hspace*{0.em}}{\qed\\}
\newenvironment{proof-of-corollary}[1][{}]{\noindent\textbf{Proof of Corollary {#1}.}
  \hspace*{0.em}}{\qed\\}
\newenvironment{proof-of-theorem}[1][{}]{\noindent\textbf{Proof of Theorem {#1}.}
  \hspace*{0.em}}{\qed\\}
\newenvironment{proof-attempt}{\noindent\textbf{Proof Attempt}
  \hspace*{0.em}}{\qed\bigskip\\}

\fi
\newtheorem{theorem}{Theorem}[section]
\newtheorem{lemma}{Lemma}[section]
\newtheorem{corollary}{Corollary}[section]
\newtheorem{proposition}{Proposition}[section]
\newtheorem{assumption}{Assumption}[section]
\newtheorem{remark}{Remark}[section]

\usetikzlibrary{calc}

\allowdisplaybreaks
\usepackage[numbers,round]{natbib}
\renewcommand*{\backref}[1]{\ifx#1\relax \else Page #1 \fi}
\renewcommand*{\backrefalt}[4]{%
  \ifcase #1 \footnotesize{(Not cited.)}%
  \or        \footnotesize{(Cited on page~#2.)}%
  \else      \footnotesize{(Cited on pages~#2.)}%
  \fi
}

\newcommand*{\colorboxed}{}
\def\colorboxed#1#{%
  \colorboxedAux{#1}%
}
\newcommand*{\colorboxedAux}[3]{%
  \begingroup
    \colorlet{cb@saved}{.}%
    \color#1{#2}%
    \boxed{%
      \color{cb@saved}%
      #3%
    }%
  \endgroup
}
\numberwithin{equation}{section}
\usepackage{xcolor} 
\usepackage{marginnote}

\newcommand{\todol}[2][]{{%
 \let\marginpar\marginnote
 \reversemarginpar
 \renewcommand{\baselinestretch}{0.8}%
 \todo[color=yellow]{#2}}}

\begin{document}

\title{An Augmented Lagrangian Value Function Method for Lower-level Constrained Stochastic Bilevel Optimization}

\author{
Hantao Nie\footnote{School of Mathematical Science, Peking University, Beijing, China. nht@pku.edu.cn} \and 
Jiaxiang Li\footnote{Department of Electrical and Computer Engineering, University of Minnesota, Minneapolis, MN, USA. li003755@umn.edu} \and 
Zaiwen Wen\footnote{School of Mathematical Science, Peking University, Beijing, China. wenzw@pku.edu.cn}
}

\date{}

\maketitle

\begin{abstract}
   Recently, lower-level constrained bilevel optimization has attracted increasing attention. However, existing methods mostly focus on either deterministic cases or problems with linear constraints. The main challenge in stochastic cases with general constraints is the bias and variance of the hyper-gradient, arising from the inexact solution of the lower-level problem. 
   In this paper, we propose a novel stochastic augmented Lagrangian value function method for solving stochastic bilevel optimization problems with nonlinear lower-level constraints. Our approach reformulates the original bilevel problem using an augmented Lagrangian-based value function and then applies a penalized stochastic gradient method that carefully manages the noise from stochastic oracles. We establish an equivalence between the stochastic single-level reformulation and the original constrained bilevel problem and provide a non-asymptotic rate of convergence for the proposed method. The rate is further enhanced by employing variance reduction techniques. Extensive experiments on synthetic problems and real-world applications demonstrate the effectiveness of our approach. 
\end{abstract}

\section{Introduction}
We consider the stochastic lower-level constrained bilevel optimization (stochastic LC-BLO) problem. The lower-level problem is defined as
\begin{equation}
\label{eq: LL}
   \begin{aligned}
    \min_{y \in Y} & \quad  G(x, y) = \mathbb{E}_{\xi \sim \mathcal{D}_{\xi}}[g(x, y; \xi)]  \\
    \text{s.t.} &\quad 
    H_i(x, y)   \leq 0, \quad i = 1, ..., p,
\end{aligned}
\end{equation}
where $Y \subseteq \mathbb{R}^n$ is a convex compact set, $\xi$ is a random variable in the space $\Omega_{\xi}$, and $\mathcal{D}_{\xi}$ is the distribution of $\xi$.
$G(x, y)$ and $H(x, y) = [H_1(x, y), ..., H_{p}(x, y)]^T$ are the lower-level objective function and constraint function, respectively. We also denote the feasible set of the lower-level problem by $\mathcal{Y}(x) := \{y \in Y \mid H(x, y) \leq 0\}$.
The bilevel optimization (BLO) problem is
\begin{equation}
\label{eq: BLO}
\begin{aligned}
   \min_{x \in X} &\quad  F(x, y^*(x)) = \mathbb{E}_{\zeta \sim \mathcal{D}_{\zeta}}[f(x, y^*(x); \zeta)]  \\
   \text{s.t.} &\quad  y^*(x) \in \arg \min_{y \in \mathcal{Y}(x)} G(x, y),
\end{aligned}
\end{equation}
where $X \subseteq \mathbb{R}^m$ is a convex compact set, $\zeta \in \Omega_{\zeta}$ is a random variable, and $\mathcal{D}_{\zeta}$ is the distribution of $\zeta$. 
We assume the gradient oracles $\nabla g(x, y; \xi)$ and $\nabla f(x, y; \zeta)$ have unavoidable noise. This framework includes deterministic LC-BLO as a special case.

As depicted in~\eqref{eq: BLO}, this hierarchical structure captures a learning‑to‑learn philosophy that underpins numerous modern machine‑learning pipelines. Hence,
BLO plays a critical role in various machine learning tasks, including hyperparameter optimization~\cite{bennett2006model, franceschi2018bilevel, mackay2019self, sinha2020gradient}, model-agnostic meta-learning~\cite{franceschi2018bilevel,ji2020provably, qin2023bi, zhu2020personalized}, and reinforcement learning~\cite{stadie2020learning,shen2024principled,yang2024bilevel, kang2023bi}. Recently, the lower-level constrained BLO (LC-BLO) has attracted increasing attention due to its wide applications such as transportation~\cite{jiang2024primal}, kernelized SVM~\cite{yao2024constrained}, meta-learning~\cite{xu2023efficient}, and data hyper-cleaning
~\cite{xu2023efficient, yao2024overcoming}.

Several methods have been proposed to solve the deterministic LC-BLO problem. The methodologies for solving LC-BLO can be broadly categorized into implicit gradient-based (IG) approaches and lower-level value function-based (LLVF) techniques. Implicit gradient methods primarily focus on computing the hyper-gradient
with some implicit gradient approximation of $\frac{d}{dx} y^*(x)$. The key issue here is that $y^*(x)$ may not be differentiable when the lower-level problem has constraints. Some research has discussed the conditions under which the hyper-gradient exists. For the linearly constrained case, IG-AL~\cite{tsaknakis2022implicit} discussed the smoothness of $y^*(x)$ by introducing the Lagrangian multiplier of the lower-level problem. SIGD~\cite{khanduri2023linearly} develops a smoothing approximation implicit gradient method to handle the non-differentiability point. Recent works~\cite{tsaknakis2023implicit,jiang2024barrier} also consider a barrier reformulation of LC-BLO to write the lower-level constraints into the objective function. Even with these existing discussions, the main challenge of designing implicit gradient methods is that second-order derivative of the lower-level objective is required to compute the hyper-gradient, which is computationally expensive.

To overcome the high computational cost of implicit gradient computation, penalty methods and value function-based methods have been proposed as a Hessian-free alternative. In~\cite{lu2024first}, the authors demonstrate that the LC-BLO can be approximately soved as minimax optimization through first-order penalty methods.
Value function-based methods introduce the value function of the lower-level problem and then replace the optimality condition $y \in \mathcal{Y}(x)$ in~\eqref{eq: BLO} with an inequality condition on the value function. The original BLO~\eqref{eq: BLO} is then equivalently reformulated as a single-level problem. Different value-function formulations have been widely studied in the literature, including the value function~\cite{jiang2024primal}, a Moreau envelope-based value function~\cite{gao2023moreau}, the proximal Lagrangian function~\cite{yao2024constrained}, and the regularized gap function~\cite{yao2024overcoming}.

To the best of our knowledge, no existing work analyzes the BLO with nonlinear lower-level constraints in the stochastic setting, which has broad applications in real-world machine learning tasks. Solving stochastic LC-BLO presents two fundamental challenges. The first is the nonsmoothness of the hyper-objective function, which arises from the coupling between the upper-level and lower-level problems. The second is the bias of the hyper-gradient due to the inexact solution of the lower-level problem. 
To address these challenges, we introduce an augmented Lagrangian function and its Moreau envelope to reformulate the bilevel problem as a single-level problem, while ensuring that the solution remains close to the optimal solution of the original problem. We then propose a novel stochastic value function-based method for stochastic LC-BLO, which carefully controls the bias of the gradient oracle to achieve convergence.

\subsection{Main contribution}
Our main contributions are as follows:

\noindent
1.
We introduce a novel reformulation of the stochastic LC-BLO by leveraging the stochastic augmented Lagrangian function and its Moreau envelope (see~\eqref{eq: augmented dual function}). This reformulation transforms the bilevel problem into a single-level problem, effectively addressing noise arising from the inexact solution of the lower-level problem (see~\eqref{eq: single level},~\eqref{eq: stochastic single level approximation}). Notably, we also ensure that the solution of the reformulated problem remains close to the optimal solution of the original bilevel problem (see Theorem~\ref{thm: maintext equivalence of single level},~\ref{thm: maintext equivalence between BLO and stochastic single level approximation}), providing a practical yet theoretically grounded approach to stochastic LC-BLO.

\noindent
2. 
We propose a novel Hessian-free method based on the stochastic reformulation for solving the stochastic LC-BLO (see Algorithm~\ref{alg: main}). 
Our work provides the first convergence analysis of value function-based algorithms for nonlinear LC-BLO in the stochastic setting.
The issue of biased gradients is mitigated by controlling the bias through the accuracy of lower-level solutions.
We derive a non-asymptotic convergence rate, proving that our method achieves $(\widetilde{O}(c\epsilon^{-2}), \widetilde{O}(c c_1^2\epsilon^{-2}))$ sample complexity on $(\zeta, \xi)$, where $c_1, c_2$ denotes the penalty parameter in the reformulations and $c = \max(c_1, c_2)$ (see Theorem~\ref{thm: maintext convergence of SALVF}, Remark~\ref{rem: maintext convergence of SALVF}). The sample complexity on $\zeta$ is further improved to $\widetilde{O}(c^{1.5}\epsilon^{-1.5})$ by employing variance reduction techniques (see Theorem~\ref{thm: maintext convergence of SALVF-VR}, Remark~\ref{rem: maintext convergence of SALVF-VR}).
In Table~\ref{tab: comparison}, we briefly summarize the existing approaches compared to the proposed methods.

\begin{table}[t]
    \caption{Comparison of methods for solving BLO with nonlinear convex lower-level constraints. 
    ``H.-free'', ``Sto'', ``Iteration'', ``Sample''
    stand for Hessian-free, stochastic, iteration complexity, and sample complexity, respectively. Iteration complexity is the number of outer-loop iterations needed to achieve the target accuracy $\epsilon$, measured by the squared gradient norm. For stochastic algorithms, sample complexities denote the number of samples for the upper and lower stochastic variables $(\zeta, \xi)$ needed to achieve the target accuracy $\epsilon$.
   The proposed methods are marked bold as \textbf{SALVF} (Algorithm~\ref{alg: main}) and \textbf{SALVF-VR} (Algorithm~\ref{alg: variance reduce}). We compare with four existing value function-based deterministic methods: BLOCC~\cite{jiang2024primal}, LV-HBA~\cite{yao2024constrained}, BiC-GAFFA~\cite{yao2024overcoming}, and GAM~\cite{xu2023efficient}.
   } 
   \centering
   \renewcommand{\arraystretch}{1.2} 
   \setlength{\tabcolsep}{4pt} 
   \resizebox{0.7\columnwidth}{!}{
   \begin{tabular}{c|c|c|c|c}
   \hline 
    & H-free & Sto & Iteration & Sample (upper, lower)  \\
   \hline
   \hline
   BLOCC & yes & no & $\tilde{\mathcal{O}}(c\epsilon^{-1})$ & --\\
   \hline
   LV-HBA & yes & no & $\mathcal{O}(\epsilon^{-2p}), p>0.5$ & --\\
   \hline
   BiC-GAFFA & yes & no & $\mathcal{O}(\epsilon^{-2p}), p>0.5$ & --\\
   \hline
   GAM & no & no & -- & -- \\
   \hline 
   \textbf{SALVF}  & yes & yes & $\tilde{\mathcal{O}}(c \epsilon^{-1})$ & $\tilde{\mathcal{O}}(c\epsilon^{-2}), \tilde{\mathcal{O}}(c c_1^2 \epsilon^{-2})$ \\
   \hline
   \textbf{SALVF-VR} & yes & yes & $\tilde{\mathcal{O}}(c^{1.5}\epsilon^{-1.5}) $& $\tilde{\mathcal{O}}(c^{1.5}\epsilon^{-1.5}), \tilde{\mathcal{O}}(c^{1.5} c_1^{2}\epsilon^{-2.5})$ \\
   \hline
   \end{tabular}
   }
   \label{tab: comparison}
\end{table}

\subsection{Notation}
For multivariate function $f(x_1, ..., x_k)$, its partial derivative with respect to the $i$-th variable is denoted by $\nabla_i f(x_1, ..., x_k)$. 
Given random variable $x(\xi) \in \mathbb{R}^n, \xi \sim \mathcal{D}_{\xi}$,
we represent its expectation and covariance matrix by $\mathbb{E}_{\xi \sim \mathcal{D}_{\xi}}[x(\xi)]$ and $\mathrm{Var}_{\xi \sim \mathcal{D}_{\xi}}[x(\xi)]$, respectively. 
The trace of its covariance matrix is denoted by
$
   \mathbb{V}_{\xi \sim \mathcal{D}_{\xi}}[x(\xi)] = \mathrm{Tr} (\mathrm{Var}_{\xi \sim \mathcal{D}_{\xi}}[x(\xi)]) 
   = \mathbb{E}_{\xi \sim \mathcal{D}_{\xi}}[ \| x - \mathbb{E}_{\xi \sim \mathcal{D}_{\xi}} [x] \|^2 ].
$
If the distribution $\mathcal{D}_{\xi}$ is clear, we abbreviate the above notations as $\mathbb{E}_{\xi} [x(\xi)]$, $\mathrm{Var}_{\xi}[x(\xi)]$ and $\mathbb{V}_{\xi}[x(\xi)]$, respectively.
Denote the normal cone of a convex set $C$ at $y$ as $\mathcal{N}_C(y) = \{v \mid \langle v, w - y \rangle \leq 0, \forall w \in C\}$.
For a scalar $a$, we define $[a]_+ = \max\{0, a\}$. For a vector $v$, we write $[v]_+ = ([v_1]_+, ..., [v_n]_+)^T$
and $[v]_+^2 = ([v_1]_+^2, ..., [v_n]_+^2)^T$.
For $s$ independent identically distributed random variables $\mathbf{\xi} = (\xi_1, ..., \xi_s)$, we denote the joint distribution by $\mathcal{D}_{\mathbf{\xi}}^s$. Further given some function $g(x, \xi)$, we define the empirical average as
$ g(x, \mathbf{\xi}) = \frac{1}{s} \sum_{i=1}^{s} g(x, \xi_i)$.




\section{Stochastic augmented value function-based method}

In this section, we propose a novel value function-based reformulation for~\eqref{eq: BLO} and a stochastic value function method. 

\subsection{Stochastic augmented Lagrangian reformulation}
We propose a novel reformulation using a stochastic augmented Lagrangian value function that transforms~\eqref{eq: BLO} into a single level problem. 
First, the constraints \eqref{eq: LL} are addressed by the augmented Lagrangian function and its corresponding stochastic version. The bilevel optimization problem~\eqref{eq: BLO} is then transformed into a single-level problem via this augmented Lagrangian function-based formulation. 

For the lower-level problem~\eqref{eq: LL}, an augmented Lagrangian penalty term is introduced by penalizing the constraints using
$$
\mathcal{A}_{\gamma_1}(x, y, z) = \frac{1}{2{{\gamma_1}}} \sum_{i=1}^{p} [{\gamma_1} z_i + H_{i}(x, y)]_+^2, 
$$
where $z_{1}, ..., z_{p} \geq 0$ and ${\gamma_1} > 0$ is a penalty parameter. The augmented Lagrangian function and its stochastic oracle are defined by adding the penalty term to the objective function and its stochastic oracle, respectively, that is,
\begin{equation}
   \label{eq: AL}
   \begin{aligned}
    \mathcal{L}_{{\gamma_1}}(x, y, z) &= G(x, y) 
    + \mathcal{A}_{\gamma_1}(x, y, z), \\
      \mathcal{L}_{{\gamma_1}}(x, y, z; \xi) &= g(x, y; \xi) 
    + \mathcal{A}_{\gamma_1}(x, y, z).
\end{aligned}
\end{equation}
The augmented dual function and its Moreau envelope are then defined as
\begin{equation}
   \label{eq: augmented dual function}
   \begin{aligned}
   {D}_{\gamma_1}(x, z) &= \min_{y \in Y} \mathcal{L}_{\gamma_1}(x, y, z),
  \\ 
   {E}_{\gamma_1}^{\gamma_2}(x, z) &= \max_{\lambda \in \mathbb{R}_+^p} \{ {D}_{\gamma_1}(x, \lambda) - \frac{\gamma_2}{2} \| \lambda - z\|^2 \},
   \end{aligned}
\end{equation}
where $\gamma_2 \geq 0$ is a regularization parameter. Then~\eqref{eq: BLO} is reformulated as an equivalent single-level problem
\begin{equation}
   \label{eq: single level}
\begin{aligned}
   \min_{(x, y, z) \in X \times Y \times \mathbb{R}_+^p} &  F(x, y)  \\
   \text{s.t.} \quad \quad \quad &  \mathcal{G}(x, y, z) := G(x, y) -{E}_{\gamma_1}^{\gamma_2}(x, z) \leq 0, \\
   & \quad H(x, y) \leq 0.
\end{aligned}
\end{equation}
In this paper, $\gamma_1, \gamma_2$ are fixed parameters. We omit the subscript ${\gamma_1}$ and the superscript ${\gamma_2}$ in $D_{\gamma_1}$ and $E_{\gamma_1}^{\gamma_2}$ to simplify notation. 
The envelope-based value function reformulation~\eqref{eq: single level} contains the value function-based reformulation 
\begin{equation}
\begin{aligned}
   \min_{(x, y) \in X \times Y} &\quad  F(x, y)  \\
   \text{s.t.} \quad  &\quad G(x, y) - \min_{y \in \mathcal{Y}(x)} G(x, y) \leq 0, \\
   &\quad  H(x, y) \leq 0, 
\end{aligned}
\end{equation}
as a special case when $\gamma_2 = 0$. The relationship between them is discussed in Appendix~\ref{sec: reformulation}. 
By introducing the auxiliary variable $z$, the advantage of the augmented function-based reformulation is that the subproblem~\eqref{eq: minimax subproblem} becomes strongly-convex-strongly-concave, which ensures faster convergence of the inner loop.

To evaluate the $E(x, z)$ in~\eqref{eq: augmented dual function}, we need to estimate 
\begin{subequations}
    \begin{align}
    (w^*, \lambda^*) &= \arg 
    \max_{\lambda \in \mathbb{R}_+^p}\min_{w \in Y}
     \left\{ \ell_{\gamma}(x, z, w, \lambda) \right\},
  \label{eq: minimax subproblem} \\
   \text{where}\quad\ell_{\gamma}(x, z, w, \lambda) &:= \mathcal{L}_{\gamma_1}(x, w, \lambda) - \frac{\gamma_2}{2} \| \lambda - z \|^2. \label{eq: tilde L}
   \end{align}
\end{subequations}
However, in the stochastic setting, the exact solution is inaccessible due to unavoidable noise in the gradient oracles.
To address this, we consider approximating the solution of~\eqref{eq: minimax subproblem} using stochastic algorithms. Specifically, let $\mathbf{\xi} = (\xi_1, ..., \xi_{s}) \in \Omega_{\xi}^s$ be the samples from $\mathcal{D}_{\xi}^{s}$. 
Denote $\mathcal{P}_{w}, \mathcal{P}_{\lambda}$ as the space of random variables mapping $\mathbf{\xi}$ to $Y$ and $\mathbb{R}_+^p$, respectively, that is, 
$$
\begin{aligned}
   \mathcal{P}_{w} &= \{\hat w: \Omega_{\xi}^s \to Y \mid \hat w \text{ is measurable}\}, \\ 
   \mathcal{P}_{\lambda} &= \{\hat \lambda: \Omega_{\xi}^s \to \mathbb{R}_+^p \mid \hat \lambda \text{ is measurable}\}.
\end{aligned}
$$
Assume $(\hat w, \hat \lambda) \in \mathcal{P}_{w} \times \mathcal{P}_{\lambda}$ is a stochastic algorithm solving the subproblem~\eqref{eq: minimax subproblem} using samples $\mathbf{\xi}$ and $(\hat w( \mathbf{\xi}), \hat \lambda(\mathbf{\xi}))$ is a pair of approximate solution of~\eqref{eq: minimax subproblem}. We are interested in the subset of ``good enough'' algorithms that can provide a sufficiently accurate solution of the subproblem~\eqref{eq: minimax subproblem} as 
\begin{equation}
   \label{eq: P delta}
\begin{aligned}
   \mathcal{P}(\delta) &= \left\{(\hat w, \hat \lambda) \in \mathcal{P}_{w} \times \mathcal{P}_{\lambda} \mid 
   \left| \mathbb{E}_{\mathbf{\xi}} \left[  \ell_{\gamma}(x, z, \hat w(\mathbf{\xi}), \hat \lambda(\mathbf{\xi}))  \right]  - E(x, z) \right| \leq \delta \right\}.
\end{aligned}
\end{equation}

With this estimator $(\hat w(\mathbf{\xi}), \hat \lambda(\mathbf{\xi}))$, we approximate the envelope function with
$\ell_{\gamma}(x, z, \hat w(\mathbf{\xi}), \hat \lambda(\mathbf{\xi}))$.
This gives the following stochastic value function-based reformulation:
\begin{equation}
   \label{eq: stochastic single level approximation}
\begin{aligned}
   \min_{
(x, y, z) \in X \times Y \times \mathbb{R}_+^p } &\quad F(x, y) \\
   \text{s.t.}\quad  \quad  &
    \hat{\mathcal{G}} (x, y, z; \mathbf{\xi}) \leq \epsilon_1,   \quad \frac{1}{2} \sum_{i=1}^p [H_i(x, y)]_+^2 \leq \epsilon_2^2, 
\end{aligned}
\end{equation}
where $\hat{\mathcal{G}}(x, y, z; \mathbf{\xi}) = G(x, y) - \ell_{\gamma}(x, z, \hat w(\mathbf{\xi}), \hat \lambda(\mathbf{\xi}))$, $\epsilon_1, \epsilon_2$ are the target accuracy of the lower-level objective function and constraints violation, respectively. The equivalence between~\eqref{eq: BLO} and~\eqref{eq: stochastic single level approximation} is established in Theorem~\ref{thm: maintext equivalence between BLO and stochastic single level approximation}. Compared to~\eqref{eq: single level}, ~\eqref{eq: stochastic single level approximation}  incorporates inexactness of the lower-level solution into the formulation, making it more practical in the stochastic setting.

\begin{algorithm}[ht]
   \caption{ {\fontsize{8}{11}\selectfont $(w^k, \lambda^k) = \text{SALM}(x^{k-1}, z^{k-1}, s, \gamma_1, \gamma_2,  \eta, \rho; \mathbf{\xi}^k)$}}
   \label{alg: ALM}
   \begin{algorithmic}[1]
   \State Input: $x^{k-1}, z^{k-1}$, iteration count $s$, primal step size $\eta_j$, dual step size $\rho_j$ for $0 \leq j \leq s-1$.
   \State Initialize $w^{k, 0}$ and $ \lambda^{k, 0} = 0$.
   \For{$j=0$ {\bfseries to} $s-1$}
   \State Update $(w^{k, j+1}, \lambda^{k, j+1})$ by~\eqref{eq: inner algorithm primal step} and~\eqref{eq: inner algorithm dual step}.
   \EndFor
   \State Output $(w^k, \lambda^k) = (w^{k, s}, \lambda^{k, s})$.
   \end{algorithmic}
\end{algorithm}


\subsection{Value function-based penalized method}
In this subsection, we develop a stochastic value function-based penalized method for the reformulation~\eqref{eq: stochastic single level approximation}. 
At $k$-th iteration, the stochastic gradient ascent descent method is applied to solve 
an approximation solution from the subproblem~\eqref{eq: minimax subproblem} with fixed $(x^{k-1}, z^{k-1})$.
More specifically, the primal and dual variables are updated by 
   \begin{subequations}
      \begin{align}
        w^{k, j+1} &= 
      \mathrm{Proj}_{Y}\left(w^{k, j} - \eta_{j}\nabla_w \ell_{\gamma}(x^{k-1}, z^{k-1}, w^{k, j}, \lambda^{k, j}; \xi^k_j) \right) , \label{eq: inner algorithm primal step} \\
  \lambda^{k, j+1} &= \mathrm{Proj}_{\mathbb{R}_+^p}\left(\lambda^{k, j} + \rho_{j} \nabla_{\lambda} \ell_{\gamma}(x^{k-1}, z^{k-1}, w^j, \lambda^{k, j}; \xi^k_j) \right)\label{eq: inner algorithm dual step},
  \end{align}
  \end{subequations}
where $\eta_{p}^j, \eta_{d}^j$ denote the primal and dual step sizes, respectively. The complete procedure is shown in Algorithm~\ref{alg: ALM}.
After that we consider a augmented Lagrangian-based penalty reformulation:
\begin{equation}
   \label{eq: maintext stochastic penalized single level}
\begin{aligned}
   \min_{(x, y, z) \in X \times Y \times Z} &\quad \mathbb{E}_{\mathbf{\xi}} [\Psi(x, y, z; \mathbf{\xi})] , \\
   \text{where}\quad\Psi(x, y, z; \mathbf{\xi}) &:= F(x, y)  + c_1 \hat{\mathcal{G}}(x, y, z; \mathbf{\xi}) + \frac{c_2}{2} \sum_{i=1}^{p} [H_i(x, y)]_+^2,
\end{aligned}
\end{equation}
where
$Z = [0, p^{-0.5} B]^p$ is the domain of $z$, $B$ is a constant and $c_1, c_2 > 0$ are the penalty parameters\footnote{It can be shown that the optimal $z$ for~\eqref{eq: maintext stochastic penalized single level} is contained in the domain $Z$ with appropriate regularity assumptions (see Assumption~\ref{ass: LICQ} and Lemma~\ref{lem: equivalence between the  value function-based reformulation and the dual augmented value function reformulation})}.
Different from the standard duality, the penalty parameter $\rho$ is also treated as a variable in the above saddle reformulation. The advantage of this saddle point refomulation is that the convexity of objective functin in~\eqref{eq: stochastic single level approximation} is not required.
The equivalence between~\eqref{eq: stochastic single level approximation} and~\eqref{eq: maintext stochastic penalized single level} is established in Theorem~\ref{thm: maintext equivalence of single level}.
We consider a stochastic gradient descent ascent method for solving~\eqref{eq: maintext stochastic penalized single level}. 
Denote all variables in~\eqref{eq: maintext stochastic penalized single level} and its region as
$$
\mathbf{u} = (x, y, z), \quad \mathcal{U} = X \times Y \times Z,
$$
for simplicity.
The gradient oracle of the objective function of~\eqref{eq: maintext stochastic penalized single level} is given by
\begin{equation}
   \label{eq: stochastic gradient of Psi}
   \begin{aligned}
   \nabla \Psi(\mathbf{u}; \mathbf{\zeta}, \mathbf{\xi}, \mathbf{\tilde{\xi}}) &= \nabla f(x, y; \mathbf{\zeta}) + c_1 \nabla \hat{\mathcal{G}}(\mathbf{u}; \mathbf{\tilde{\xi}})  
    + c_2 \sum_{i=1}^{p} [H_i(x, y)]_+  \nabla H_i(x, y),
\end{aligned}
\end{equation}
with the stochastic oracles in mini-batched as
\begin{equation}
   \label{eq: stochastic gradient of Psi mini-batched}
    \begin{aligned}
   \nabla f(x, y; \mathbf{\zeta}) &= \frac{1}{r} \sum_{j=1}^{r} \nabla f(x, y; \zeta_j), \\
   \nabla \hat{\mathcal{G}}(\mathbf{u}; \mathbf{\xi}, \mathbf{\tilde{\xi}}) &=\frac{1}{q} \sum_{j=1}^{q} \nabla g(x, y; \tilde \xi_j)- \frac{1}{q} \sum_{j=1}^{q} \nabla \tilde{\mathcal{L}} (x, z, \hat w(\mathbf{\xi}), \hat \lambda(\mathbf{\xi}); \tilde \xi_j).
\end{aligned}
\end{equation}
By substituting $r= r_k, q = q_k$ in~\eqref{eq: stochastic gradient of Psi},~\eqref{eq: stochastic gradient of Psi mini-batched} and conditioned on $\mathbf{\xi}^k$ we rewrite
\begin{equation}
   \label{eq: Psi k}
   \begin{aligned}
   \hat{\mathcal{G}}^k(\mathbf{u}) = \hat{\mathcal{G}}(x, y, &z; \mathbf{\xi}^k), \quad 
\Psi^k(\mathbf{u}) = \Psi(x, y, z; \mathbf{\xi}^k), \\
\nabla \Psi^k(\mathbf{u}; \mathbf{\zeta}^k, \mathbf{\tilde{\xi}}^k) &= \nabla \Psi(x, y, z; \mathbf{\zeta}^k, \mathbf{\xi}^k, \mathbf{\tilde{\xi}}^k).
\end{aligned}
\end{equation}

The complete procedure is summarized as follows. At the $k$-th iteration, an estimator $(w^k, \lambda^k) = (\hat w(x^{k-1}, z^{k-1}; \mathbf{\xi}^k), \hat \lambda(x^{k-1}, z^{k-1}; \mathbf{\xi}^k))$ is computed utilizing Algorithm~\ref{alg: ALM} with samples $\mathbf{\xi}^k = (\xi_1^k, ..., \xi_{s_k}^k) \sim \mathcal{D}_{\xi}^{s_k}$. 
Then the stochastic gradient descent method is applied to solve~\eqref{eq: maintext stochastic penalized single level} with samples $\mathbf{\zeta}^k = (\zeta_1^k, ..., \zeta_{r_k}^k) \sim \mathcal{D}_{\zeta}^{r_k}$, $\mathbf{\tilde{\xi}}^k = (\tilde \xi_1^k, ..., \tilde \xi_{q_k}^k) \sim \mathcal{D}_{\xi}^{q_k}$. 
After $K$ iterations, we output $(x^R, y^R)$ where the index $R$ is randomly chosen according to the probability mass function
\begin{equation}
   \label{eq: probability mass function}
   \mathrm{Prob}(R=k) = \frac{\alpha_k}{\sum_{k=0}^{K-1} \alpha_k}, \quad k=0, \ldots, K-1 .
\end{equation}
Besides, an extra SALM loop 
\begin{equation}
   \label{eq: feasible projection}
(y', z') = \textbf{SALM}(x^R, 0, s^K, {\gamma_1}, 0, \eta, \rho; \mathbf{\xi}^K),
\end{equation}
can be applied to guarantee the feasibility of the final output.
The complete procedure is shown in Algorithm~\ref{alg: main}.


\begin{algorithm}[ht]
   \caption{SALVF}
   \label{alg: main}
   \begin{algorithmic}[1]
   \State Input penalty parameters $c_1, c_2$, iteration number $K$, sample sizes $s_k, r_k, q_k$ and step sizes $\eta_{p}^k, \eta_{d}^k, \alpha_k$.
   \State Initialize $x^0, y^0, z^0$.
   \For{$k=0$ {\bfseries to} $K-1$}
      \State Run 
      Algorithm~\ref{alg: ALM} with samples $\mathbf{\xi}^k$ to compute $$(\hat w^{k}, \hat \lambda^{k})
      = \text{SALM}(x^{k-1}, z^{k-1}, s_k, {\gamma_1}, \gamma_2, \eta_k, \rho_k; \mathbf{\xi}^k).$$ 
      \State Sample $\mathbf{\zeta}^k \sim \mathcal{D}_{\zeta}^{r_k}$ and $\mathbf{\tilde{\xi}}^k \sim \mathcal{D}_{\xi}^{q_k}$.\;
      \State Compute direction
      $
      d^k = \nabla \Psi^k(\mathbf{u}^k; \mathbf{\zeta}^k, \mathbf{\tilde{\xi}}^k) 
      $ by~\eqref{eq: stochastic gradient of Psi}.
      \\
      \State Update 
      $
      \mathbf{u}^{k+1} = \mathrm{Proj}_{\mathcal{U}}(\mathbf{u}^k - \alpha_k d^k )$.
      
   \EndFor
   \State Choose index $R$ with probability mass function~\eqref{eq: probability mass function}.
   Output $(x^R, y^R)$.
   \State (Optional) Compute~\eqref{eq: feasible projection} and
   Output $(x^R, y')$.
   \end{algorithmic}
\end{algorithm}


\subsection{Variance reduced SALVF method}
When sampling on $\zeta$  is significantly more expensive than sampling on $\xi$,  a natural question arises:  is it possible to reduce the sample size of  $\zeta$, thereby allowing an increase in the sample size of  $\xi$?
In this subsection, we apply a variance reduction technique with a similar update rule to STORM~\cite{cutkosky2019momentum} to reduce the sample complexity on $\zeta$.
At each iteration, the direction $d^k$ is updated by
\begin{equation}
   \label{eq: direction update in variance reduction}
   d^{k} =\nabla \Psi^k(\mathbf{u}^{k}; \mathbf{\zeta}^{k}, \mathbf{\tilde{\xi}}^k) + (1 - \beta_{k})(d^{k-1} -\nabla \Psi^k(\mathbf{u}^{k-1}; \mathbf{\zeta}^{k}, \mathbf{\tilde{\xi}}^k)).
   \end{equation}
Unlike the STORM method, our approach deals with a different scenario, where the main challenge arises from the biased gradient oracle $\nabla \Psi^k$ due to the inexactness of estimator $(\hat w^k, \hat \lambda^k)$. This challenge is addressed by
carefully handling the extra bias term and designing proper coefficients $\beta_k$.
The complete procedure is summarized in Algorithm~\ref{alg: variance reduce} and the convergence guarantee is provided in Theorem~\ref{thm: maintext convergence of SALVF-VR}.

\begin{algorithm}[ht]
   \caption{SALVF-VR}
   \label{alg: variance reduce}
   \begin{algorithmic}[1]
    \State Input penalty parameters $c_1, c_2$, iteration number $K$, sample sizes $s_k, r_k, q_k$, step sizes $\eta_{p}^k, \eta_{d}^k, \alpha_k$ and coefficient $\beta_k$.
   \State Initialize $x^0, y^0, z^0$.\;
   \For{$k=0$ {\bfseries to} $K-1$}
      \State Run 
      Algorithm~\ref{alg: ALM} with samples $\mathbf{\xi}^k$ to compute $$(\hat w^{k}, \hat \lambda^{k})
      = \text{SALM}(x^{k-1}, z^{k-1}, s_k, {\gamma_1}, \gamma_2, \eta_k, \rho_k; \mathbf{\xi}^k).$$ 
      \State Sample $\mathbf{\zeta}^k$ from $\mathcal{D}_{\zeta}^{r_k}$ and $\mathbf{\tilde{\xi}}^k$ from $\mathcal{D}_{\xi}^{q_k}$.

      \If{$k=0$}
        \State Compute 
      $
      d^0 = \nabla \Psi(\mathbf{u}^0; \mathbf{\zeta^0}, \mathbf{\tilde{\xi}}^0) .
      $
      \Else
      \State Update the direction by~\eqref{eq: direction update in variance reduction}.
      \EndIf
      \State Update 
      $
      \mathbf{u}^{k+1} = \mathrm{Proj}_{\mathcal{U}}(\mathbf{u}^k - \alpha_k d^k ) .
      $
   \EndFor
   \State Choose index $R$ with probability mass function ~\eqref{eq: probability mass function}.
   Output $(x^R, y^R)$.
   \State (Optional) Compute~\eqref{eq: feasible projection} and
   Output $(x^R, y')$.
   \end{algorithmic}
\end{algorithm}



\section{Theoretical analysis}
\subsection{Basic assumptions}
In this section, we inspect the properties of the penalty reformulated problem~\eqref{eq: maintext stochastic penalized single level}, also provide non-asymptotic convergence guarantee for the proposed algorithms (Algorithm~\ref{alg: main} and~\ref{alg: variance reduce}).
First some basis assumptions in the literature of stochastic bilevel optimization~\cite{hong2023two,chen2021tighter} are introduced, including assumptions on the smoothness, convexity, boundedness, and the stochastic oracle associated with the objective and constraints.

\begin{assumption}   
   \textbf{(Lipschitz continuity)}
   \label{ass: Lipschitz continuity}
   Assume the $\nabla F(x, y)$, $\nabla G(x, y)$, $\nabla H(x, y)$ are $L_F, L_G, L_H$-Lipschitz continuous, respectively.
\end{assumption}

\begin{assumption}   
   \textbf{(Convexity)}
\label{ass: convexity}
   Assume $G(x, y)$ is $\mu_{G}$-strongly convex in $y$ for any $x \in X$, $H(x, y)$ is convex in $y$ for any $x \in X$.
\end{assumption}

\begin{remark}
Assumption~\ref{ass: convexity} implies that $y^*(x)$ defined in~\eqref{eq: BLO} is unique for any $x \in X$. 
\end{remark}
\begin{assumption}
    \textbf{(Boundedness)}
    \label{ass: boundedness}
    Assume $\nabla G(x, y)$, $H(x, y)$ and $\nabla H(x, y)$ are bounded, that is,
   \begin{equation*}
   \begin{aligned}
      &\|\nabla G(x, y)\| \leq M_{G, 1}, \\
        |H_i(x, y) | \leq &M_{H, 0},  
        ~\| \nabla H_i(x, y) \| \leq M_{H, 1}, 1\leq i\leq p. 
        \end{aligned} 
    \end{equation*}
\end{assumption}

\begin{assumption}   
   \textbf{(Stochastic derivative)}
   \label{ass: stochastic oracle}
   The stochastic oracle $\nabla f(x, y; \zeta)$, $\nabla g(x, y; \xi)$ are unbiased estimator of $\nabla F(x, y)$, $\nabla G(x, y)$, respectively, and their variances are bounded by $\sigma_{f}^2, \sigma_{g}^2$, respectively. 
\end{assumption}

To ensure the strong duality and regularity of the optimal points, we assume Slater's condition and linear independence constraint qualification (LICQ) hold for the lower-level constants, which are common in nonlinear optimization analysis.
\begin{assumption}
   \textbf{(LL Slater's condition)}
   \label{ass: LL Slater's condition}
   For any fixed $x \in X$, Slater's condition holds for~\eqref{eq: LL}, that is, there exist $\epsilon_0(x) > 0$ and $y_0(x)$ such that 
   \[
   H_i(x, y_0(x)) < -\epsilon_0(x), \quad i=1,...,p .
   \]
\end{assumption}

\begin{assumption}
   \textbf{(LICQ)}
   \label{ass: LICQ}
   For any $x \in X$ and $y = y^*(x)$, $\{\nabla_y H_i(x, y) | H_i(x, y) =0 \} $ is linearly independent. 
   Denote the matrix $\mathcal{C}(x, y) =[\nabla_y H_i(x, y) ]_{i \in \{i | H_i(x, y) =0 \}} $.
   Since $X, Y$ are compact sets, we further assume the smallest singular value satisfies
   $$
   \sigma_{\min} (\mathcal{C}(x, y)\mathcal{C}(x, y)^\top
   ) \geq \sigma_0^2 > 0, (x, y) \in X \times Y.
   $$
\end{assumption}

\subsection{Equivalence of reformulations}

In this subsection, the equivalence between the reformulations~\eqref{eq: maintext stochastic penalized single level} and the original BLO formulation in~\eqref{eq: BLO} is established. 
We take  $B = p^2 \sigma_0^2 M_{H,1}(M_{G,1} + p M_{H,1})$ as provided in~\eqref{eq: maintext stochastic penalized single level} throughout the analysis. The deterministic case is first analyzed as a special case. The equivalence is then extended to the stochastic case, emphasizing the key improvements introduced by the stochastic reformulation.

\subsubsection{Deterministic case}
The following theorem establishes the equivalence between the penalized form of~\eqref{eq: single level} and the original BLO.

\begin{theorem}
   \label{thm: maintext equivalence of single level}
   Suppose that Assumptions~\ref{ass: Lipschitz continuity} ~\ref{ass: convexity} and~\ref{ass: LICQ} holds and $\gamma_1, \gamma_2 > 0$ are fixed parameters. 

   \noindent
   1. Assume $(x^*, y^*)$ is a global solution to~\eqref{eq: BLO} and $c_1\geq   \frac{L}{ 2 \mu_G} \epsilon^{-1}, c_2 \geq  (c_1)^2 B^2\epsilon^{-1}$. 
   There exists $z^* \in \mathbb{R}_+^p$ such that $(x^*, y^*, z^*)$ is a $\epsilon$-global-minima of the following penalized form
   \begin{equation}
      \label{eq: maintext penalized single level 2}
   \begin{aligned}
      \min_{(x, y, z) \in X \times Y \times Z} &\quad  \Psi(x, y, z) = F(x, y)  + c_1 \mathcal{G}(x, y, z) + \frac{c_2}{2} \sum_{i=1}^p [H_i(x, y)]_+^2.
   \end{aligned}
   \end{equation}
   \noindent 
   2. By taking $c_1 = c_1^* + 2 :=  \frac{L}{ 2 \mu_G} \epsilon^{-1} + 2, c_2 = c_2^* + 2 := (c_1^*)^2 B^2\epsilon^{-1}$, 
   any $\epsilon$-global-minima of~\eqref{eq: maintext penalized single level 2} is an $\epsilon$-global-minima the following approximation of BLO
   \begin{equation}
      \label{eq: maintext approximation of BLO 2}
      \begin{aligned}
         \min_{ (x, y, z) \in X \times Y \times Z} &\quad  F(x, y) 
         \quad \text{s.t.} \quad  G(x, y) - E(x, z) \leq \epsilon_1, \quad  \frac{1}{2} \sum_{i=1}^p [H_i(x, y)]_+^2 \leq \epsilon_2^2,
      \end{aligned}
   \end{equation}
   with some $\epsilon_1 \leq \epsilon$, respectively.
\end{theorem}

This theorem indicates that the penalized reformulation can approximate the original bilevel optimization problem within a controlled error bound.

\subsubsection{Stochastic case}
The major difference between the deterministic and stochastic reformulation is the inexact solution of the subproblem~\eqref{eq: minimax subproblem}. By controlling the inexactness, we design an approximated stochastic reformulation for the stochastic bilevel optimization problem.
Denote the penalized form as
$$
\begin{aligned}
&\Psi(x, y, z, \hat w, \hat \lambda; \mathbf{\xi}) 
= F(x, y)  + c_1 (G(x, y)  
  -\ell_{\gamma}(x, z, \hat w(\mathbf{\xi}), \hat \lambda(\mathbf{\xi})) ) + \frac{c_2}{2} \sum_{i=1}^{p} [H_i(x, y)]_+^2.
\end{aligned}
$$
The following theorem shows the equivalence between this penalized form and~\eqref{eq: BLO} (see Theorem~\ref{thm: equivalence of stochastic single level} for proofs). 

\begin{theorem}
   \label{thm: maintext equivalence between BLO and stochastic single level approximation}
   Suppose that Assumptions~\ref{ass: Lipschitz continuity},~\ref{ass: convexity} and~\ref{ass: LL Slater's condition} holds and $\gamma_1, \gamma_2 > 0$ are fixed parameters. 

\noindent
1. 
Assume $(x^*, y^*)$ is a global solution to~\eqref{eq: BLO}. 
If $\mathcal{P}(\delta)$ defined in~\eqref{eq: P delta} is nonempty for any $(x, z) \in X \times \mathbb{R}_+^p$,
then for any $(\hat w, \hat \lambda) \in \mathcal{P}(\delta)$,
there exists $z^*$ such that 
$(x^*, y^*, z^*)$ is a $\epsilon$-global-minima 
of the following penalized form 
there exists $z^*$ such that 
$(x^*, y^*, z^*)$ is a $\epsilon$-global-minima 
of the following penalized form 
\begin{equation}
   \label{eq: maintext stochastic penalized single level 2}
\begin{aligned}
\min_{
(x, y, z) \in X \times Y \times Z } &\quad  \mathbb{E} [\Psi(x, y, z, \hat w, \hat \lambda; \mathbf{\xi})] 
\end{aligned}
\end{equation}
with any
$c_1 \geq  \frac{2 L}{3 \mu_G} \epsilon^{-1}$, $c_2 \geq \frac{3}{2} (c_1)^2 B^2\epsilon^{-1}$ and $\delta \leq \frac{\epsilon}{6 c_1}$.

\noindent 
2.  
By taking $c_1 = c_1^* + 2 :=  \frac{2 L}{3 \mu_G} \epsilon^{-1} + 2, 
c_2 = c_2^* + 2 := \frac{3}{2} (c_1^*)^2 B^2\epsilon^{-1} + 2$ and
$\delta \leq \frac{\epsilon}{6c_1}$,
for any $(\hat w, \hat \lambda) \in \mathcal{P}(\delta)$,
the $\epsilon$-global-minima of~\eqref{eq: maintext stochastic penalized single level 2} is a $\epsilon$-global-minima of the following approximation of BLO:
\begin{equation}
   \label{eq: maintext stochastic approximation of BLO 2}
   \begin{aligned}
      \min_{
         (x, y, z) \in X \times Y \times Z 
         } &\quad  F(x, y)  \\
      \text{s.t.} \quad \quad & \quad  G(x, y) - \mathbb{E} [\ell_{\gamma}(x, z, \hat w(\mathbf{\xi}), \hat \lambda(\mathbf{\xi}))] \leq \epsilon_1,\\  &\quad \frac{1}{2} \sum_{i=1}^p [H_i(x, y)]_+^2 \leq \epsilon_2^2,
   \end{aligned}
\end{equation}
with some $\epsilon_1, \epsilon_2 \leq\frac{13}{12}\epsilon$.
\end{theorem}

\subsection{Convergence analysis}
Denote $\mathcal{F}_k$ and $\tilde{\mathcal{F}}_k$ as the $\sigma$-algebra generated by $\{\mathbf{ \xi_l}\}_{l=0}^{k} \cup \{\mathbf{\zeta}^l, \mathbf{\tilde \xi}^l \}_{l=0}^{k-1}$ and $\{\mathbf{\xi}^l\}_{l=0}^{k} \cup \{\mathbf{\zeta}^l, \mathbf{\tilde \xi}^l\}_{l=0}^{k}$ respectively. Then 
$$
\emptyset = \tilde{\mathcal{F}}_0 \subset \mathcal{F}_1 \subset \tilde{\mathcal{F}}_1 \subset \cdots \subset \mathcal{F}_k \subset \tilde{\mathcal{F}}_k \subset \cdots
$$ 
is the filtration generated by the random variables in Algorithms~\ref{alg: main} and~\ref{alg: variance reduce}.

\begin{theorem}

   \label{thm: convergence of ALM}
   Suppose 
    Assumptions~\ref{ass: Lipschitz continuity}-\ref{ass: LL Slater's condition}
    holds.
      By taking the step sizes in Algorithm~\ref{alg: ALM} as 
      \begin{equation}
      \eta_j = \frac{\eta}{j+1}, \quad \rho_j = \frac{\rho}{j+1},
      \end{equation}
      there exist constants $\bar{\phi}_1, \bar{\phi}_2 > 0$ such that the output pair $(w^{k, s}, \lambda^{k ,s})$ satisfying
      \begin{equation}
         \label{eq: maintext point convergence of inner loop}
         \begin{aligned}
          \mathbb{E}_{\mathbf{\xi}} \| w^{k, s} - w^*(x^{k-1}, z^{k-1})\|^2 &\leq \bar{\phi}_1 \frac{1+\log(s)}{s} ,\\
          \mathbb{E}_{\mathbf{\xi}} \| \lambda^{k, s} - \lambda^*(x^{k-1}, z^{k-1})\|^2 &\leq \bar{\phi}_2 \frac{1+\log(s)}{s} .
         \end{aligned}
      \end{equation}
\end{theorem}

Define the bias of the gradient oracle of $\Psi$ as
   \begin{equation}
       \label{eq: bk}
   b^k = \nabla \Psi^k(\mathbf{u}^k; \mathbf{\zeta}^k, \mathbf{\tilde{\xi}}^k) - \nabla \Psi(\mathbf{u}^k).
   \end{equation}
   We establish the following lemma to control the bias of the gradient oracle  in terms of conditional expectation.
   \begin{lemma}
      \label{lem: maintext bound of bias}
       The bias $b^k$ has a controllable bound as 
       \begin{equation}
          \label{eq: maintext bound of bias}
          \begin{aligned}
           \mathbb{E} [ \| b^k \|^2 | \tilde{\mathcal{F}}_{k-1}] 
           &\leq   
            2 \left(\frac{\sigma_{f}^2}{r_k} + c_1^2 \frac{\sigma_{\mathcal{G}}^2}{q_k} + c_1^2(\epsilon_{\mathcal{G}}^k)^2  \right).
      \end{aligned}
       \end{equation}
   \end{lemma}
 Here $\mathbb{E}[\cdot] $ is the abbreviation of $\mathbb{E}_{\mathbf{\zeta}^k, \mathbf{\tilde{\xi}}^k, \mathbf{\xi}^k}[\cdot]$ and $\epsilon^k_{\mathcal{G}}, (\sigma_{\mathcal{G}}^k)^2$ denote the upper bounds of the bias and variance of $\nabla \mathcal{G}^k(\mathbf{u}; \mathbf{\zeta}^k, \mathbf{\tilde{\xi}}^k)$, respectively. 
 $\epsilon^k_{\mathcal{G}}, (\sigma_{\mathcal{G}}^k)^2$ are constants conditioned on $\tilde{\mathcal{F}_{k-1}}$ and can be further bounded by polynomials of $\mathbb{E}_{\mathbf{\xi}^k} \| w^{k, s} - w^*(x^{k-1}, z^{k-1})\|^2$, $\mathbb{E}_{\mathbf{\xi}^k} \| \lambda^{k, s} - \lambda^*(x^{k-1}, z^{k-1})\|^2$. Therefore we can control the bias $b^k$ by enhancing the accuracy of Algorithm~\ref{alg: ALM}.
 
Denote $c = \max(c_1, c_2)$.
With the convergence results of Algorithm~\ref{alg: ALM} and the boundedness of gradient oracle, we maintain  the convergence of Algorithm~\ref{alg: main}.

\begin{theorem}
   \label{thm: maintext convergence of SALVF}
   Suppose Assumptions~\ref{ass: Lipschitz continuity}-\ref{ass: 
   LICQ}
   Take constant step size $\alpha_k = {\alpha} < \frac{1}{2L_{\Psi}}$ constant sample sizes as $r_k = r, q_k = q$, and $s_k = s$.
   Then the sequence $\{\mathbf{u}^k\}_{k=0}^{K}$ generated by Algorithm~\ref{alg: main} satisfies
   \begin{equation*}
\begin{aligned}
\mathbb{E}\left[ \frac{1}{K}\sum_{k=0}^{K-1} \frac{1}{\alpha_k}\|\mathbf{u}^{k+1} - \mathbf{u}^k \|^2 \right]
    &\leq \mathcal{O}\left(\frac{1}{\alpha K} +  \frac{1}{r} + \frac{c_1^2}{q} + \frac{c_1^2}{s} \right).
\end{aligned}
\end{equation*}
\end{theorem}

The detailed proof is available in Theorem~\ref{thm: convergence of outer loop} and Corollary~\ref{cor: rate of outer loop}. We summarize the proof idea of Theorem~\ref{thm: maintext convergence of SALVF} as follows. By combining Theorem~\ref{thm: convergence of ALM} with Lemma~\ref{lem: maintext bound of bias}, the bias and variance are bounded with $\widetilde{O}(s_k)$. Further analyzing the biased stochastic gradient descent methods provides the desired result.
    
\begin{corollary}
   \label{rem: maintext convergence of SALVF}
    If the rate is measured by $\mathrm{dist}(0, \nabla \Psi(\mathbf{u}) + \mathcal{N}_{\mathcal{U}}(\mathbf{u}))^2 $, we have the following equivalent conclusion (see Theorem~\ref{thm: gradient bound of Psi} and corollary~\ref{cor: rate of outer loop} for detailed proof):
    \begin{equation*}
        \begin{aligned}
 \mathbb{E} [\mathrm{dist}(0, \nabla \Psi(\mathbf{u}^R) + \mathcal{N}_{\mathcal{U}}(\mathbf{u}^{R}))^2]  
    &\leq \mathcal{O}\left(\frac{1}{\alpha K} + 
    \frac{1}{r} + \frac{c_1^2}{q} + \frac{c_1^2}{s}\right).
\end{aligned}
\end{equation*}
\end{corollary}

\begin{remark}
   The step size condition $\alpha_k < \frac{1}{2L_{\Psi}}$ and~\eqref{eq: Lipschitz continuity of Psi} implies $\alpha_k$ is a most $\widetilde{\mathcal{O}}(c^{-1})$. 
   With $\alpha \sim \mathcal{O}(c^{-1})$, 
   $r \sim \mathcal{O}(\epsilon^{-1})$,
   $q \sim \mathcal{O}(c_1^2 \epsilon^{-1})$,  $s \sim {\mathcal{O}}(c_1^2 \epsilon^{-1})$, $K \sim {\mathcal{O}}(c \epsilon^{-1})$, the right side of the above inequality is $\widetilde{\mathcal{O}}(\epsilon)$.
  Then the sample complexity on $(\zeta, \xi)$ is $(\widetilde{\mathcal{O}}(c  \epsilon^{-2}), \widetilde{\mathcal{O}}(c c_1^2 \epsilon^{-2}))$.
\end{remark}

\begin{remark}
   Theorem~\ref{thm: maintext equivalence between BLO and stochastic single level approximation} shows that~\eqref{eq: maintext stochastic penalized single level} is equivalent to the original problem~\eqref{eq: BLO} in the sense of $\epsilon$-accuracy by taking $c_1 \sim \mathcal{O}(\epsilon^{-1})$, $c_2 \sim \mathcal{O}(\epsilon^{-3})$ and $\delta \sim \mathcal{O}(\epsilon^{-2})$. Under this condition, the sample complexity on $(\zeta, \xi)$ is $(\widetilde{\mathcal{O}}( \epsilon^{-5}), \widetilde{\mathcal{O}}(\epsilon^{-7}))$.
\end{remark}

By introducing the following averaged Lipschitz assumption, we can further improve the convergence rate utilizing variance reduction techniques.

\begin{assumption}
   \label{ass: averaged Lipschitz}
   Assume $\nabla f(x, y; \zeta), \nabla g(x, y; \xi)$ are averaged Lipschitz continuous, that is, 
   \begin{equation}
      \begin{aligned}
         & \mathbb{E}_{\zeta} [\| \nabla f(x_1, y_1; \zeta) - \nabla f(x_2, y_2; \zeta) \|^2] 
         \leq L_f^2 \| (x_1, y_1) - (x_2, y_2) \|^2, \\
         & \mathbb{E}_{\xi} [\| \nabla g(x_1, y_1; \xi) - \nabla g(x_2, y_2; \xi) \|^2] 
          \leq L_g^2 \| (x_1, y_1) - (x_2, y_2) \|^2.
      \end{aligned}
   \end{equation}
\end{assumption}

\begin{theorem}
   \label{thm: maintext convergence of SALVF-VR}
   Suppose Assumptions~\ref{ass: convexity}-\ref{ass: LICQ} and~\ref{ass: averaged Lipschitz} hold.
   Take $\alpha_k = \alpha (k+1)^{-\frac{1}{3}}$ and $\beta_{k+1} = \beta \alpha_k^2$ in the outer loop and take constant sample sizes as $r_k = r, q_k = q$, and $s_k = s$.
   Then the sequence $\{\mathbf{u}^k\}_{k=0}^{K}$ generated by Algorithm~\ref{alg: variance reduce} satisfies
\begin{equation}
      \begin{aligned}
         &\quad \frac{1}{\sum_{k=0}^K \alpha_k}\sum_{k=0}^{K-1} \mathbb{E}\left[ \frac{1}{\alpha_k}\|\mathbf{u}^{k+1} - \mathbf{u}^k\|^2 \right] 
         \leq\widetilde{\mathcal{O}} \left(\frac{1}{\alpha K^{\frac{2}{3}}} + (c_1^2 + \frac{K^{\frac{2}{3}}}{r})(\frac{1}{q} + \frac{1}{s})
           + \frac{1}{K^{\frac{2}{3}}r} \right).
       \end{aligned}
\end{equation}

\end{theorem}
The detailed proof is available in Theorem \ref{thm: decrease of variance reduction} and Corollary~\ref{cor: rate of variance reduction}. The proof sketch is summarized as follows: first we show the error of the direction $e^k = d^k - \mathbb{E}[\Psi^k(\mathbf{u}^k; \mathbf{\zeta}^k, \mathbf{\tilde{\xi}}^k)]$ is bounded by the linear functions of the previous error $e^{k-1}$. By proving the reduction of the merit function
$
\Psi(\mathbf{u}^k) + \theta_k \|e^k\|
$ with some coefficients $\theta_k$, we establish the convergence of Algorithm~\ref{alg: variance reduce}.

\begin{remark}
   \label{rem: maintext convergence of SALVF-VR}
Further take  $K \sim \mathcal{O}(c^{1.5}\epsilon^{-1.5}  )$, $r \sim \mathcal{O}(1)$, $q \sim \mathcal{O}(c_1^2 \epsilon^{-1})$, $s \sim \mathcal{O}(c_1^2 \epsilon^{-1})$, the right side is
$\widetilde{\mathcal{O}}(\epsilon)$.  The sample complexity on $(\zeta, \xi)$ is $(\widetilde{\mathcal{O}}(c^{1.5} \epsilon^{-1.5}), \widetilde{\mathcal{O}}(c^{1.5} c_1^2 \epsilon^{-2.5}))$. From the Lemma~\ref{thm: convergence of ALM},~\ref{lem: maintext bound of bias} we know the upper bound of the $\|b^k\|^2$ is $\widetilde{\mathcal{O}}(\frac{c_1^2}{s})$, hence to achieve an $\epsilon$-optimal solution by biased gradient-based approach, the condition $\widetilde{\mathcal{O}}(\frac{c_1^2}{s}) \leq \epsilon$ cannot be be further improved. That is, the sample complexity on $\xi$ in each iteration is at least $\mathcal{O}(c_1^2\epsilon^{-1})$.
To this extent, current analysis requires a larger sample complexity on the lower-level to reduce the upper-level complexity, and it remains an open question whether variance-reduction could reduce the sample complexity fo upper- and lower-level at the same time. 

\end{remark}



\section{Numerical experiments}
This section presents numerical experiments to demonstrate the effectiveness of the proposed algorithms, compared with baselines including LV-HBA~\cite{yao2024constrained}, GAM~\cite{xu2023efficient}, and BLOCC~\cite{jiang2024primal}. 
\subsection{Toy example}
\begin{figure}[ht]
   \centering
   \begin{subfigure}[b]{0.48\textwidth}  
       \centering
       \includegraphics[width=\textwidth]{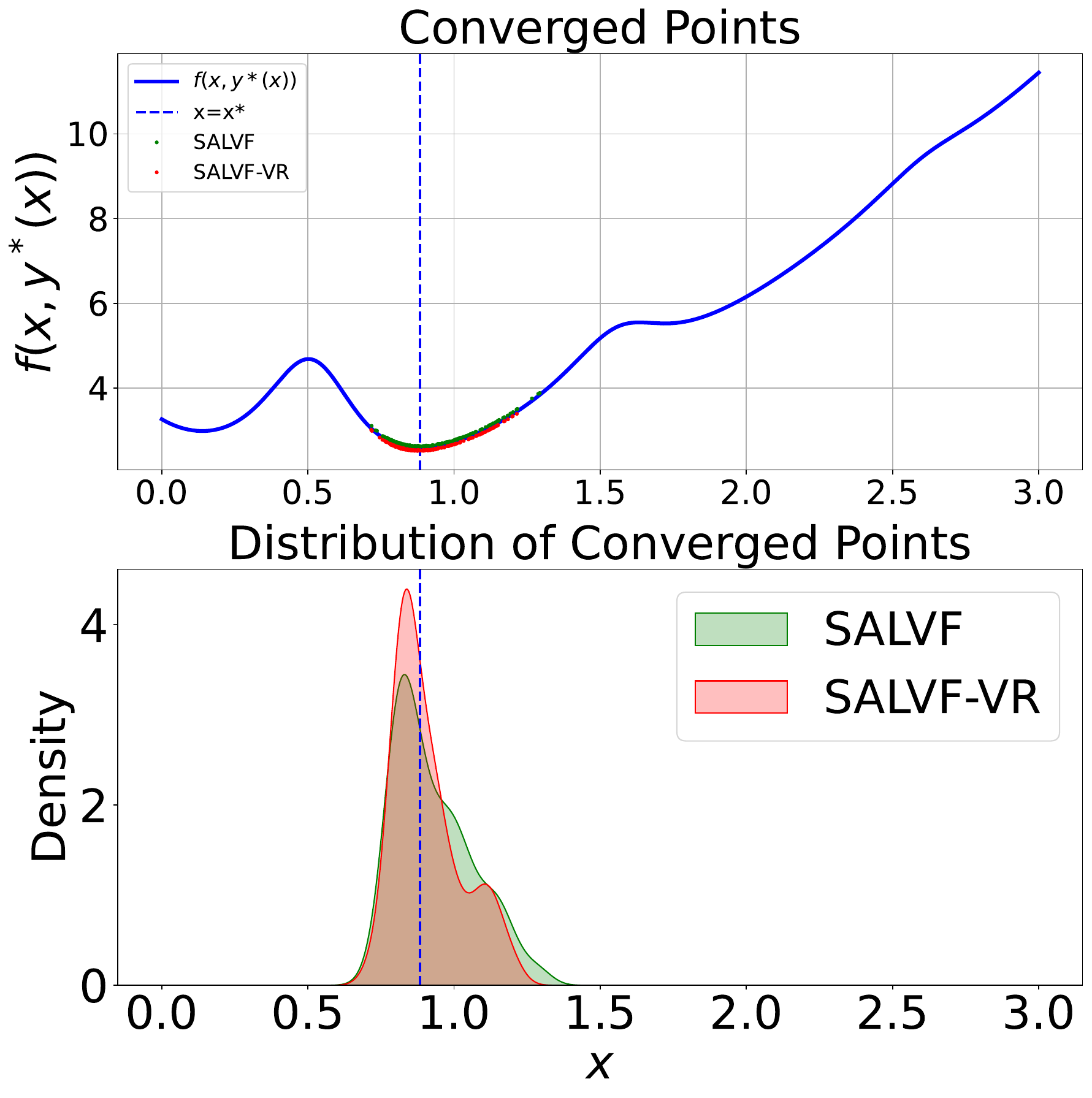}
       \caption{$x$ and its distribution}
       \label{fig: toyexample density}
   \end{subfigure}
   \hfill
   \begin{subfigure}[b]{0.48\textwidth}  
       \centering
       \includegraphics[width=\textwidth]{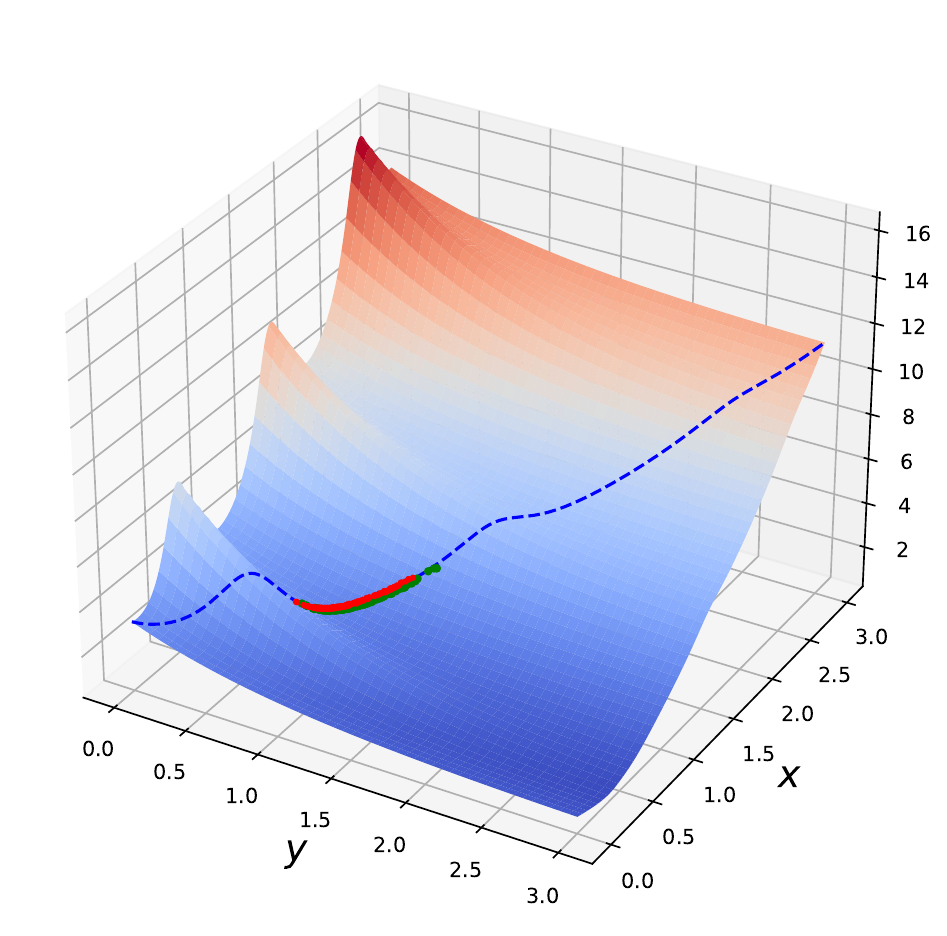}
       \caption{3D plot of convergent points}
       \label{fig: toyexample 3d}
   \end{subfigure}
   \caption{The converged points of Algorithm~\ref{alg: main}.}
   \label{fig: toy_example}
\end{figure}

Consider the following example from~\citet{jiang2024primal}:
$$
\begin{aligned}
& \min _{x \in[0,3]} F\left(x, y^*\right)=\frac{e^{-y^*(x)+2}}{2+\cos (6 x)}+\frac{1}{2} \log \left((4 x-2)^2+1\right) \\
& \text { s.t. } y^*(x) \in \arg \min _{y \in \mathcal{Y}(x)} G(x, y)=(y-2 x)^2,
\end{aligned}
$$
where $\mathcal{Y}(x)=\{y \in [0, 3]| H(x, y)\leq 0\}$ and $H(x, y)=y - x$. The lower-level problem has closed form solution $y^*(x) = x$.
Now assume the gradient oracle of $F$ and $G$ have Gaussian noise with variance $\sigma=0.1$, that is,
$$
   \nabla f(x, y; \zeta) = \nabla F(x, y) + \zeta, 
   \nabla g(x, y; \zeta) = \nabla G(x, y) + \xi,
$$
with $\zeta \sim N(0, \sigma^2I), \xi \sim N(0, \sigma^2I)$.
We pick 200 random points $(x^0, y^0) \in [0, 3] \times [0, 3]$ as initial points and allow the maximum sample on $\zeta$ as $2500$. The final iterated points of
Algorithm~\ref{alg: main} and~\ref{alg: variance reduce}
are collected in Figure~\ref{fig: toy_example}. Figure~\ref{fig: toyexample density} plots the points $(x, y)$ projected on the $y = y^*(x)$ and the distribution of the output $x$ and Figure~\ref{fig: toyexample 3d} shows a 3D plot of the output $x$ and $y$. 
As shown in the figures, the converged points of both algorithms
are close to the global optimal solution $x^*$ and form an approximate Gaussian distribution. 
The distribution of SALVF-VR is more concentrated than SALVF, which demonstrates the acceleration effect of variance reduction techniques.

\subsection{Hyperparameter tuning for SVM}
\begin{figure}[htbp]
   \centering
   \begin{subfigure}[b]{0.48\textwidth}
       \centering
       \includegraphics[width=\textwidth]{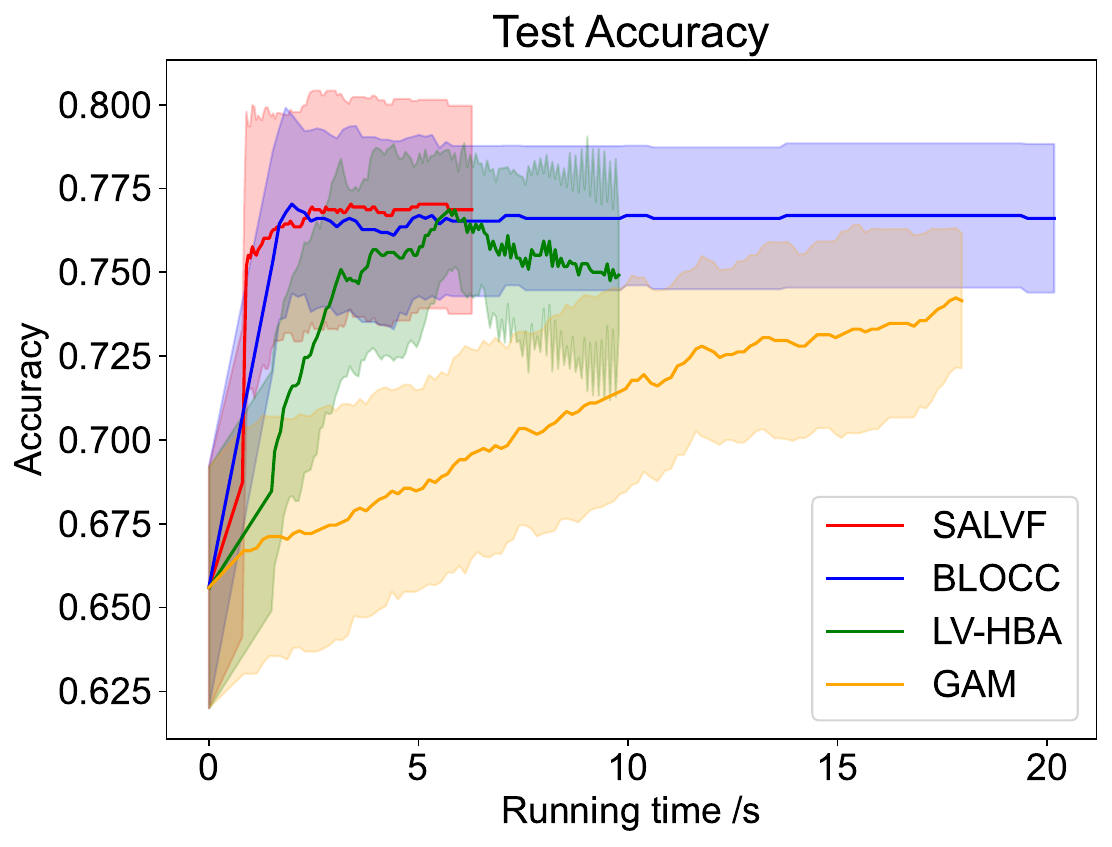}
       \caption{Diabetes: test acc. v.s. time}
       \label{fig: SVM_test_acc_vs_time}
   \end{subfigure}
   \hfill
   \begin{subfigure}[b]{0.48\textwidth}
       \centering
       \includegraphics[width=\textwidth]{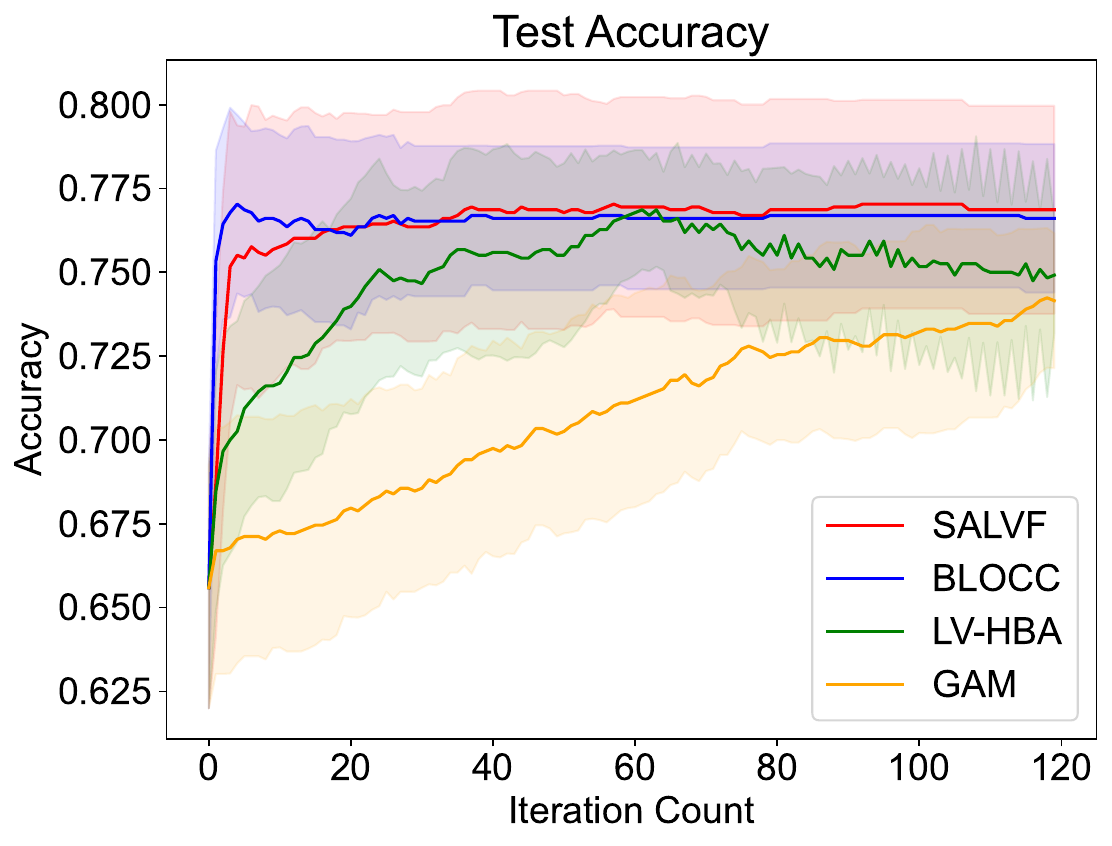}
       \caption{Diabetes: test acc. v.s. iter.}
       \label{fig: SVM_test_acc_vs_iter}
   \end{subfigure}
   \hfill
   \begin{subfigure}[b]{0.48\textwidth}
       \centering
       \includegraphics[width=\textwidth]{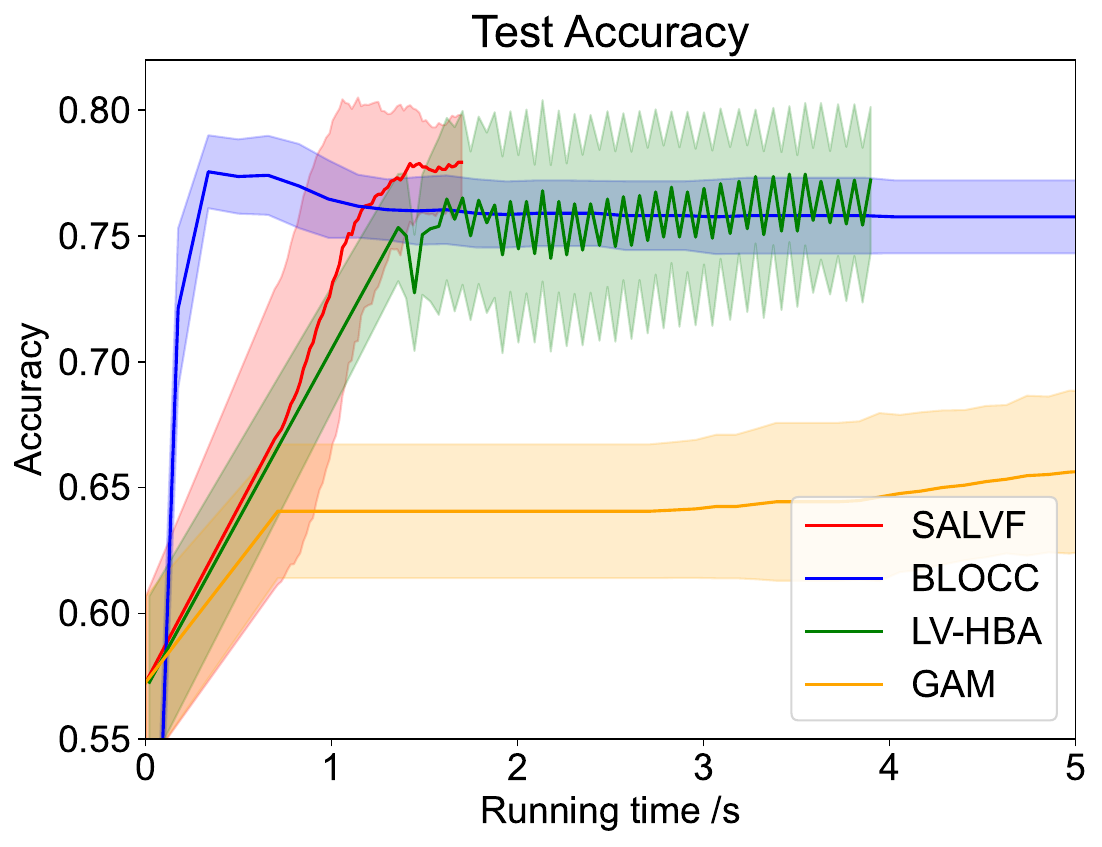}
       \caption{Fourclass: test acc. v.s. time}
       \label{fig: SVM_test_acc_vs_time_fourclass}
   \end{subfigure}
   \hfill
   \begin{subfigure}[b]{0.48\textwidth}
       \centering
       \includegraphics[width=\textwidth]{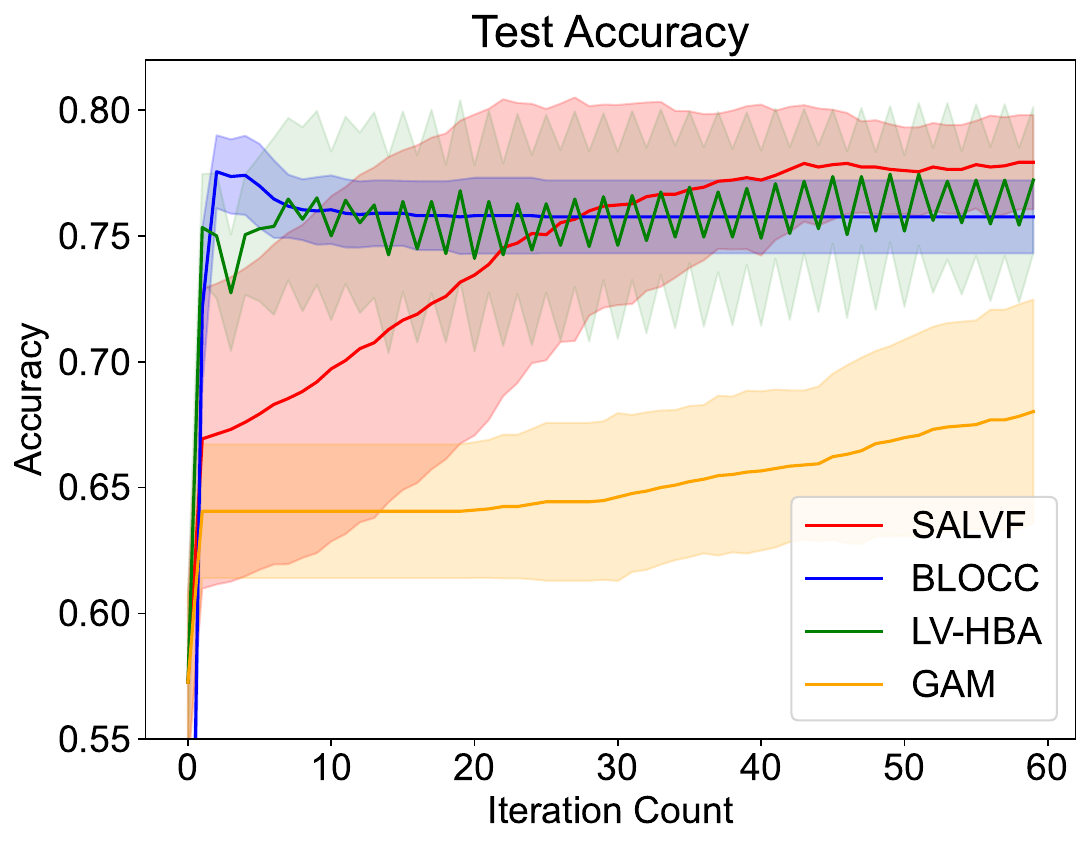}
       \caption{Fourclass: test acc. v.s. iter.}
       \label{fig: SVM_test_acc_vs_iter_fourclass}
   \end{subfigure}
   \caption{The performance of SALVF compared with baselines on SVM hyperparameter optimization. The abbreviations ``test acc.'' and ``iter.'' stand for test accuracy and iterations, respectively. The curves are averaged over 10 random seeds. The curves in Figure~\ref{fig: SVM_test_acc_vs_time},~\ref{fig: SVM_test_acc_vs_time_fourclass} are clipped at the maximum iteration 120 and 60, respectively.}
   \label{fig: SVM}
\end{figure}

Consider the following support vector machine (SVM) problem as the lower-level problem:
$$
\begin{aligned}
\min_{w ,b ,\xi} &~ 
\frac{1}{2}\|w\|^2+\frac{1}{2} \frac{1}{\sum_{i=1} \exp(c_i)}\sum_{i=1}^{N} \exp({c_i})  \xi_i  \\
\text{s.t. } &~ l_i(z_i^\top w +b) \geq 1 - \xi, \quad \forall (z_i, l_i) \in {D}_{tr},
\end{aligned}
$$
where $c_1 = (c_1, ..., c_N)$ is the hyperparameters, $\mathcal{D}_{tr}$ is the training set, $N = |\mathcal{D}_{tr}|$ and $l_i$ is the label of the $i$-th sample. The upper-level problem is to minimize the validation error with respect to $c_1$:
$$
\min_c 
{\mathbb{E}}_{{(z, l) \sim \mathcal{D}_{val}}} [\exp\left( 
1-l (z^\top w^*+b^*) 
\right)] .
$$
In Figure~\ref{fig: SVM},
we compare the performance of SALVF with baselines, LV-HBA~\cite{yao2024constrained}, GAM~\cite{xu2023efficient} and BLOCC~\cite{jiang2024primal}. 
Since the lower-level problem is deterministic, we allow all algorithms to access the exact optimal solution of the lower-level problem by calling the ECOS~\cite{domahidi2013ecos} solver and we set $\gamma_2 = 0$. We extend BLOCC, LV-HBA and GAM to their stochastic versions by replacing (projected) gradient descent with (projected) stochastic gradient descent in the upper-level update.
Figure~\ref{fig: SVM} shows the test accuracy of SALVF versus time and iterations on the Diabetes and Fourclass datasets. SALVF achieves a better accuracy over the baselines. Although BLOCC also has an approximate peak accuracy, we can see that iterations of SALVF are more time-efficient. This is because SALVF requires a double loop
 while BLOCC requires a triple loop, which is more computationally expensive.

\subsection{Weight decay tuning}
\begin{figure}[htbp]
   \centering
   \begin{subfigure}[b]{0.48\textwidth}
       \centering
       \includegraphics[width=\textwidth]{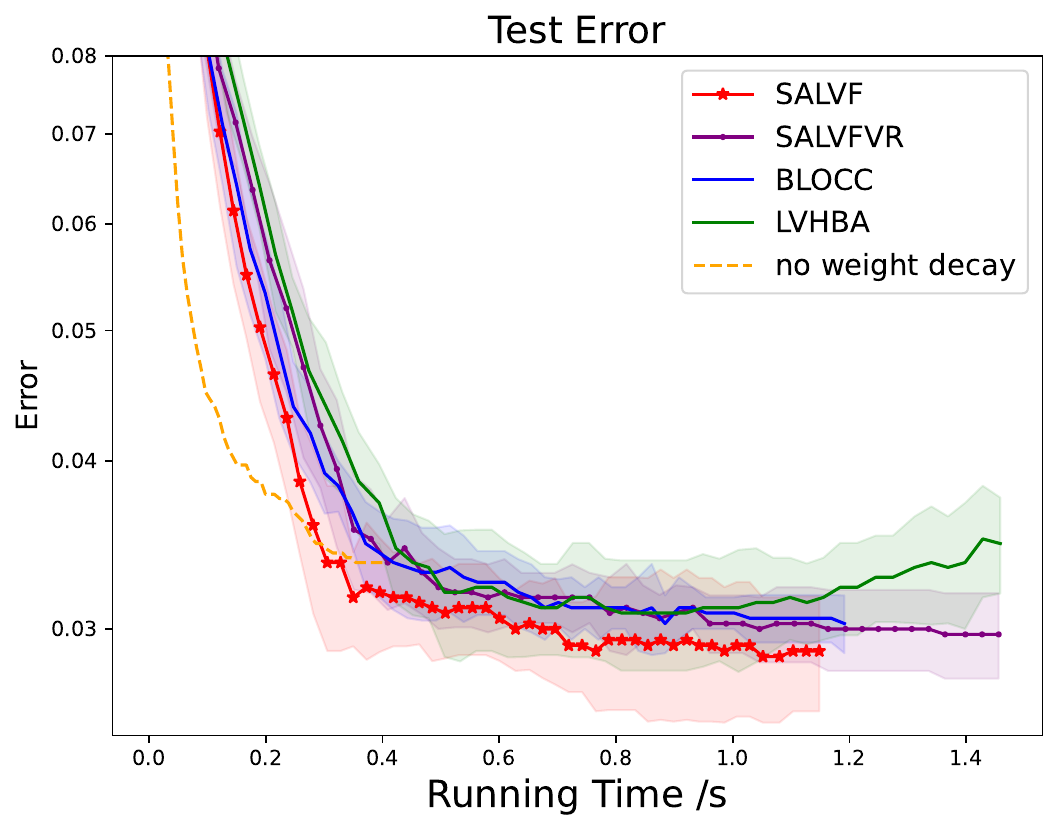}
       \caption{Test acc. v.s. time}
       \label{fig: weight_decay_test_error_vs_time}
   \end{subfigure}
   \hfill
   \begin{subfigure}[b]{0.48\textwidth}
       \centering
       \includegraphics[width=\textwidth]{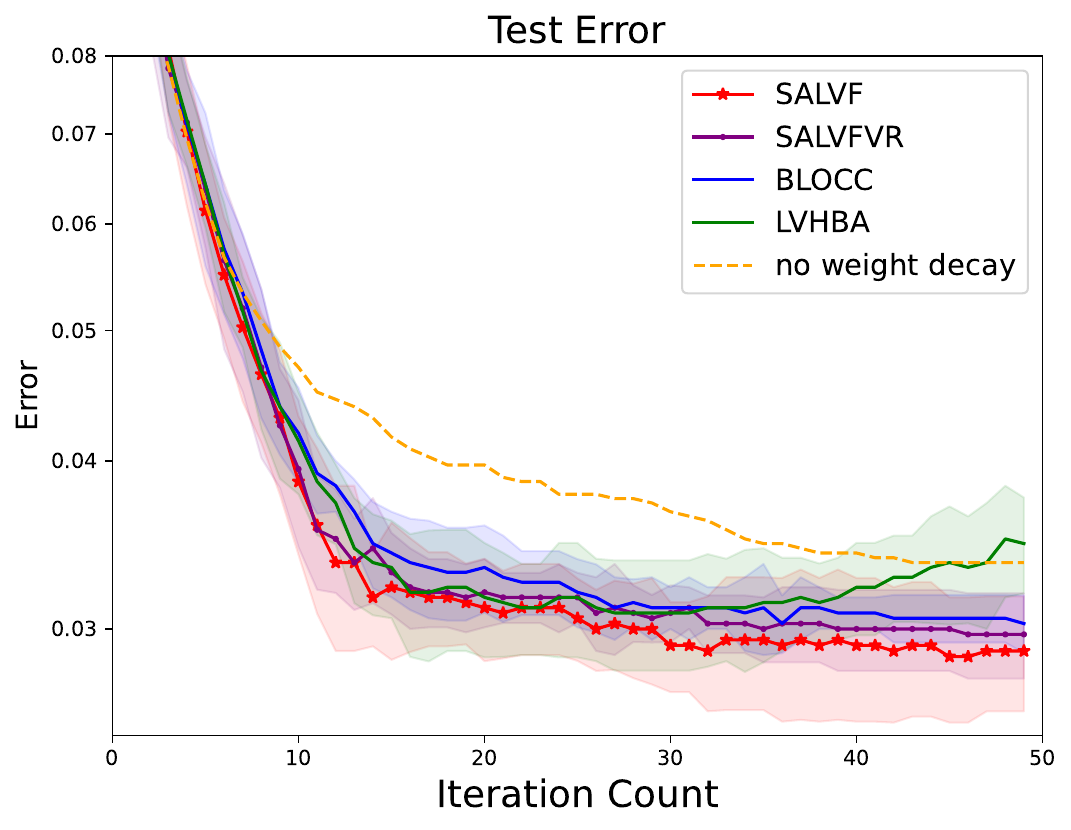}
       \caption{Test acc. v.s. iter.}
       \label{fig: weight_decay_test_error_vs_iter}
   \end{subfigure}
   \caption{The performance of SALVF compared with baselines on digit dataset. The curves are averaged over 10 random seeds. The curves in Figure~\ref{fig: weight_decay_test_error_vs_time} are clipped at the maximum iteration 50.}
   \label{fig: weight_decay}
\end{figure}

Given a neural network $f(w, x)$ where $w$ are the weight and bias parameters in each layer, the goal is to optimize the weight decay parameter $\lambda$ for a neural network model. 
To improve the generalization performance, the weight decay parameters $C$ are introduced, which impose the constraint $\|w\| \leq C$. This can be formulated as the following stochastic BLO:
$$
\begin{aligned}
\min _{C>0} \quad &\mathbb{E}_{(x, y) \sim \mathcal{D}_{val}}[\ell(y, f(w^*(C), x))] \\
\text{s.t.} \quad &w^*(C)\in\arg \min _{\|w\| \leq C} \mathbb{E}_{(x, y) \sim \mathcal{D}_{tr}} \ell\left(y_i, f\left(w, x_i\right)\right).
\end{aligned}
$$
The upper level focuses on performance on a validation set, while the lower level involves constrained classifier training. 
We compare the performance of SALVF, SALVF-VR and BLOCC on the digit dataset~\cite{alpaydin1996pen} with a two-layer MLP as the base model. The results are shown in Figure~\ref{fig: weight_decay}. The ``no weight decay'' curve represents the model's performance without weight decay. By incorporating weight decay, all bilevel methods exhibit improved performance and reduce overfitting. From Figure~\ref{fig: weight_decay_test_error_vs_time} we see that SALVF is the most time-efficient, thanks to the simplicity of each step in its double-loop iteration process.

\section{Conclusion}

In this paper, we propose a Hessian-free  method for solving stochastic LC-BLO problems.
We present the first non-asymptotic rate analysis on value function-based algorithms for nonlinear LC-BLO in the stochastic setting.
The sample complexity of our algorithm is $
(\widetilde{\mathcal{O}}(c_1\epsilon^{-2}), \widetilde{\mathcal{O}}(c_1^3\epsilon^{-2}))
$ on upper- and lower-level random variables, respectively. The sample complexity of upper-level variables is further improved to $\widetilde{\mathcal{O}}(c_1^{1.5}\epsilon^{-1.5})$ using variance reduction techniques.
Numerical experiments on synthetic and real-world data demonstrate the effectiveness of the proposed approach.



%
%
%
\appendix
    \section{Convergence analysis}

    \subsection{Properties of Moreau envelope}
    \label{sec: Moreau envelope}
    Given function $\phi(z)$, the Moreau envelope is defined as
    $
    e^{\gamma} (z) = \min_{\lambda} \phi(\lambda) + \frac{\gamma}{2} \| \lambda - z\|^2.
    $To make a comprehensive explanation on the enveloped value function in \eqref{eq: augmented dual function}, we introduce some well-known properties of the Moreau envelope. The detailed proof can be referred to literature such as 
    \cite{rockafellar2009variational, grimmer2023landscape}.
    \begin{proposition}
       \label{prop: properties of Moreau envelope}
    
       \noindent
       1.
          The Moreau envelope $e^{\gamma}(z)$ is a continuous lower approximation of $\phi(z)$, i.e.,
          $
             \lim_{{\gamma} \to 0} e^\gamma(z) \leq  \phi(z), \forall z
          $
          and
          $
             \lim_{{\gamma} \to 0} e^\gamma (z) =  \phi(z)
          $.
          If $\phi(z)$ is $l$-Lipschitz continuous, $\rho$-weakly convex in $z$, then 
           this difference is at most $O({\gamma})$, i.e.,
          \[
          |e^\gamma(z) - \phi(z)| \leq O({\gamma}) .
          \]
    
       \noindent 2.
          If $\phi$ is $\rho$-strongly convex in $z$, $\rho \geq 0$, then $e^\gamma(z)$ is $(\frac{1}{\gamma} + \frac{1}{\rho})^{-1}$-strongly convex.
    
       \noindent 3.
          The Moreau envelope $e^{\gamma}(z)$ has the same minimizer as $\phi(z)$, i.e.,
          \[
             e^\gamma(z) = \phi(z), \forall z \in \arg \min_{z} \phi(z) = \arg \min_{z} e^\gamma(z) .
          \]
    
       \noindent 4.
          Its gradient at $z$ is given by
          $$
             \nabla e^\gamma(z) = \nabla \phi(\mathrm{prox}_{\gamma \phi}(z)) =
                {{\gamma}}(z - \mathrm{prox}_{{\gamma} \phi}(z)), 
          $$ 
          where $\mathrm{prox}_{\frac{1}{\gamma} \phi}(z) = \arg \min_{\lambda} \phi(\lambda) + \frac{{\gamma}}{2} \| \lambda - z\|^2$ is the proximal operator.
          Therefore $e^\gamma(z)$ is Lipschitz smooth given that $\phi(z)$ is Lipschitz smooth.
    \end{proposition}
    In our reformulation, we introduce $E(x, z)$ in~\eqref{eq: augmented dual function}. Then $- E(x, z)$ is the Moreau envelope of the augmented dual function $- D(x, z)$. Utilizing the properties of the Moreau envelope, we transform the original strongly-convex–concave saddle problem $\max_{z \in \mathbb{R}_+^p} \min _{y \in Y} \mathcal{L}_{\gamma_1}(x, y, z)$ into a strongly-convex–strongly-concave saddle-point problem 
    $\max_{w \in \mathbb{R}_+^p} \min _{\lambda \in Y} \tilde{\mathcal{L}}_{\gamma}(x, z, w, \lambda)$
    without changing the optimal solution.

    \subsection{Properties of the augmented Lagrangian function} \label{sec: properties of AL}
    
    In this section we review some important properties of the augmented Lagrangian function defined in~\eqref{eq: AL}. 
    \subsection{Augmented Lagrangian duality}
    In this ssubsection we introduce the augmented Lagrangian duality for general constrained optimization problems. The augmented Lagrangian duality is a powerful tool for solving constrained optimization problems, especially when the constraints are non-convex. The notation introduced in this subsection is only used within this subsection.
    
    For a general constrained optimization problems
    \begin{equation}
       \label{eq: general constrained optimization}
       \min_{y \in Y} \quad G(y) \quad \text{s.t. } \quad H(y) \leq 0,
    \end{equation}
    The augmented Lagrangian penalty term is defined as
    \begin{equation}
       \label{eq: augmented Lagrangian}
       \mathcal{A}_{\gamma}(y, z) = G(y) + \frac{1}{\gamma }\sum_{i=1}^{p} \left( [\gamma z_i + H_i(y)]_+^2 - \gamma^2 z_i^2 \right) .
    \end{equation}
    and the augmented Lagrangian dual function is defined as
    \begin{equation}
       \mathcal{L}_{\gamma}(y, z) = G(y) + \mathcal{A}_{\gamma}(y, z),
    \end{equation}
    where $z_i$ is the dual variable associated with the $i$-th constraint $H_i(y) \leq 0$ and $\gamma$ is the penalty parameter. The augmented Lagrangian dual function contains the Lagrangian function as a special case when $\gamma = + \infty$.
    
    Under the convexity assumption and Slater's condition, we have the following proposition about the augmented Lagrangian duality.
    \begin{proposition}
    \label{prop: strong duality}
    \textbf{(strong duality)}
          Suppose $G, H$ are convex and Slater's condition holds, i.e., there exists $y^* \in Y$ such that $H(y^*) < 0$. Then the following statements hold:
    
    \noindent
    1.
          The 
          dual variables $z$ exists and the strong duality of augmented Lagrangian holds, i.e.,
          \begin{equation*}
             \min_{y} \max_{z \in \mathbb{R}_+^p} \mathcal{L}_{{\gamma}}(y, z) = \max_{z \in \mathbb{R}_+^p} \min_{y} \mathcal{L}_{{\gamma}}(y, z) .
          \end{equation*}
    \noindent      
    2.
     The strong duality of the regularized augmented Lagrangian holds, i.e.,
          \[
             \min_{y} \max_{z \in \mathbb{R}_+^p} \mathcal{L}_{{\gamma}}(y, z) - \frac{\sigma}{2} \|z - z'\|^2  = \max_{z \in \mathbb{R}_+^p} \min_{y} \mathcal{L}_{{\gamma_1}}(y, z) - \frac{\sigma}{2} \|z - z'\|^2 .
          \]
          holds for any given $\sigma > 0$ and $z' \in \mathbb{R}_+^p$.
    \end{proposition}
    The proof of Proposition~\ref{prop: strong duality} is provided in Chapter 17 of~\citet{nocedal1999numerical}. A direct consequence of Proposition~\ref{prop: strong duality} is that  minimax in~\eqref{eq: minimax subproblem} is interchangeable.
    
    
    \subsection{Gradient oracles}
    
    The gradient of $\mathcal{L}_{{\gamma_1}}(x, y, z)$ is given by
    \begin{equation}
       \label{eq: derivative of L}
       \begin{aligned}
          \nabla_x \mathcal{L}_{{\gamma_1}}(x, y, z) &= \nabla_x g(x, y)   
          + \frac{1}{{\gamma_1}} \sum_{i=1}^{p}  [  {\gamma_1} z_i + H_{i}(x, y)]_+\nabla_x H_{i}(x, y) , \\
       \nabla_y \mathcal{L}_{{\gamma_1}}(x, y, z) &= \nabla_y g(x, y)   
       + \frac{1}{{\gamma_1}} \sum_{i=1}^{p}  [  {\gamma_1} z_i + H_{i}(x, y)]_+\nabla_y H_{i}(x, y) , \\
       \nabla_z \mathcal{L}_{{\gamma_1}}(x, y, z) &=  [{\gamma_1} z + H(x, y)]_+ - {\gamma_1} z = \max(- {\gamma_1} z, H(x, y)) .
       \end{aligned}
    \end{equation}
    With some simple computation, it can be shown that these gradients are bounded by linear functions of $\|z\|$.
    \begin{lemma}
       \label{lem: bound of gradient of L}
    Under Assumption~\ref{ass: boundedness}, the gradients oracle $\nabla_{1} \mathcal{L}_{{\gamma_1}}(x, y, z), \nabla_{2} \mathcal{L}_{{\gamma_1}}(x, y, z)$ are bounded by
       \begin{subequations}
          \begin{align}
          \|\nabla_x \mathcal{L}_{{\gamma_1}}(x, y, z) \|
             &\leq M_{G, 1} + p \left(\frac{M_{H, 0}M_{H, 1}}{\gamma_1} + \|z\|\right), 
        \label{eq: bound of gradient of Lk 1}\\
          \|\nabla_y \mathcal{L}_{{\gamma_1}}(x, y, z) \|
          &\leq  M_{G, 1} + p \left(\frac{M_{H, 0}M_{H, 1}}{\gamma_1} + \|z\|\right).
          \label{eq: bound of gradient of Lk 2}
          \end{align}
       \end{subequations}
    \end{lemma}
    
    {\it Proof}
    It follows from~\eqref{eq: derivative of L} that
        \begin{equation*}
            \begin{aligned}
                 \|\nabla_x \mathcal{L}_{{\gamma_1}}(x, y, z) \| &\leq \| \nabla_x g(x, y) \| 
           + \frac{1}{{\gamma_1}} \sum_{i=1}^{p}  ( {\gamma_1} | z_i | +  | H_{i}(x, y) | ) \| \nabla_y H_{i}(x, y)\|  \\
          &\leq \| \nabla G(x, y; \xi) \| + \frac{1}{{\gamma_1}} \sum_{i=1}^{p}  ( {\gamma_1} | z_i | +  M_{H, 0})M_{H, 1} 
          \\
          &\leq M_{G, 1} + p \left(\frac{M_{H, 0}M_{H, 1}}{\gamma_1} + \|z\|\right).
            \end{aligned}
        \end{equation*}
        The proof of \eqref{eq: bound of gradient of Lk 2} is the same.
    \qed
    By substituting $\nabla G(x, y)$ with $\nabla g(x, y; \xi)$ in~\eqref{eq: derivative of L}
    , we derive the gradient oracle of $\mathcal{L}_{{\gamma_1}}(x, y, z; \xi)$ as
    \begin{equation}
       \label{eq: gradient oracle of L}
       \begin{aligned}
          \nabla_x \mathcal{L}_{{\gamma_1}}(x, y, z; \xi) &= \nabla_x g(x, y; \xi)
          + \frac{1}{{\gamma_1}} \sum_{i=1}^{p}  [  {\gamma_1} z_i + H_{i}(x, y)]_+\nabla_x H_{i}(x, y) , \\
       \nabla_y \mathcal{L}_{{\gamma_1}}(x, y, z; \xi) &= \nabla_y g(x, y; \xi)
       + \frac{1}{{\gamma_1}} \sum_{i=1}^{p}  [  {\gamma_1} z_i + H_{i}(x, y)]_+\nabla_y H_{i}(x, y) , \\
       \nabla_z \mathcal{L}_{{\gamma_1}}(x, y, z; \xi) &=  [{\gamma_1} z + H(x, y)]_+ - {\gamma_1} z = \max(- {\gamma_1} z, H(x, y)) .
       \end{aligned}
    \end{equation}
    We introduce several properties of the augmented Lagrangian function’s gradient oracle that are essential for the analysis of our approach. The following lemma guarantees that the stochastic gradient remains within bounded norms, which is a crucial condition for stability.
    \begin{lemma}
       \label{lem: bound of first order oracle of L with respect to y}
    Under Assumption~\ref{ass: boundedness} and~\ref{ass: stochastic oracle}, the gradient oracles $\nabla_{1} \mathcal{L}_{{\gamma_1}}(x, y, z; \xi), \nabla_{2} \mathcal{L}_{{\gamma_1}}(x, y, z; \xi)$ are bounded by
       \begin{subequations}
          \begin{align}
             \mathbb{E}_{\xi \sim \mathcal{D}_{\xi}} [\|\nabla_x \mathcal{L}_{{\gamma_1}}(x, y, z; \xi) \|^2 ] 
             &\leq M_{\mathcal{L}, 1} + M_{\mathcal{L}, 2} \|z\|^2, 
        \label{eq: bound of first order oracle of Lk 1}\\
          \mathbb{E}_{\xi \sim \mathcal{D}_{\xi}} [\|\nabla_y \mathcal{L}_{{\gamma_1}}(x, y, z; \xi) \|^2 ]
          &\leq M_{\mathcal{L}, 1} + M_{\mathcal{L}, 2} \|z\|^2,
          \label{eq: bound of first order oracle of Lk 2}
          \end{align}
       \end{subequations}
       respectively,
    where 
    $
        M_{\mathcal{L}, 1} = (2 + p ) (M_{G, 1}^2 +\sigma_{g}^2  + \frac{1}{\gamma_1^2} M_{H, 0}) ,  M_{\mathcal{L}, 2} = (2 + p )M_{H, 0}^2M_{H, 1}^2 .
    $
    \end{lemma}
    {\it Proof}
    It follows from~\eqref{eq: gradient oracle of L} that
        \begin{equation*}
            \begin{aligned}
                 \|\nabla_x \mathcal{L}_{{\gamma_1}}(x, y, z; \xi) \| &\leq \| \nabla_x g(x, y; \xi) \| 
           + \frac{1}{{\gamma_1}} \sum_{i=1}^{p}  ( {\gamma_1} | z_i | +  | H_{i}(x, y) | ) \| \nabla_y H_{i}(x, y)\|  \\
          &\leq \| \nabla_x g(x, y; \xi) \| + \frac{1}{{\gamma_1}} \sum_{i=1}^{p}  ( {\gamma_1} | z_i | +  M_{H, 0})M_{H, 1} .
            \end{aligned}
        \end{equation*}
        By Cauchy-Schwarz inequality, it holds that 
        \begin{equation}
          \label{eq: proof of bound of first order oracle of Lk 1}
        \begin{aligned}
            &\quad \left(\| \nabla_x g(x, y; \xi) \|  + \frac{1}{{\gamma_1}} \sum_{i=1}^{p}  ( {\gamma_1} | z_i | +  M_{H, 0}) M_{H, 1} \right)^2 \\
            &\leq 
            (1 + 1 + p ) \left( \| \nabla_x g(x, y; \xi) \|^2 + \frac{1}{\gamma_1^2} M_{H, 0}^2 +  M_{H, 0}^2 M_{H, 1}^2 \sum_{i=1}^{p} z_i^2 \right).
        \end{aligned}
        \end{equation}
        Assumption~\ref{ass: boundedness} and~\ref{ass: stochastic oracle} imply 
        $$
          \mathbb{E}_{\xi \sim \mathcal{D}_{\xi}} [\|\nabla_x g(x, y; \xi)\|^2]
          = \mathbb{E}_{\xi \sim \mathcal{D}_{\xi}} [\|\nabla_x g(x, y)\|^2]
          + \mathbb{V}_{\xi \sim \mathcal{D}_{\xi}} [\nabla_x g(x, y; \xi)]
          \leq M_{G, 1}^2 + \sigma_{g}^2.
        $$
        Taking expectation on~\eqref{eq: proof of bound of first order oracle of Lk 1} gives~\eqref{eq: bound of first order oracle of Lk 1}. The proof of~\eqref{eq: bound of first order oracle of Lk 2} is similar. 
    \qed
    
    The convexity of $G(x, y)$ and $H(x, y)$ implies the convex-concavity of $\mathcal{L}_{{\gamma_1}}(x, y, z)$ as shown in the following lemma.
    \begin{lemma}
       \label{lem: convex-concav}
       \textbf{(convex-concavity)}
       The function $\mathcal{L}_{{\gamma_1}}(x, y, z)$ is $\mu_G$-strongly convex in $y$ and concave in $z$.
    \end{lemma}
    {\it Proof}
       From~\eqref{eq: derivative of L} we can compute the second order derivative with respect to $y$ as 
       \begin{equation}
       \label{eq: second order derivative of L}
          \begin{aligned}
      \nabla_y^2 \mathcal{L}_{{\gamma_1}}(x, y, z) = &\nabla_y^2 G(x, y) + \frac{1}{{\gamma_1}} \sum_{i=1}^{p} \{ [{\gamma_1} z_i + H_{i}(x, y)]_+ \nabla_y^2 H_{i}(x, y) \\
       & + \mathbb{I}_{\{ {\gamma_1} z_i + H_{i}(x, y) > 0\}} \nabla_y H_{i}(x, y) \nabla_y H_{i}(x, y)^T \}   \succeq \mu_G I .
       \end{aligned}
       \end{equation}
       This implies that $\mathcal{L}_{{\gamma_1}}(x, y, z)$ is $\mu_G$-strongly convex in $y$. Additionally~\eqref{eq: derivative of L}
       shows that $\nabla_z \mathcal{L}_{\gamma_1}(x, y, z)$ is monotonically decreasing with respect to $z$. Therefore $ \mathcal{L}_{\gamma_1}(x, y, z)$ is concave in $z$.
    \qed

    \begin{lemma}
       \label{lem: convex-concav of tilde L}
       The function $\ell_{\gamma}(x, y, w, \lambda)$ is $\mu_G$-strongly convex in $w$ and $\gamma_2$-strongly concave in $\lambda$.
    \end{lemma}
    {\it Proof}
       Combining Lemma~\ref{lem: convex-concav} and \eqref{eq: tilde L}, the conclusion follows.
    \qed
    Besides, its gradients can be computed as
    \begin{equation}
       \label{eq: derivative of tilde L}
       \begin{aligned}
       \nabla_w \ell_{\gamma}(x, z, w, \lambda) &= \nabla_y g(x, w)   
       + \frac{1}{{\gamma_1}} \sum_{i=1}^{p}  [  {\gamma_1} \lambda_i + H_{i}(x, w)]_+\nabla_y H_{i}(x, w)  , \\
       \nabla_{\lambda} \ell_{\gamma}(x, z, w, \lambda) &=  [{\gamma_1} \lambda + H(x, y)]_+ - {\gamma_1} \lambda - \gamma_2 (\lambda - z) = \max(- {\gamma_1} \lambda, H(x, y)) - \gamma_2 (\lambda - z) .
       \end{aligned}
    \end{equation}
    
    \subsubsection{Comparison between Lagrangian function and augmented Lagrangian function} \label{sec: comparison of Lagrangian and augmented Lagrangian}
    In this section, we provide a comprehensive comparison between the Lagrangian function and the augmented Lagrangian function. The key points are summarized as follows: first, the Lagrangian function is a special case of the augmented Lagrangian function when  $\gamma = +\infty$; second, the augmented Lagrangian function-based reformulation offers a tighter upper bound for the variance estimation when computing the upper-level gradient.
    The Lagrangian function $\mathcal{L}(x, y, z)$ of~\eqref{eq: LL}
    is defined as
    \begin{equation}
       \label{eq: Lagrangian function}
          \mathcal{L}(x, y, z) =  G(x, y) + \sum_{i=1}^{p} z_i H_i(x, y), \quad (x, y, z) \in X \times Y \times \mathbb{R}_+^p.
    \end{equation}
    The inequality relation between objective function, Lagrangian function and augmented Lagrangian function is given in the following proposition.
    \begin{proposition}
       \label{prop: augmented Lagrangian}
       If $H(x, y) \leq 0, z \geq 0$, then it holds that
       \[
          \mathcal{L}(x, y, z) \leq \mathcal{L}_{{\gamma_1}}(x, y, z) \leq G(x, y) .
       \]
    \end{proposition}
    {\it Proof}
    For any $i \in \{1, ..., p\}$, there are two cases:
    
    \noindent
    1.
       If ${\gamma_1} z_i + H_i(x, y) \leq 0$, then 
           $
       \frac{1}{2{\gamma_1}} ([{\gamma_1} z_i + H_i(x, y)]_+^2 - {\gamma_1}^2 z_i^2) = - \frac{1}{2} {\gamma_1} z_i^2 
       $.
       It holds that $  z_i H_i(x, y) \leq - \frac{1}{2} {\gamma_1} z_i^2  \leq 0$.
       
    \noindent 
    2.
          If  ${\gamma_1} z_i + H_i(x, y) > 0$, then
       $
       \frac{1}{2{\gamma_1}} ([{\gamma_1} z_i + H_i(x, y)]_+^2 - {\gamma_1}^2 z_i^2) = z_i H_i(x, y) + \frac{1}{2{\gamma_1}} H_i^2(x, y)
       $. Hence $z_i H_i(x, y) \leq z_i H_i(x, y) + \frac{1}{2{\gamma_1}} H_i^2(x, y) \leq 0$.
    
       Combining the above two cases, we have $\mathcal{L}(x, y, z) \leq \mathcal{L}_{{\gamma_1}}(x, y, z) \leq G(x, y)$. This completes the proof.
    \qed

    From equations~\eqref{eq: AL} and~\eqref{eq: Lagrangian function}, it is evident that the Lagrangian function is the limit of the augmented Lagrangian function as $ \gamma_1 = +\infty$. During the optimization process, it is desirable for the variance of the gradient oracle to be small in order to ensure the algorithm's stability. The second-order moment serves as an upper bound for the variance. The following lemma demonstrates that the gradient of the augmented Lagrangian term for each constraint has a smaller second-order moment compared to that of the Lagrangian function. This property is a key reason for using the augmented Lagrangian function in our algorithm.
    From equations~\eqref{eq: AL} and~\eqref{eq: Lagrangian function}, it is evident that the Lagrangian function is the limit of the augmented Lagrangian function as $ \gamma_1 = +\infty$. During the optimization process, it is desirable for the variance of the gradient oracle to be small in order to ensure the algorithm's stability. The second-order moment serves as an upper bound for the variance. The following lemma demonstrates that the gradient of the augmented Lagrangian term for each constraint has a smaller second-order moment compared to that of the Lagrangian function. This property is a key reason for using the augmented Lagrangian function in our algorithm.

    \begin{lemma}
       Assume $x$ is fixed and the random variable $(\hat w(\mathbf{\xi}), \hat \lambda(\mathbf{\xi})) \in \mathcal{Y}(x) \times \mathbb{R}_+^p$, $\mathbf{\xi} \in \Omega_{\xi}^s$.
      Then it holds that
       \begin{equation}
          \begin{aligned}
             \mathbb{E}_{\mathbf{\xi}} \left[ \left\| \frac{1}{{\gamma_1}}  [  {\gamma_1} \hat \lambda_i(\mathbf{\xi}) + H_{i}(x, \hat w(\mathbf{\xi}))]_+ \nabla H_{i}(x, \hat w(\mathbf{\xi})) \right\|^2 \right] \leq \mathbb{E}_{\mathbf{\xi}} \left[\left\| \hat \lambda_i(\mathbf{\xi}) \nabla H_{i}(x, \hat w(\mathbf{\xi})) \right\|^2 \right].
          \end{aligned}
       \end{equation}
    \end{lemma}
    
    {\it Proof}
       The condition
    $(\hat w(\mathbf{\xi}), \hat \lambda(\mathbf{\xi})) \in \mathcal{Y}(x) \times \mathbb{R}_+^p$ indicates that $H(x, \hat w(\mathbf{\xi})) \leq 0$ and $\hat \lambda(\mathbf{\xi}) \geq 0$. Then we have
    \begin{equation*}
       \begin{aligned}
          0 \leq \frac{1}{{\gamma_1}}  [  {\gamma_1} \hat \lambda_i(\mathbf{\xi}) + H_{i}(x, \hat w(\mathbf{\xi}))]_+  \leq \hat \lambda_i(\mathbf{\xi}), 
          \quad i = 1, ..., p.
       \end{aligned}
    \end{equation*}
    Taking square and then taking expectation give the desired result.
    \qed
    \begin{lemma}
       Assume $x$ is fixed and the random variable $(\hat w(\mathbf{\xi}), \hat \lambda(\mathbf{\xi})) \in \mathcal{Y}(x) \times \mathbb{R}_+^p$, $\mathbf{\xi} \in \Omega_{\xi}^s$.
      Then it holds that
       \begin{equation}
          \begin{aligned}
             \mathbb{E}_{\mathbf{\xi}} \left[ \left\| \frac{1}{{\gamma_1}}  [  {\gamma_1} \hat \lambda_i(\mathbf{\xi}) + H_{i}(x, \hat w(\mathbf{\xi}))]_+ \nabla H_{i}(x, \hat w(\mathbf{\xi})) \right\|^2 \right] \leq \mathbb{E}_{\mathbf{\xi}} \left[\left\| \hat \lambda_i(\mathbf{\xi}) \nabla H_{i}(x, \hat w(\mathbf{\xi})) \right\|^2 \right].
          \end{aligned}
       \end{equation}
    \end{lemma}
    
    {\it Proof}
       The condition
    $(\hat w(\mathbf{\xi}), \hat \lambda(\mathbf{\xi})) \in \mathcal{Y}(x) \times \mathbb{R}_+^p$ indicates that $H(x, \hat w(\mathbf{\xi})) \leq 0$ and $\hat \lambda(\mathbf{\xi}) \geq 0$. Then we have
    \begin{equation*}
       \begin{aligned}
          0 \leq \frac{1}{{\gamma_1}}  [  {\gamma_1} \hat \lambda_i(\mathbf{\xi}) + H_{i}(x, \hat w(\mathbf{\xi}))]_+  \leq \hat \lambda_i(\mathbf{\xi}), 
          \quad i = 1, ..., p.
       \end{aligned}
    \end{equation*}
    Taking square and then taking expectation give the desired result.
    \qed
    


    \subsection{Analysis on the reformulation} \label{sec: reformulation}
    
    In this section, we show the equivalence between the reformulations and the original BLO~\eqref{eq: BLO}. First we consider the deterministic case as a special case of the stochastic case. Then we extend the equivalence to the stochastic case and show the major improvements of the stochastic reformulation.

    \subsection{Deterministic case}
    To further analyze the optimization problem, we establish a critical property of the lower-level objective function $G(x, y)$. 
    Define $v(x) = \min_{y \in \mathcal{Y}(x)} G(x, y)$.
    The following proposition shows that $G(x, y)$ satisfies a quadratic growth condition, ensuring the uniqueness of solutions in the lower-level problem.

    \begin{proposition}
       \label{prop: quadratic growth}
       \textbf{(Quadratic growth)}
       1.
       Suppose Assumption~\ref{ass: convexity} holds.
       For any $x \in X$, 
       $G(x, y)$ has $\mu_G$-quadratic growth with respect to $y \in \mathcal{Y}(x)$, namely,
       \begin{equation*}
          \begin{aligned}
          G(x, y) - v(x) &\geq \frac{\mu_G}{2} \| y - y^*(x) \|^2, \quad  y \in \mathcal{Y}(x). 
          \end{aligned}
       \end{equation*}
       \noindent
       2. Suppose Assumption~\ref{ass: convexity} and~\ref{ass: boundedness} holds.
       For any $x \in X$ and $y \in Y$ satisfying $\mathrm{dist}(y, \mathcal{Y}(x)) \leq \epsilon$, we have
       \begin{equation*}
          G(x, y) - v(x) \geq \frac{\mu_G}{4} \|y - y^*(x)\|^2  - M_{G, 1} \epsilon - \frac{\mu_G}{2} \epsilon^2 . 
       \end{equation*}
    
       3. Suppose Assumption~\ref{ass: convexity} and~\ref{ass: boundedness} holds. For any $x \in X$ and $y \in Y$, it holds that
       it holds that
       \begin{equation*}
          G(x, y) - v(x) \geq     \frac{\mu_G}{2} \|y - y^*(x)\|^2 - p^{\frac{1}{2}} B M_{H, 1} \| y - y^*(x) \|,
       \end{equation*}
       where $B$ is defined in Lemma~\ref{lem: equivalence between the  value function-based reformulation and the dual augmented value function reformulation}.3.
    Further assume $y$ satisfys $\frac{1}{2} \sum_{i=1}^{p} [H_i(x, y)]_+^2 \leq \epsilon^2$, then 
       \begin{equation*}
          G(x, y) - v(x) \geq   \frac{\mu_G}{2} \|y - y^*(x)\|^2 - \sqrt{2} B \epsilon .
       \end{equation*}
    
    \end{proposition}
    {\it Proof}
       For any $x \in X$, denote
       $$
       (y^*(x), \lambda^*(x)) = \arg \min_{y \in Y} \max_{\lambda \in \mathbb{R}_+^p} G(x, y) + \lambda^\top H(x, y).
       $$
       By Assumption~\ref{ass: convexity}, the strong convexity in $y$ implies its quadratic growth, i.e., 
       $$
       G(x, y) + (\lambda^*(x))^\top H(x, y) \geq G(x, y^*(x)) + \lambda^*(x)^\top H(x, y^*(x)) + \frac{\mu_G}{2} \|y - y^*(x)\|^2 .
       $$
       Furthermore, the complementary slackness condition implies $\lambda^*_i H_i(x, y^*(x)) = 0$ for $i = 1, ..., p$. Therefore, we have
       \begin{equation*}
          \begin{aligned}
            &\quad  G(x, y) - v(x) \\
              &\geq G(x, y) + (\lambda^*(x))^\top H(x, y) - G(x, y^*(x)) \quad \text{(since  $\lambda^*(x) \geq 0$ and $H(x, y) \leq 0$)} \\
             &= G(x, y) + (\lambda^*(x))^\top H(x, y) - G(x, y^*(x)) - \lambda^*(x)^\top H(x, y^*(x)) \\
             &\geq \frac{\mu_G}{2} \|y - y^*(x)\|^2 , \quad y \in \mathcal{Y}(x).
          \end{aligned}
       \end{equation*}
       2. 
       For any $y \in Y$ satisfying $\mathrm{dist}(y, \mathcal{Y}(x)) \leq \epsilon$, denote 
       $y' = \mathrm{Proj}_{\mathcal{Y}(x)}(y)$ as the projection of $y$ onto $\mathcal{Y}(x)$. Then $\| y - y' \| \leq \epsilon$. From the results of 1, we have 
       $$
       G(x, y') - v(x) \geq  \frac{\mu_G}{2}  \|y' - y^*(x)\|^2 .
       $$
       Assumption~\ref{ass: boundedness} implies
       $$
       G(x, y) - G(x, y') \geq  - M_{G, 1} \| y - y'\|.
       $$
       Combining  the above two inequalities and the fact $\frac{1}{2} \| y - y^*(x) \|^2  \leq \| y' - y^*(x)\|^2 + \| y - y'\|^2$, we have
       $$
       \begin{aligned}
          G(x, y) -v(x) &=  G(x, y) - G(x, y')  + G(x, y') - v(x)  \\
          &\geq - M_{G, 1} \| y - y'\| + \frac{\mu_G}{2}  \|y' - y^*(x)\|^2 \\
          &\geq - M_{G, 1} \epsilon + \frac{\mu_G}{2} \left( \frac{1}{2} \|y - y^*(x)\|^2 - \| y - y'\|^2 \right) \\
          &\geq  \frac{\mu_G}{4} \|y - y^*(x)\|^2  - M_{G, 1} \epsilon - \frac{\mu_G}{2} \epsilon^2 .
       \end{aligned}
       $$
    
       3.
       By Assumption~\ref{ass: boundedness} and Lemma~\ref{lem: equivalence between the  value function-based reformulation and the dual augmented value function reformulation}, for any $x \in X$ and $y \in Y$, it holds that
       $$
       \lambda^*(x) H(x, y)  \leq \| \lambda^*(x) \| [\| H(x, y)]_+ \| \leq p^{\frac{1}{2}} B M_{H, 1} \| y - y^*(x) \|.
       $$
       Then following the proof in Lemma~\ref{prop: quadratic growth}.1, we have
       \begin{equation*}
          \label{eq: proof of quadratic growth 3}
          \begin{aligned}
             G(x, y) - v(x) 
             &\geq \frac{\mu_G}{2} \|y - y^*(x)\|^2 - \lambda^*(x) H(x, y)  \\
             &\geq \frac{\mu_G}{2} \|y - y^*(x)\|^2 - p^{\frac{1}{2}} B M_{H, 1} \| y - y^*(x) \| .
          \end{aligned}
       \end{equation*}
    
       Futher assume $y \in Y$ satisfys $\frac{1}{2}\sum_{i=1}^{p} [H_i(x, y)]_+^2 \leq \epsilon^2$, then 
       \begin{equation*}
          \begin{aligned}
             G(x, y) - v(x) 
             &\geq \frac{\mu_G}{2} \|y - y^*(x)\|^2 - \lambda^*(x) H(x, y)  \\
             &\geq \frac{\mu_G}{2} \|y - y^*(x)\|^2 - \| \lambda^*(x) \| [\| H(x, y)]_+ \|   \\
             &\geq \frac{\mu_G}{2} \|y - y^*(x)\|^2 - \sqrt{2}B \epsilon .
          \end{aligned}
       \end{equation*}
    This completes the proof.
    \qed

    The following lemma shows the relationship between $v(x)$ and the enveloped value $E(x, z)$ defined in \eqref{eq: augmented dual function}.
    \begin{lemma}
       \label{lem: equivalence between the  value function-based reformulation and the dual augmented value function reformulation}
       Suppose Assumptions~\ref{ass: convexity},~\ref{ass: LL Slater's condition} and~\ref{ass: LICQ} hold.
    
       \noindent
       1. It holds that $E(x, z) \leq v(x), \forall z \in \mathbb{R}_+^{p}$. The equality holds if and only if there exist $y \in \mathcal{Y}(x)$ such that $(y, z) \in \arg \min_{w \in Y} \max_{\lambda \in \mathbb{R}_+^p} \mathcal{L}_{\gamma_1}(x, w, \lambda)$.
       
       \noindent
       2. Assume $y \in \mathcal{Y}(x)$, then 
       $y \in \arg \min_{\mathcal{Y}(x)} G(x, y)$ is equivalent to existing $z \in \mathbb{R}_+^p$ such that $G(x, y) - E(x, z) \leq 0$.
    
       3. For any $x \in X$, there exists $z$ such that $(y, z) \in \arg \min_{w \in Y} \max_{\lambda \in \mathbb{R}_+^p}\mathcal{L}_{\gamma_1}(x, y, z)$ satisfying
       $|z_i| \leq p^{-0.5} B$, where $B = p^{2} \sigma_0^{-2} M_{H, 1}( M_{G, 1} + pM_{H, 1})$.
       This statement also holds for $\gamma_1 = +\infty$ .
    
    \end{lemma}
    {\it Proof}
       1. For any fixed $(x, y)$, we have
       $$
       \frac{1}{2{{\gamma_1}}} ([{\gamma_1} z_i + H_{i}(x, y)]_+^2 - {\gamma_1}^2 z_i^2) = \begin{cases}
          z_i H_i(x, y) + \frac{1}{2{{\gamma_1}}} H_{i}(x, y)^2,
          & \text{if } z_i \geq -\frac{1}{\gamma_1} H_{i}(x, y), \\
          - \frac{1}{2} \gamma_1 z_i^2, & \text{if } z_i < -\frac{1}{\gamma_1} H_{i}(x, y).
       \end{cases}
       $$
       We consider maximizing with respect to $z_i \in \mathbb{R}_+$. If $H_{i}(x, y) > 0$, then 
       $\max_{z_i \in \mathbb{R}_+} \frac{1}{2{{\gamma_1}}} ([{\gamma_1} z_i + H_{i}(x, y)]_+^2 - {\gamma_1}^2 z_i^2) = \max_{z_i \in \mathbb{R}_+} z_i H_i(x, y) + \frac{1}{2{{\gamma_1}}} H_{i}(x, y)^2 = +\infty$.
       If $H_{i}(x, y) \leq 0$, then
       $$
       \begin{aligned}
       &\quad \max_{z_i \in \mathbb{R}_+}\frac{1}{2{{\gamma_1}}} ([{\gamma_1} z_i + H_{i}(x, y)]_+^2 - {\gamma_1}^2 z_i^2) \\
       &= 
       \max \left\{ \max_{0 \leq z_i \leq -\frac{1}{\gamma_1} H_{i}(x, y)} - \frac{1}{2{{\gamma_1}}} \gamma_1^2 z_i^2, 
       \max_{z_i \geq -\frac{1}{\gamma_1} H_i(x, y)} z_i H_i(x, y) + \frac{1}{2{{\gamma_1}}} H_{i}(x, y)^2 \right\} \\
       &= \max \left\{ 0, - \frac{1}{2{{\gamma_1}}} H_{i}(x, y)^2 \right\} = 0.
       \end{aligned}
       $$
       Combining the above two cases, we have
       $$
       \begin{aligned}
       \max_{z_i \in \mathbb{R}_+}\frac{1}{2{{\gamma_1}}} ([{\gamma_1} z_i + H_{i}(x, y)]_+^2 - {\gamma_1}^2 z_i^2) 
       &= \begin{cases}
          +\infty, & \text{if } H_{i}(x, y) > 0, \\
          0 , & \text{if } H_{i}(x, y) \leq 0.
       \end{cases} 
       \end{aligned}
       $$
       This implies that
       \begin{equation}
       \label{eq: max of augmented Lagrangian with respect to z}
       \max_{z \in \mathbb{R}_+^p} \mathcal{L}_{\gamma_1}(x, y, z) = \begin{cases}
          +\infty, & \text{if } \exists H_i(x, y) > 0, \\
          G(x, y) , & \text{if } H_i(x, y) \leq 0,  i=1, ..., p.
       \end{cases}
    \end{equation}

    From the definition of $E(x, z)$, for $\forall z \in \mathbb{R}_+^p$, we have 
       \begin{equation}
          \begin{aligned}
             E(x, z) 
             &= \max_{\lambda \in \mathbb{R}_+^p} \min_{y \in Y} \left\{ \mathcal{L}_{\gamma_1}(x, y, \lambda) - \frac{\gamma_2}{2} \| \lambda - z \|^2 \right\} \\
             &= \min_{y \in Y} \max_{\lambda \in \mathbb{R}_+^p} \left\{ \mathcal{L}_{\gamma_1}(x, y, \lambda) - \frac{\gamma_2}{2} \| \lambda - z \|^2 \right\}  \quad \text{(by Proposition~\ref{prop: strong duality})}\\
             &\leq \min_{y \in Y} \max_{\lambda \in \mathbb{R}_+^p} \mathcal{L}_{\gamma_1}(x, y, \lambda) \\
             &= \min_{y \in Y} \left\{ G(x, y) | H(x, y) \leq 0 \right\} \quad \text{(by~\eqref{eq: max of augmented Lagrangian with respect to z})} \\
             &= v(x) .
          \end{aligned}
       \end{equation}
       The equality in the inequality holds if and only if $(y, z) \in \arg \min_{y \in Y} \max_{z \in \mathbb{R}_+^p} \mathcal{L}_{\gamma_1}(x, y, z)$.

       \noindent
       2. Since $E(x, z) \leq v(x) \leq G(x, y), \forall z \in \mathbb{R}_+^p, \forall y \in \mathcal{Y}(x)$, 
       hence the condition $G(x, y) - E(x, z) \leq 0$ holds if and only if
       $G(x, y) = v(x) = E(x, z)$, which implies $y \in \arg \min_{y \in \mathcal{Y}(x)} G(x, y)$. \\
       Conversely, suppose $y \in \arg \min_{\mathcal{Y}(x)} G(x, y)$,  Take $z \in \arg \max_{\lambda \in \mathbb{R}_+^p} \mathcal{L}_{\gamma_1}(x, y, \lambda)$, then
       \begin{equation}
          \begin{aligned}
             E(x, z) 
             &= \min_{w \in Y} \max_{\lambda \in \mathbb{R}_+^p} \left\{ \mathcal{L}_{\gamma_1}(x, w, \lambda) - \frac{\gamma_2}{2} \| \lambda - z \|^2 \right\}\\
             &\geq \min_{w \in Y} \mathcal{L}_{\gamma_1}(x, w, z) \quad \text{(taking $\lambda = z$)}\\
             &= \min_{w \in \mathcal{Y}(x)} G(x, w) = G(x, y) .
          \end{aligned}
       \end{equation}
    
       3. The optimal condition of $(y, z) \in \arg  \min_{w \in Y} \max_{\lambda \in \mathbb{R}_+^p} \mathcal{L}_{\gamma_1}(x, y, z)$ implies that $H(x, y) \leq 0$ and
       $\nabla_y \mathcal{L}_{\gamma_1}  (x, y, z) = 0$, that is,
       \begin{equation}
       \label{eq: proof of equivalence between the  value function-based reformulation and the dual augmented value function reformulation 3}
       \nabla_y g(x, y) + \frac{1}{{\gamma_1}} \sum_{i=1}^{p}  [  {\gamma_1} z_i + H_{i}(x, y)]_+\nabla_y H_{i}(x, y) = 0.
    \end{equation}
       From the proof of 1, we know that if $H_i(x, y) < 0$, then $z_i = 0$. 
       If $H_i(x, y) = 0$, then $\nabla_y H_{i}(x, y) \neq 0$. Otherwise suppose $H_i(x, y) = 0$ and $\nabla_y H_{i}(x, y) = 0$ hold simultaneously. The convexity of $H_i(x, y)$ implies $H_i(x, y) \geq 0$ for all $y \in Y$, which contradicts Assumption~\ref{ass: LL Slater's condition}. Substituting these cases in~\eqref{eq: proof of equivalence between the  value function-based reformulation and the dual augmented value function reformulation 3} yields
       \begin{equation}
          \label{eq: proof of equivalence between the  value function-based reformulation and the dual augmented value function reformulation 5}
          \nabla_y g(x, y) + \sum_{i: H_i(x, y) = 0}  z_i\nabla_y H_{i}(x, y) +  \sum_{i: H_i(x, y) < 0}  [  H_{i}(x, y)]_+\nabla_y H_{i}(x, y) = 0.
       \end{equation}
       Note that the above conditions are equivalent to the KKT conditions.
       Since the lower-level problem is convex and Slater's conditions holds, this is a sufficient condition for the optimal point. Assume matrix $\mathcal{C}$ is the concatenation of $\nabla_y H_i(x, y) = 0, i \in \{i, H_i(x, y)\}$, then $z$ with the minimal norm satisfying \eqref{eq: proof of equivalence between the  value function-based reformulation and the dual augmented value function reformulation 5} has a bound as 
       $$
       | z_i | \leq \left\| (\mathcal{C}^\top \mathcal{C})^{-1} \mathcal{C}^\top \left(\nabla_y g(x, y) + \sum_{i: H_i(x, y) < 0}  [  H_{i}(x, y)]_+\nabla_y H_{i}(x, y)\right) \right\|
       \leq p^{1.5} \sigma_0^2 M_{H, 1}( M_{G, 1} + p M_{H,0}M_{H, 1})
       $$
       The inequality uses Assumption~\ref{ass: LICQ}.
       As a special case of $\gamma_1 = +\infty$, 
       $$
       | z_i | \leq \left\| (\mathcal{C}^\top \mathcal{C})^{-1} \mathcal{C}^\top \nabla_y G(x, y)  \right\|
       \leq p^{1.5} \sigma_0^{-2} M_{H, 1} M_{G, 1} .
       $$
       This completes the proof.
    \qed
    
    \begin{remark}
       Under the strong duality condition, $v(x) = \max_{z \in \mathbb{R}_+^p} D(x, z)$. Proposition~\ref{prop: properties of Moreau envelope}.1 and~\ref{prop: properties of Moreau envelope}.3 imply that $E(x, z) \leq \max_{z \in \mathbb{R}_+^p} E(x, z) = \max_{z \in \mathbb{R}_+^p} D(x, z) = v(x)$. This shows that the enveloped value function $E(x, z)$ is a lower approximation of $v(x)$ and does not change the optimal solution.
    \end{remark}

    The following theorem shows the equivalence between~\eqref{eq: BLO} and the penalty reformulation in the  deterministic setting. This theorem generalizes Theorem 1 of~\cite{shen2023penalty}; indeed, choosing $\epsilon_2 = 0$ reduces our statement to exactly that result.

    \begin{theorem}
       \label{lem: equivalence of single level 1}
       Suppose that Assumptions~\ref{ass: Lipschitz continuity} ~\ref{ass: convexity} and~\ref{ass: LICQ} holds and $\gamma_1, \gamma_2 > 0$ are fixed parameters. 
    
    \noindent
    1. Assume $(x^*, y^*)$ is a global solution to~\eqref{eq: BLO} and $c_1\geq  c_1^* = \frac{L}{\mu_G} \epsilon^{-1}$. 
    Then $(x^*, y^*)$ is a $\epsilon$-global-minima of the following penalized form
    \begin{equation}
       \label{eq: partly penalized single level}
    \begin{aligned}
       \min_{(x, y) \in X \times Y} &\quad F(x, y)  + c_1 (G(x, y) - v(x)) \quad  \text{s.t.}\quad  \frac{1}{2} \sum_{i=1}^p [H_i(x, y)]_+^2 \leq \epsilon_2^2.
    \end{aligned}
    \end{equation}
    Furthermore, there exist $z^* \in \mathbb{R}_+^p$ such that $(x^*, y^*, z^*)$ is a $\epsilon$-global-minima of the following penalized form
    \begin{equation}
       \label{eq: partly penalized single level 2}
    \begin{aligned}
       \min_{(x, y, z) \in X \times Y \times Z} &\quad  F(x, y)  + c_1 \mathcal{G}(x, y, z) \quad  \text{s.t.}\quad  \frac{1}{2} \sum_{i=1}^p [H_i(x, y)]_+^2 \leq \epsilon_2^2 ,
    \end{aligned}
    \end{equation}
    with $\epsilon_2 \leq \frac{\epsilon}{2 \sqrt{2} c_1 B}$.
    
    \noindent 
    2. By taking $c_1 = c_1^* + 2$, 
    any $\epsilon$-global-minima of~\eqref{eq: partly penalized single level} and~\eqref{eq: partly penalized single level 2} is an $\epsilon$-global-minima of the following two approximations of BLO:
    \begin{equation}
       \label{eq: approximation of BLO}
       \begin{aligned}
          \min_{ (x, y) \in X \times Y} &\quad  F(x, y) 
           \quad \text{s.t.} \quad   G(x, y) - v(x) \leq \epsilon_1, \quad  \frac{1}{2} \sum_{i=1}^p [H_i(x, y)]_+^2 \leq \epsilon_2^2.
       \end{aligned}
    \end{equation}
    and
    \begin{equation}
       \label{eq: approximation of BLO 2}
       \begin{aligned}
          \min_{ (x, y, z) \in X \times Y \times Z} &\quad  F(x, y) 
          \quad \text{s.t.} \quad  G(x, y) - E(x, z) \leq \epsilon_1, \quad  \frac{1}{2} \sum_{i=1}^p [H_i(x, y)]_+^2 \leq \epsilon_2^2,
       \end{aligned}
    \end{equation}
    with some $\epsilon_1, \epsilon_2 \leq \epsilon$.
    
    \end{theorem}
    {\it Proof}
    \noindent
    1.
    Lemma~\ref{lem: equivalence between the value function-based reformulation and the dual augmented value function reformulation} shows that~\eqref{eq: partly penalized single level 2} is equivalent to~\eqref{eq: partly penalized single level}
    by minimizing $z$ first.
    Hence it suffices to show the equivalence between~\eqref{eq: BLO} and~\eqref{eq: partly penalized single level}. Let $\mathrm{p}(x, y) = G(x, y) - v(x)$. Proposition~\ref{prop: quadratic growth}.3
    implies that $\mathrm{p}(x, y) \geq \frac{\mu_G}{2} \|y - y^*(x)\|^2 - \sqrt{2} B\epsilon_2 $ for any $y$ satisfying $\frac{1}{2} \sum_{i=1}^p [H_i(x, y)]_+^2 \leq \epsilon_2^2$. 
    By the Lipschitz property of $F(x, \cdot)$, we have
    \begin{equation}
       \label{eq: proof of equivalence of single level 1 1}
       \begin{aligned}
          F(x, y)  + c_1 \mathrm{p}(x, y)- F(x, y^*(x)) 
          &\geq - L_F \| y - y^*(x) \|   + \frac{c_1 \mu_G}{2} \|y - y^*(x)\|^2  - c_1  \sqrt{2} B \epsilon_2 \\
          &\geq -\frac{L_F}{2 c_1 \mu_G} - c_1  \sqrt{2} B \epsilon_2 .
       \end{aligned}
    \end{equation}
    By taking $c_1 \geq c_1^* =  \frac{L}{\mu_G} \epsilon^{-1}$ and $\epsilon_2 \leq \frac{\epsilon}{2 \sqrt{2} c_1 B}$,
    the right side of the above inequality $\geq -\epsilon$.
    For the optimal solution $(x^*, y^*)$ of~\eqref{eq: BLO} and $\forall (x, y) \in (X, Y)$, it holds that 
    \begin{equation}
    \label{eq: proof of equivalence of single level 1 2} 
    F(x^*, y^*)  + c_1 \mathrm{p}(x^*, y^*)  = F(x^*, y^*) \leq F(x, y^*(x)) \leq F(x, y)  + c_1 \mathrm{p}(x, y) + \epsilon.
    \end{equation}
    The first inequality follows from the optimality of $x^*$ and $y^* = y^*(x^*)$. The
    second inequality follows from~\eqref{eq: proof of equivalence of single level 1 1}.
    The inequality~\eqref{eq: proof of equivalence of single level 1 2} implies that $(x^*, y^*)$ is a $\epsilon$-optimal solution of~\eqref{eq: partly penalized single level}. Further take $z^* \in \arg \max_{\lambda \in \mathbb{R}_+^p} \mathcal{L}_{\gamma_1}(x^*, y^*, \lambda)$, then $(x^*, y^*, z^*)$ is a $\epsilon$-optimal solution of~\eqref{eq: partly penalized single level 2}. 
    By Lemma~\ref{lem: equivalence between the value function-based reformulation and the dual augmented value function reformulation}, we can find such $z^*$ satisfying $z \in Z$.
    This completes the proof.
       
    \noindent
    2. 
    Lemma~\ref{lem: equivalence between the value function-based reformulation and the dual augmented value function reformulation} shows that
    \eqref{eq: approximation of BLO 2} is equivalent to~\eqref{eq: approximation of BLO} by minimizing $z$ first.
    For any $\epsilon$-optimal solution $(\hat x, \hat y)$ of~\eqref{eq: partly penalized single level}, it holds that
    $$ 
    F (\hat x, \hat y) + c_1 \mathrm{p}(\hat x, \hat y) \leq F(x^*, y^*) + c_1 \mathrm{p}(x^*, y^*) + \epsilon \leq F(\hat x, \hat y) + c_1^*\mathrm{p}(\hat x, \hat y) + 2 \epsilon .$$
    The second inequality follows from~\eqref{eq: proof of equivalence of single level 1 1} and $\mathrm{p}(x^*, y^*) = 0$.
    Then we have
    \begin{equation}
       \mathrm{p}(\hat x, \hat y) \leq \frac{2\epsilon}{c_1 - c_1^*} = \epsilon.
    \end{equation}
    Take $\epsilon_1 = p(\hat x, \hat y), \epsilon_2 = \sqrt{ \frac{1}{2}\sum_{i=1}^p [H_i(\hat x, \hat y)]_+^2 }$.
    Then  $(\hat x, \hat y)$ is feasible for~\eqref{eq: approximation of BLO}. For any feasible solution $(x, y)$ of~\eqref{eq: approximation of BLO}, the $\epsilon$-optimality of $(\hat x, \hat y)$ implies
    \begin{equation}
       \begin{aligned}
          F(\hat x, \hat y) - F(x, y) 
          &\leq c_1 \left(\mathrm{p}(x, y) - \mathrm{p}(\hat x, \hat y) \right)  + \epsilon \\
          &= c_1 \left(\mathrm{p}(x, y) - \epsilon_1\right) + \epsilon \leq \epsilon.
       \end{aligned}
    \end{equation}
    The last inequality follows from the feasibility of $(x, y)$.
    Therefore $(\hat x, \hat y)$ is a $\epsilon$-optimal solution of~\eqref{eq: approximation of BLO}. This completes the proof.
    \qed

    \begin{theorem}
       \label{thm: equivalence of single level}
       Suppose that Assumptions~\ref{ass: Lipschitz continuity} ~\ref{ass: convexity} and~\ref{ass: LICQ} holds and $\gamma_1, \gamma_2 > 0$ are fixed parameters. 
    
    \noindent
    1. Assume $(x^*, y^*)$ is a global solution to~\eqref{eq: BLO} and $c_1\geq  c_1^* = \frac{L}{ 2 \mu_G} \epsilon^{-1}, c_2 \geq c_2^* = (c_1^*)^2 B^2\epsilon^{-1}$. 
    Then $(x^*, y^*)$ is a $\epsilon$-global-minima of the following penalized form
    \begin{equation}
       \label{eq: penalized single level}
    \begin{aligned}
       \min_{(x, y) \in X \times Y} &\quad  \Psi(x, y) = F(x, y)  + c_1 (G(x, y) - v(x)) + \frac{c_2}{2} \sum_{i=1}^p [H_i(x, y)]_+^2. 
    \end{aligned}
    \end{equation}
    Furthermore, there exist $z^* \in \mathbb{R}_+^p$ such that $(x^*, y^*, z^*)$ is a $\epsilon$-global-minima of the following penalized form
    \begin{equation}
       \label{eq: penalized single level 2}
    \begin{aligned}
       \min_{(x, y, z) \in X \times Y \times Z} &\quad  \Psi(x, y, z) = F(x, y)  + c_1 \mathcal{G}(x, y, z) + \frac{c_2}{2} \sum_{i=1}^p [H_i(x, y)]_+^2.
    \end{aligned}
    \end{equation}
    \noindent 
    2. By taking $c_1 = c_1^* + 2, c_2 = c_2^* + 2$, 
    any $\epsilon$-global-minima of~\eqref{eq: penalized single level} and~\eqref{eq: penalized single level 2} is an $\epsilon$-global-minima of 
    ~\eqref{eq: approximation of BLO} and~\eqref{eq: approximation of BLO 2}
    with some $\epsilon_1, \epsilon_2 \leq \epsilon$, respectively.
    
    \end{theorem}

    {\it Proof}
    \noindent
    1.
    Lemma~\ref{lem: equivalence between the value function-based reformulation and the dual augmented value function reformulation} shows that
    \eqref{eq: penalized single level 2} is equivalent to 
    \eqref{eq: penalized single level}
    by minimizing $z$ first.
    Hence it suffices to show the equivalence between~\eqref{eq: BLO} and~\eqref{eq: penalized single level}. Let $\mathrm{p}(x, y) = G(x, y) - v(x)$. Proposition~\ref{prop: quadratic growth}.3
    implies that $\mathrm{p}(x, y) \geq {\mu_G} \|y - y^*(x)\|^2 - \lambda^*(x) H(x, y) $.
    By the Lipschitz property of $F(x, \cdot)$, we have
    \begin{equation}
       \label{eq: proof of equivalence of single level 1}
       \begin{aligned}
          &\quad \Psi(x, y) - F(x, y^*(x)) \\
          &=F(x, y) - F(x, y^*(x)) + c_1 \mathrm{p}(x, y) + \frac{c_2}{2} \sum_{i=1}^{p} [H_i(x, y)]_+^2  \\
          &\geq - L_F \| y - y^*(x) \| +{c_1}\frac{\mu_G}{2} \|y - y^*(x)\|^2  - c_1 \lambda^*(x) H(x, y) + \frac{c_2}{2} \sum_{i=1}^{p} [H_i(x, y)]_+^2 \\
          &\geq -\frac{L_F}{2 c_1 \mu_G} - \frac{c_1^2}{2 c_2} \| \lambda^*(x)\|^2\\
          &\geq -\frac{L_F}{2 c_1 \mu_G} - \frac{c_1^2}{2 c_2} B^2 
           , \quad \forall (x, y) \in (X, Y).
       \end{aligned}
    \end{equation}
    By taking $c_1 \geq c_1^* =  \frac{L}{2 \mu_G} \epsilon^{-1}$ and $c_2 \geq c_2^* = (c_1^*)^2 B^2\epsilon^{-1}$, the right side of the above inequality $\geq -\epsilon$.
    For the optimal solution $(x^*, y^*)$ of~\eqref{eq: BLO}, it holds that 
    \begin{equation}
    \label{eq: proof of equivalence of single level 2}   
    \Psi(x^*, y^*) = F(x^*, y^*) \leq F(x, y^*(x)) \leq \Psi(x, y) + \epsilon, \quad \forall (x, y) \in (X, Y).
    \end{equation}
    The first inequality follows from the optimality of $x^*$ and the definition of $v(x)$. The second inequality follows from~\eqref{eq: proof of equivalence of single level 1}.
    This implies that $(x^*, y^*)$ is a $\epsilon$-optimal solution of~\eqref{eq: penalized single level}. Further take $z^* \in \arg \max_{\lambda \in \mathbb{R}_+^p} \mathcal{L}_{\gamma_1}(x^*, y^*, \lambda)$, then $(x^*, y^*, z^*)$ is a $\epsilon$-optimal solution of~\eqref{eq: penalized single level 2}. 
    By Lemma~\ref{lem: equivalence between the value function-based reformulation and the dual augmented value function reformulation}, we can find such $z^*$ satisfying $z \in Z$.
    This completes the proof.
       
    \noindent
    2. 
    Lemma~\ref{lem: equivalence between the value function-based reformulation and the dual augmented value function reformulation} shows that
    \eqref{eq: approximation of BLO 2} is equivalent to~\eqref{eq: approximation of BLO} by minimizing $z$ first.
    For any $\epsilon$-optimal solution $(\hat x, \hat y)$ of~\eqref{eq: penalized single level}, it holds that
    $$
    F(\hat x, y)  + c_1 \mathrm{p}(\hat x, \hat y) + \frac{c_2}{2} \sum_{i=1}^p [H_i(\hat x, \hat y)]_+^2 \leq \Psi(x^*, y^*) + \epsilon \leq F(\hat x, \hat y) + c_1^*\mathrm{p}(\hat x, \hat y) + \frac{c_2^*}{2} \sum_{i=1}^p [H_i(\hat x, \hat y)]_+^2 + 2 \epsilon .
    $$
    The second inequality follows from~\eqref{eq: proof of equivalence of single level 1}.
    From the selection of $c_1, c_2$, we have
    \begin{equation}
       \mathrm{p}(\hat x, \hat y) + \frac{1}{2}\sum_{i=1}^p [H_i(\hat x, \hat y)]_+^2 \leq \frac{2\epsilon}{2} = \epsilon.
    \end{equation}
    Take $\epsilon_1 = \mathrm{p}(\hat x, \hat y)\leq \epsilon, \epsilon_2 = \frac{1}{2}\sum_{i=1}^p [H_i(\hat x, \hat y)]_+^2 \leq \epsilon$.
    This implies $(\hat x, \hat y)$ is feasible for~\eqref{eq: approximation of BLO}. For any feasible solution $(x, y)$ of~\eqref{eq: approximation of BLO}, the $\epsilon$-optimality of $(\hat x, \hat y)$ implies
    \begin{equation}
       \begin{aligned}
          F(\hat x, \hat y) - F(x, y) 
          &\leq c_1 \left(\mathrm{p}(x, y) - \mathrm{p}(\hat x, \hat y) \right) + c_2 \left( \frac{1}{2} \sum_{i=1}^p [H_i(x, y)]_+^2 - \frac{1}{2}\sum_{i=1}^p [H_i(\hat x, \hat y)]_+^2 \right) + \epsilon \\
          &= c_1 \left(\mathrm{p}(x, y) - \epsilon_1\right) + c_2 \left(\frac{1}{2} \sum_{i=1}^p [H_i(x, y)]_+^2 - \epsilon_2 \right) + \epsilon \leq \epsilon.
       \end{aligned}
    \end{equation}
    The last inequality follows from the feasibility of $(x, y)$.
    Therefore $(\hat x, \hat y)$ is a $\epsilon$-optimal solution of~\eqref{eq: approximation of BLO}. This completes the proof.
    \qed
    
    \subsection{Stochastic case}
    Given $s$ samples $\mathbf{\xi} = (\xi_1, ..., \xi_s) \sim \mathcal{D}_{\xi}^s$. 
    Denote $\mathcal{P}_{w}, \mathcal{P}_{\lambda}$ as the space of random variables mapping $\mathbf{\xi}$ to $Y$ and $\mathbb{R}_+^p$,
     respectively, that is, 
    $$
    \begin{aligned}
       \mathcal{P}_{w} &= \{w: \Omega_{\xi}^s \to Y \mid w \text{ is measurable}\}, \quad 
       \mathcal{P}_{\lambda} = \{\lambda: \Omega_{\xi}^s \to \mathbb{R}_+^p \mid \lambda \text{ is measurable}\}.
    \end{aligned}
    $$
    
    Assume $(\hat w, \hat \lambda) \in \mathcal{P}_{w} \times \mathcal{P}_{\lambda}$ are two random variables depending on $\mathbf{\xi}$.
     In this section, $\mathbb{E}[\cdot]$ is the abbreviation of $\mathbb{E}_{\mathbf{\xi} \sim \mathcal{D}_{\xi}^s}[\cdot]$.
    The following theorem shows the equivalence between~\eqref{eq: BLO} and the penalty reformulation in the stochastic setting.

    \begin{theorem}
       \label{thm: equivalence of partial penalized stochastic single level}
       Suppose that Assumptions~\ref{ass: Lipschitz continuity},~\ref{ass: convexity},~\ref{ass: stochastic oracle} and~\ref{ass: LL Slater's condition} holds and $\gamma_1, \gamma_2 > 0$ are fixed parameters. 
    
    \noindent

    \noindent
    1. Assume $(x^*, y^*)$ is a global solution to~\eqref{eq: BLO}. 
    If $\mathcal{P}(\delta)$ defined in~\eqref{eq: P delta} is nonempty for any $(x, z) \in X \times Z$,
    then for any $(\hat w, \hat \lambda) \in \mathcal{P}(\delta)$,
    there exists $z^*$ such that 
    $(x^*, y^*, z^*)$ is a $\epsilon$-global-minima 
    of the following penalized form
    \begin{equation}
    \label{eq: stochastic partial penalized single level}
        \begin{aligned}
       \min_{(x, y, z)\in X \times Y \times Z} &\quad  \mathbb{E} \left[F(x, y)  + c_1 \left(G(x, y) - \ell_{\gamma}(x, z, \hat w(\mathbf{\xi}), \hat \lambda(\mathbf{\xi})) \right)\right], \\
       \text{s.t.}\quad \quad  & \quad G(x, y) - \mathbb{E} [\ell_{\gamma}(x, z, \hat w(\mathbf{\xi}), \hat \lambda(\mathbf{\xi})) ] \leq \epsilon_1, \quad \frac{1}{2} \sum_{i=1}^{p} [H_i(x, y)]_+^2 \leq \epsilon_2^2,
    \end{aligned}
    \end{equation}
    with any
    $c_1 \geq  \frac{L}{\mu_G} \epsilon^{-1}$, $\epsilon_2 \leq \frac{\epsilon}{4 \sqrt{2} c_1 B}$ and $\delta \leq \frac{\epsilon}{8c_1}$.
    

    \noindent 
    2.  
    By taking $c_1 = c_1^* + 2 := \frac{L}{\mu_G}\epsilon^{-1} + 2$, $\epsilon_2 \leq \frac{\epsilon}{4 \sqrt{2} c_1 B}$ and $\delta \leq \frac{\epsilon}{8c_1}$,
    for any $(\hat w, \hat \lambda) \in \mathcal{P}(\delta)$,
    the $\epsilon$-global-minima of~\eqref{eq: stochastic partial penalized single level} is a $\epsilon$-global-minima of the following approximation of BLO:
    \begin{equation}
    \label{eq: stochastic approximation of BLO 2}
        \begin{aligned}
       \min_{(x, y, z)\in X \times Y \times Z} &\quad  F(x, y) \\
       \text{s.t.}\quad \quad  & \quad G(x, y) - \mathbb{E} [\ell_{\gamma}(x, z, \hat w(\mathbf{\xi}), \hat \lambda(\mathbf{\xi})) ] \leq \epsilon_1, \quad \frac{1}{2} \sum_{i=1}^{p} [H_i(x, y)]_+^2 \leq \epsilon_2^2,
    \end{aligned}
    \end{equation}
    with some $\epsilon_1\leq \frac{17}{16} \epsilon$.
    
    \end{theorem}

    {\it Proof}
    \noindent
       1. By Lemma~\ref{lem: equivalence between the  value function-based reformulation and the dual augmented value function reformulation}, we have $E(x, z) \leq v(x), \forall z$. Then for $\forall (x, y, z) \in X \times Y \times Z, (\hat w, \hat \lambda) \in \mathcal{P}(\delta)$,  it holds that
       \begin{equation}
       \label{eq: proof of equivalence of partial penalized stochastic single level 3}
          \begin{aligned}
             &\quad \mathbb{E} \left[F(x, y)  + c_1 \left(G(x, y) - \ell_{\gamma}(x, z, \hat w(\mathbf{\xi}), \hat \lambda(\mathbf{\xi})) \right)\right] - F(x^*, y^*) \\
             &\geq \mathbb{E} \left[F(x, y)  + c_1 \left(G(x, y) - \ell_{\gamma}(x, z, \hat w(\mathbf{\xi}), \hat \lambda(\mathbf{\xi})) \right)\right] - F(x, y^*(x)) \\
             &=F(x, y) - F(x, y^*(x)) + c_1 (G(x, y) - E(x, z))
             - c_1 \left(\mathbb{E} [\ell_{\gamma}(x, z, \hat w(\mathbf{\xi}), \hat \lambda(\mathbf{\xi}))] - E(x, z)\right)  \\
             &\geq F(x, y) - F(x, y^*(x)) + c_1 (G(x, y) - v(x))
             - c_1 \delta  \quad \text{(by $E(x, z) \leq v(x)$ and~\eqref{eq: P delta})} \\
             &\geq -\frac{L_F}{2 c_1 \mu_G} - c_1  \sqrt{2} B \epsilon_2 - c_1 \delta \quad \text{(by~\eqref{eq: proof of equivalence of single level 1 1})} .
    \end{aligned}
       \end{equation}
    By~\eqref{eq: P delta}, it holds that
       \begin{equation}
          \label{eq: proof of equivalence of partial penalized stochastic single level 4}
       \begin{aligned}
          \mathbb{E} [G(x^*, y^*) - \ell_{\gamma}(x^*, z^*, \hat w(\mathbf{\xi}), \hat \lambda(\mathbf{\xi})) ] \geq -\delta.
       \end{aligned}
    \end{equation}
       Combining this with~\eqref{eq: proof of equivalence of partial penalized stochastic single level 3} gives
       $$
       \begin{aligned}
       &\quad \mathbb{E} \left[F(x, y)  + c_1 \left(G(x, y) - \ell_{\gamma}(x, z, \hat w(\mathbf{\xi}), \hat \lambda(\mathbf{\xi})) \right)\right] - \mathbb{E} \left[ F(x^*, y^*) - c_1 \left(G(x^*, y^*) - (x^*, z^*, \hat w(\mathbf{\xi}), \hat \lambda(\mathbf{\xi})) \right)\right] \\
       &\geq -\frac{L_F}{2 c_1 \mu_G} - c_1  \sqrt{2} B \epsilon_2 - 2 c_1 \delta.
       \end{aligned}
       $$
    
       By taking any $c_1 \geq  \frac{L}{\mu_G} \epsilon^{-1}$, $\epsilon_2 \leq \frac{\epsilon}{4 \sqrt{2} c_1 B}$ and $\delta \leq \frac{\epsilon}{8 c_1}$, the right side of the above inequality $\geq -\epsilon$.
    For the optimal solution $(x^*, y^*)$ of~\eqref{eq: BLO}
    and any $(\hat w, \hat \lambda) \in \mathcal{P}(\delta)$,
       by taking $z^* \in \arg \max_{z \in \mathbb{R}_+^p} \mathcal{L}_{\gamma_1}(x^*, y^*, z)$, we have $E(x^*, z^*) = v(x^*) = G(x^*, y^*)$. Lemma~\eqref{lem: equivalence between the  value function-based reformulation and the dual augmented value function reformulation}.3 implies that such a $z^* \in Z$ exists.
       Hence $(x^*, y^*, z^*)$ is a $\epsilon$-optimal solution of~\eqref{eq: stochastic partial penalized single level}.

    \noindent
    2.
    For any $\epsilon$-optimal solution $(\bar x, \bar y, \bar z)$ of~\eqref{eq: stochastic partial penalized single level}, it holds that
    $$
    \begin{aligned}
    &\quad \mathbb{E} \left[F(\bar x, \bar y)  + c_1 \left(G(\bar x, \bar y) - \ell_{\gamma}(\bar x, \bar z, \hat w(\mathbf{\xi}), \hat \lambda(\mathbf{\xi})) \right) \right] \\
    &\leq \mathbb{E} \left[F(x^*, y^*)  + c_1 \left(G(x^*, y^*) - \ell_{\gamma}(x^*, z^*, \hat w(\mathbf{\xi}), \hat \lambda(\mathbf{\xi})) \right) \right] + \epsilon \quad \text{(by $\epsilon$-optimality)} \\
    &\leq \mathbb{E} \left[F(x^*, y^*) \right] + c_1 \delta + \epsilon \quad \text{(by \eqref{eq: proof of equivalence of partial penalized stochastic single level 4})} \\
    &\leq \mathbb{E} \left[ F( \bar x, \bar y) + c_1^* 
    \left(G(\bar x, \bar y) - \ell_{\gamma}(\bar x, \bar z, \hat w(\mathbf{\xi}), \hat \lambda(\mathbf{\xi}))\right)\right]  + c_1\delta + 2 \epsilon ,
    \quad \text{{(by~\eqref{eq: proof of equivalence of partial penalized stochastic single level 3} with $c_1^*$ taken).}}
    \end{aligned}
    $$
    Then we have
    \begin{equation}
       G(\bar x, \bar y) - \mathbb{E} [\ell_{\gamma}(\bar x, \bar z, \hat w(\mathbf{\xi}), \hat \lambda(\mathbf{\xi}))]   \leq \frac{2\epsilon + c_1 \delta}{c_1 - c_1^*} \leq \frac{2\epsilon + \frac{\epsilon}{8}}{2} = 
       \frac{17}{16}\epsilon.
    \end{equation}
    Take $\epsilon_1 =G(\bar x, \bar y) - \mathbb{E} [\ell_{\gamma}(\bar x, \bar z, \hat w(\mathbf{\xi}), \hat \lambda(\mathbf{\xi}))]$.
    This implies $(\bar x, \bar y, \bar z)$ is feasible for~\eqref{eq: stochastic approximation of BLO 2}. For any feasible solution $(x, y, z)$ of~\eqref{eq: stochastic approximation of BLO 2}, the $\epsilon$-optimality of $(\bar x, \bar y, \bar z)$ implies
    \begin{equation}
       \begin{aligned}
          F(\bar x, \bar y) - F(x, y) 
          &\leq c_1 \left(G(x, y) - \mathbb{E} [\ell_{\gamma}(x, z, \hat w(\mathbf{\xi}), \hat \lambda(\mathbf{\xi}))]  \right) \\
          &\quad - c_1 \left( G(\bar x, \bar y) - \mathbb{E} [\ell_{\gamma}(\bar x, \bar z, \hat w(\mathbf{\xi}), \hat \lambda(\mathbf{\xi}))] 
        \right) + \epsilon \\
          &= c_1\left(G(x, y) - \mathbb{E} [\ell_{\gamma}(x, z, \hat w(\mathbf{\xi}), \hat \lambda(\mathbf{\xi}))] - \epsilon_1\right)  + \epsilon \leq \epsilon.
       \end{aligned}
    \end{equation}
    The last inequality follows from the feasibility of $(x, y)$.
    Therefore $(\bar x, \bar y, \bar z)$ is a $\epsilon$-optimal solution of~\eqref{eq: stochastic approximation of BLO 2}. 
    This completes the proof.
    \qed

    \begin{theorem}
       \label{thm: equivalence of stochastic single level}
       Suppose that Assumptions~\ref{ass: Lipschitz continuity},~\ref{ass: convexity},~\ref{ass: stochastic oracle} and~\ref{ass: LL Slater's condition} holds and $\gamma_1, \gamma_2 > 0$ are fixed parameters. 
    \noindent
    1. Assume $(x^*, y^*)$ is a global solution to~\eqref{eq: BLO}. 
    If $\mathcal{P}(\delta)$ defined in~\eqref{eq: P delta} is nonempty for any $(x, z) \in X \times Z$,
    then for any $(\hat w, \hat \lambda) \in \mathcal{P}(\delta)$,
    there exists $z^*$ such that 
    $(x^*, y^*, z^*)$ is a $\epsilon$-global-minima 
    of the following penalized form
    \begin{equation}
    \label{eq: stochastic penalized single level 2}
        \begin{aligned}
       \min_{(x, y, z)\in X \times Y \times Z} &\quad  \mathbb{E} [\Psi(x, y, z, \hat w, \hat \lambda
       ; \mathbf{\xi})] \\
       & \quad s:= \mathbb{E} \left[F(x, y)  + c_1 (G(x, y) - \ell_{\gamma}(x, z, \hat w(\mathbf{\xi}), \hat \lambda(\mathbf{\xi})) )
       + \frac{c_2}{2} \sum_{i=1}^{p} [H_i(x, y)]_+^2
       \right],
    \end{aligned}
    \end{equation}
    with any
    $c_1 \geq  \frac{2 L}{3 \mu_G} \epsilon^{-1}$, $c_2 \geq \frac{3}{2} (c_1)^2 B^2\epsilon^{-1}$ and $\delta \leq \frac{\epsilon}{6 c_1}$.
    
    \noindent 
    2.  
    By taking $c_1 = c_1^* + 2 :=  \frac{2 L}{3 \mu_G} \epsilon^{-1} + 2, 
    c_2 = c_2^* + 2 := \frac{3}{2} (c_1^*)^2 B^2\epsilon^{-1} + 2$ and
    $\delta \leq \frac{\epsilon}{6c_1}$,
    for any $(\hat w, \hat \lambda) \in \mathcal{P}(\delta)$,
    the $\epsilon$-global-minima of~\eqref{eq: stochastic penalized single level 2} is a $\epsilon$-global-minima of the~\eqref{eq: stochastic approximation of BLO 2}
    with some $\epsilon_1, \epsilon_2 \leq \frac{13}{12}\epsilon$.
    
    \end{theorem}

    {\it Proof}
       \noindent
          1. By Lemma~\ref{lem: equivalence between the  value function-based reformulation and the dual augmented value function reformulation}, we have $E(x, z) \leq v(x), \forall z$. Then for $\forall (x, y, z) \in X \times Y \times Z, (\hat w, \hat \lambda) \in \mathcal{P}(\delta)$,  it holds that
          \begin{equation}
          \label{eq: proof of equivalence of stochastic single level 3}
             \begin{aligned}
                &\quad \mathbb{E} \left[\Psi(x, y, z, \hat w, \hat \lambda; \mathbf{\xi}) \right] - F(x^*, y^*) \\
                &\geq \mathbb{E} \left[\Psi(x, y, z, \hat w, \hat \lambda; \mathbf{\xi})\right] - F(x, y^*(x)) \\
                &=F(x, y) - F(x, y^*(x)) + c_1 (G(x, y) - E(x, z))
                - c_1 \left(\mathbb{E} [\ell_{\gamma}(x, z, \hat w(\mathbf{\xi}), \hat \lambda(\mathbf{\xi}))] - E(x, z)\right) \\
                &\quad + \frac{c_2}{2} \sum_{i=1}^{p} [H_i(x, y)]_+^2 \\
                &\geq F(x, y) - F(x, y^*(x)) + c_1 (G(x, y) - E(x, z))
                - c_1 \delta +  \frac{c_2}{2} \sum_{i=1}^{p} [H_i(x, y)]_+^2  \quad \text{(by~\eqref{eq: P delta})} \\
                &\geq -\frac{L_F}{2 c_1 \mu_G} - \frac{c_1^2}{2 c_2} B^2 - c_1 \delta \quad \text{(by~\eqref{eq: proof of equivalence of single level 1})}.
       \end{aligned}
          \end{equation}
       By~\eqref{eq: P delta}, it holds that
          \begin{equation}
             \label{eq: proof of equivalence of stochastic single level 4}
          \begin{aligned}
             \mathbb{E} [G(x^*, y^*) - \ell_{\gamma}(x^*, z^*, \hat w(\mathbf{\xi}), \hat \lambda(\mathbf{\xi})) ] \geq -\delta.
          \end{aligned}
       \end{equation}
          Combining this with~\eqref{eq: proof of equivalence of stochastic single level 3} gives
          $$
          \begin{aligned}
          &\quad \mathbb{E} \left[\Psi(x, y, z, \hat w, \hat \lambda; \mathbf{\xi})\right] - \mathbb{E} \left[ \Psi(x^*, y^*, z^*, \hat w, \hat \lambda; \mathbf{\xi})\right] 
          \geq -\frac{L_F}{2 c_1 \mu_G} - \frac{c_1^2}{2c_2} B^2 - 2 c_1 \delta.
          \end{aligned}
          $$
       By taking $c_1 \geq  \frac{2 L}{3 \mu_G} \epsilon^{-1}$, $c_2 \geq \frac{3}{2} (c_1)^2 B^2\epsilon^{-1}$
          and $\delta \leq \frac{\epsilon}{6 c_1}$, the right side of the above inequality $\geq -\epsilon$.
       For the optimal solution $(x^*, y^*)$ of~\eqref{eq: BLO}
       and any $(\hat w, \hat \lambda) \in \mathcal{P}(\delta)$,
          by taking $z^* \in \arg \max_{z \in \mathbb{R}_+^p} \mathcal{L}_{\gamma_1}(x^*, y^*, z)$, we have $E(x^*, z^*) = v(x^*) = G(x^*, y^*)$. Lemma~\eqref{lem: equivalence between the  value function-based reformulation and the dual augmented value function reformulation}.3 implies that such a $z^* \in Z$ exists.
          Hence $(x^*, y^*, z^*)$ is a $\epsilon$-optimal solution of~\eqref{eq: stochastic partial penalized single level}.

       \noindent
       2.
       For any $\epsilon$-optimal solution $(\bar x, \bar y, \bar z)$ of~\eqref{eq: stochastic partial penalized single level}, it holds that
       $$
       \begin{aligned}
       &\quad \mathbb{E} \left[\Psi(\bar x, \bar y, \bar z, \hat w, \hat \lambda; \mathbf{\xi})  + c_1 \left(G(\bar x, \bar y) - \ell_{\gamma}(\bar x, \bar z, \hat w(\mathbf{\xi}), \hat \lambda(\mathbf{\xi})) \right) \right] \\
       &\leq \mathbb{E} \left[\Psi(x^*, y^*, z^*, \hat w, \hat \lambda; \mathbf{\xi}) \right] + \epsilon \quad \text{(by $\epsilon$-optimality)} \\
       &\leq \mathbb{E} \left[F(x^*, y^*) \right] + c_1 \delta + \epsilon \quad \text{(by \eqref{eq: proof of equivalence of stochastic single level 4})} \\
       &\leq \mathbb{E} \left[ F( \bar x, \bar y) + c_1^* 
       \left(G(\bar x, \bar y) - \ell_{\gamma}(\bar x, \bar z, \hat w(\mathbf{\xi}), \hat \lambda(\mathbf{\xi}))\right) + \frac{c_2^*}{2} \sum_{i=1}^{p}[H(\bar x, \bar y)]_+^2\right]  + c_1\delta + 2 \epsilon ,
       \quad \text{{(by~\eqref{eq: proof of equivalence of stochastic single level 3} with $c_1^*$ taken).}}
       \end{aligned}
       $$
       Then we have
       \begin{equation}
          G(\bar x, \bar y) - \mathbb{E} [\ell_{\gamma}(\bar x, \bar z, \hat w(\mathbf{\xi}), \hat \lambda(\mathbf{\xi}))] + \frac{1}{2} \sum_{i=1}^{p}[H(\bar x, \bar y)]_+^2   \leq \frac{2\epsilon + c_1 \delta}{2} \leq \frac{2\epsilon + \frac{\epsilon}{6}}{2} = 
          \frac{13}{12}\epsilon.
       \end{equation}
       Take $\epsilon_1 =G(\bar x, \bar y) - \mathbb{E} [\ell_{\gamma}(\bar x, \bar z, \hat w(\mathbf{\xi}), \hat \lambda(\mathbf{\xi}))]$, $\epsilon_2 = \frac{1}{2} \sum_{i=1}^{p}[H(\bar x, \bar y)]_+^2$.
       This implies $(\bar x, \bar y, \bar z)$ is feasible for~\eqref{eq: stochastic approximation of BLO 2}. For any feasible solution $(x, y, z)$ of~\eqref{eq: stochastic approximation of BLO 2}, the $\epsilon$-optimality of $(\bar x, \bar y, \bar z)$ implies
       \begin{equation}
          \begin{aligned}
             &\quad F(\bar x, \bar y) - F(x, y) \\
             &\leq c_1 \left(G(x, y) - \mathbb{E} [\ell_{\gamma}(x, z, \hat w(\mathbf{\xi}), \hat \lambda(\mathbf{\xi}))]  \right)  - c_1 \left( G(\bar x, \bar y) - \mathbb{E} [\ell_{\gamma}(\bar x, \bar z, \hat w(\mathbf{\xi}), \hat \lambda(\mathbf{\xi}))] 
           \right) 
           \\
           &\quad + \frac{c_2}{2} \sum_{i=1}^{p}[H_i(x, y)]_+^2 - \frac{c_2}{2} \sum_{i=1}^{p}[H_i(\bar x, \bar y)]_+^2
             + \epsilon \\
             &= c_1\left(G(x, y) - \mathbb{E} [\ell_{\gamma}(x, z, \hat w(\mathbf{\xi}), \hat \lambda(\mathbf{\xi}))] - \epsilon_1\right)  + c_2\left(\frac{1}{2} \sum_{i=1}^{p}[H_i(x, y)]_+^2 - \frac{1}{2} \sum_{i=1}^{p}[H_i(\bar x, \bar y)]_+^2\right)
           + \epsilon \\
             &\leq \epsilon.
          \end{aligned}
       \end{equation}
       The last inequality follows from the feasibility of $(x, y, z)$.
       Therefore $(\bar x, \bar y, \bar z)$ is a $\epsilon$-optimal solution of~\eqref{eq: stochastic approximation of BLO 2}. 
       This completes the proof.
       \qed

    \section{Analysis on the inner loop}
    
    In this section, we analyze the convergence of Algorithm~\ref{alg: ALM}. 
    Assume $(x, y, z)  = (x^{k-1}, y^{k-1}, z^{k-1})$ and ${\gamma_1}, \gamma_2$ are fixed.
    The expectation in this section denotes the expectation conditioned on $\tilde{\mathcal{F}}_{k-1}$, that is, $\mathbb{E} [\cdot] \triangleq \mathbb{E} [ \cdot | \tilde{\mathcal{F}}_{k-1}]$. Since $k$ is fixed, we abbreviate $w^{k, j}$ as $w^j$ and $\lambda^{k, j}$ as $\lambda^j$ for $j = 0, 1, ...$ in  this section.
    
    We write  $(w^*, \lambda^*) = (w^*(x, z), \lambda^*(x, z))$ as defined in~\eqref{eq: minimax subproblem}
    and $\ell_{\gamma}(x, z, w, \lambda)$ as  defined in~\eqref{eq: tilde L}.
    By Proposition~\ref{prop: strong duality}, the optimal solution $(w^*, \lambda^*)$ is a saddle point of $\ell_{\gamma}(x, z, w, \lambda)$, i.e.,
       \begin{equation}
          \label{eq: saddle point of Lagrangian}
       \ell_{\gamma}(x, z, w^*, \lambda) \leq \ell_{\gamma}(x, z,w^*, \lambda^*) \leq \ell_{\gamma}(x, z, w, \lambda^*), \quad \forall w \in Y, \lambda \in \mathbb{R}_+^p.
       \end{equation}
    
    To analyze the convergence of the inner loop, we establish the decreasing property for running a single step of the inner algorithm.  The following lemma gives the decreasing property of by taking gradient descent ascent steps. 
    \begin{lemma}
       \label{lem: inner algorithm decrease}
    The following relationships hold by taking primal step~\eqref{eq: inner algorithm primal step} and dual step~\eqref{eq: inner algorithm dual step}, respectively:
    \begin{equation}
    \label{eq: inner algorithm primal decrease}
    \begin{aligned}
         &\quad \mathbb{E}_{\xi_1, ..., \xi_j} [ \ell_{\gamma}(x, z, w^j, \lambda^j)] - \mathbb{E} [ \ell_{\gamma}(x, z, w, \lambda^j) ] + \frac{1}{2 \eta_j} \mathbb{E}_{\xi_1, ..., \xi_j} [ \|w^{j+1} - w\|^2] \\
          &\leq  \frac{1}{2} (\frac{1}{\eta_j} - \mu_G) \mathbb{E}_{\xi_1, ..., \xi_j} [ \|w^j - w\|^2] + \frac{\eta_j}{2} \mathbb{E}_{\xi_1, ..., \xi_j} [\|\nabla_w \ell_{\gamma}(x, z, w^j, \lambda^j; \xi_j)\|^2], \quad \forall w \in Y,
    \end{aligned}
    \end{equation}
    and
    \begin{equation}
    \label{eq: inner algorithm dual decrease}
    \begin{aligned}
          &\quad \ell_{\gamma}(x, z, w^j, \lambda) - \ell_{\gamma}(x, z, w^j, \lambda^j)
          +
          \frac{1}{2\rho_j}\| \lambda^{j+1} - \lambda\|^2
          \\
          &\leq \frac{1}{2}(\frac{1}{\rho_j} - \gamma_2) \| \lambda^{j} - \lambda \|^2 +  
          \frac{\rho_j}{2} \|\nabla_{\lambda} \ell_{\gamma}(x, z, w^j, \lambda^j)  \|^2, \quad \forall \lambda \in \mathbb{R}_+^p.
       \end{aligned}       
    \end{equation}
    Here the notation $\mathbb{E}_{\xi_1, ..., \xi_j}[\cdot]$ is the abbreviation of $\mathbb{E}_{\xi_1 \sim \mathcal{D}_{\xi}, ..., \xi_j \sim \mathcal{D}_{\xi}}[\cdot]$.
    \end{lemma}
    {\it Proof}
       The projection gradient step~\eqref{eq: inner algorithm primal step} gives
       \begin{equation}
          \label{eq: proof of inner algorithm decrease 1}
          \langle w^{j+1} - w, w^{j+1} - (w^j - \eta_j \nabla_y \mathcal{L}_{{\gamma_1}}(x, w^j, \lambda^j; \xi_j)) \rangle \leq 0.
       \end{equation}
       This implies
       \begin{equation}
          \label{eq: proof of inner algorithm decrease 2}
          \langle w^{j+1} - w, \nabla_y \mathcal{L}_{\gamma_1}(x, w^j, \lambda^j; \xi_j) \rangle \leq - \frac{1}{\eta_j} \langle w^{j+1} - w^j, w^{j+1} - w\rangle.
       \end{equation}
       By Young's inequality, we have
       \begin{equation}
          \label{eq: inner algorithm decrease 3}
          \begin{aligned}
          \langle w^{j+1} - w^j, \nabla_y \mathcal{L}_{\gamma_1}(x, w^j, \lambda^j; \xi_j) \rangle  
          &\geq - \frac{1}{2 \eta_j} \|w^{j+1} - w^j\|^2 - \frac{\eta_j}{2} \|\nabla_y \mathcal{L}_{\gamma_1}(x, w^j, \lambda^j; \xi_j)\|^2 .\\
          \end{aligned}
       \end{equation}
       By the strong convexity of $\mathcal{L}_{\gamma_1}(x, w, \lambda)$ with respect to $w$, we have
       \begin{equation}
          \label{eq: inner algorithm decrease 4}
          \langle w^j - w, \nabla_y \mathcal{L}_{\gamma_1}(x, w^j, \lambda^j) \rangle \geq 
          \mathcal{L}_{\gamma_1}(x, w^j, \lambda^j) - \mathcal{L}_{\gamma_1}(x, w, \lambda^j) + \frac{\mu_G}{2} \|w^j - w\|^2.
       \end{equation}
       Note that $w^j$ is independent of $\xi_j$, and hence 
       \begin{equation}
       \begin{aligned}
        &\quad  \mathbb{E}_{\xi_1, ..., \xi_j}[\langle w^j - w, \nabla_y \mathcal{L}_{\gamma_1}(x, w^j, \lambda^j; \xi_j)\rangle ]   \\
        &= \mathbb{E}_{\xi_1, ..., \xi_{j-1}} [\mathbb{E}_{\xi_j} [\langle w^j - w, \nabla_y \mathcal{L}_{\gamma_1}(x, w^j, \lambda^j; \xi_j)\rangle | \xi_1, ... ,\xi_{j-1}]] \\
        &= \mathbb{E}_{\xi_1, ..., \xi_j}[ \langle w^j - w, \nabla_y \mathcal{L}_{\gamma_1}(x, w^j, \lambda^j)\rangle ] .\\
        \end{aligned}
        \end{equation}
       Summing up ~\eqref{eq: inner algorithm decrease 3},\eqref{eq: inner algorithm decrease 4} and then taking expectation gives
       \begin{equation*}
          \label{eq: inner algorithm decrease 5}
          \begin{aligned}
          &\quad \mathbb{E}_{\xi_1, ..., \xi_j} [\langle w^{j+1} - w, \nabla_y \mathcal{L}_{\gamma_1}(x, w^j, \lambda^j; \xi_j) \rangle ] \\
          &\geq - \frac{1}{2 \eta_j} \mathbb{E}_{\xi_1, ..., \xi_j} [\|w^{j+1} - w^j\|^2 ]- \frac{\eta_j}{2} \mathbb{E}_{\xi_1, ..., \xi_j} [\|\nabla_y \mathcal{L}_{\gamma_1}(x, w^j, \lambda^j; \xi_j)\|^2 ]\\
          & +\mathbb{E}_{\xi_1, ..., \xi_j} [ \mathcal{L}_{\gamma_1}(x, w^j, \lambda^j) ]- \mathbb{E}_{\xi_1, ..., \xi_j} [\mathcal{L}_{\gamma_1}(x, w, \lambda^j)]  + \frac{\mu_G}{2} \mathbb{E}_{\xi_1, ..., \xi_j} [\|w^j - w\|^2] \\
           &= - \frac{1}{2 \eta_j} \mathbb{E}_{\xi_1, ..., \xi_j} [ \|w^{j+1} - w^j\|^2 ]- \frac{\eta_j}{2} \mathbb{E}_{\xi_1, ..., \xi_j} [ \|\nabla_w \ell_{\gamma}(x, z, w^j, \lambda^j; \xi_j)\|^2] \\
          & + \mathbb{E}_{\xi_1, ..., \xi_j} [ \ell_{\gamma}(x, z, w^j, \lambda^j) ]- \mathbb{E}_{\xi_1, ..., \xi_j} [ \ell_{\gamma}(x, z, w, \lambda^j) ] + \frac{\mu_G}{2} \mathbb{E}_{\xi_1, ..., \xi_j} [ \|w^j - w\|^2 ].
          \end{aligned}
       \end{equation*}
       
       Utilizing the identity 
       $
       \|w^{j+1} - w^j\|^2 + \|w^{j+1} - w\|^2 - \|w^j - w\|^2 - 2 \langle w^{j+1} - w, w^{j+1} - w^j\rangle = 0
       $
       and~\eqref{eq: proof of inner algorithm decrease 2},
       we have
       \begin{equation*}
       \begin{aligned}
          &\quad \mathbb{E}_{\xi_1, ..., \xi_j} [ \ell_{\gamma}(x, z, w^j, \lambda^j)] - \mathbb{E}_{\xi_1, ..., \xi_j} [ \ell_{\gamma}(x, z, w, \lambda^j) ] \\
           &\leq \mathbb{E}_{\xi_1, ..., \xi_j} [ \langle w^{j+1} - w, \nabla_y \mathcal{L}_{\gamma_1}(x, w^j, \lambda^j; \xi_j) \rangle ] + \frac{1}{2 \eta_j} \mathbb{E}_{\xi_1, ..., \xi_j} [ \|w^{j+1} - w^j\|^2 ] \\
           &\quad + \frac{\eta_j}{2} \mathbb{E}_{\xi_1, ..., \xi_j} [ \|\nabla_w \ell_{\gamma}(x, z, w^j, \lambda^j; \xi_j)\|^2]
           - \frac{\mu_G}{2}\mathbb{E}_{\xi_1, ..., \xi_j} [ \|w^j - w\|^2] \\
             &\leq - \frac{1}{ \eta_j} \mathbb{E}_{\xi_1, ..., \xi_j} [ \langle w^{j+1} - w^j, w^{j+1} - w\rangle ] + \frac{1}{2\eta_j} \mathbb{E}_{\xi_1, ..., \xi_j} [ \|w^{j+1} - w^j\|^2 ] \\
             &\quad - \frac{\mu_G}{2}\mathbb{E}_{\xi_1, ..., \xi_j} [ \|w^j - w\|^2] 
              + \frac{\eta_j}{2} \mathbb{E}_{\xi_1, ..., \xi_j} [ \|\nabla_w \ell_{\gamma}(x, z, w^j, \lambda^j; \xi_j)\|^2] \\
             &=\frac{1}{2} (\frac{1}{\eta_j} - \mu_G) \mathbb{E}_{\xi_1, ..., \xi_j} [ \|w^j - w\|^2] -\frac{1}{2 \eta_j} \mathbb{E}_{\xi_1, ..., \xi_j} [ \|w^{j+1} - w\|^2] \\
             &\quad + \frac{\eta_j}{2} \mathbb{E}_{\xi_1, ..., \xi_j} [\|\nabla_w \ell_{\gamma}(x, z, w^j, \lambda^j; \xi_j)\|^2].
       \end{aligned}
       \end{equation*}
       This gives~\eqref{eq: inner algorithm primal decrease}. 
       $\mathcal{L}(x, w, \lambda)$ is concave in $\lambda$ implies 
       $\ell_{\gamma}(x, z, w, \lambda)$ is $\gamma_2$-strongly concave in $\lambda$.
       Similarly, we can obtain~\eqref{eq: inner algorithm dual decrease} by taking the dual step~\eqref{eq: inner algorithm dual step}. This completes the proof.
    \qed

       Combining the primal and dual steps,
       the following lemma shows the  decreasing property with respect to the duality gap.
       \begin{corollary} For any $(w, \lambda) \in Y \times \mathbb{R}_+^p$ it holds that 
       \begin{equation}
             \label{eq: inner algorithm decrease corollary}
             \begin{aligned}
             &\quad \mathbb{E}_{\xi_1, ..., \xi_j} [\ell_{\gamma}(x, z, w^j, \lambda) - \ell_{\gamma}(x, z, w, \lambda^j)] + \frac{1}{2 \eta_j} \mathbb{E}_{\xi_1, ..., \xi_j} [\|w^{j+1} - w\|^2] \\
             &\quad + \frac{1}{2 \rho_j} \mathbb{E}_{\xi_1, ..., \xi_j} [\| \lambda^{j+1} - \lambda\|^2] \\
             &\leq  \frac{1}{2} (\frac{1}{\eta_j} - \mu_G)\mathbb{E}_{\xi_1, ..., \xi_j}[\|w^j - w\|^2] + \frac{1}{2}(\frac{1}{\rho_j} - \gamma_2) \mathbb{E}_{\xi_1, ..., \xi_j}[\| \lambda^j - \lambda\|^2] \\
             &\quad + \frac{\eta_j}{2} (M_{\mathcal{L}, 1} + M_{\mathcal{L}, 2} \mathbb{E}_{\xi_1, ..., \xi_j} [\|\lambda^j\|^2]) + \frac{\rho_j}{2} \left(4 \gamma_2 \|z \|^2 + 4 M_{H, 0}^2 +  2(\gamma_1 + \gamma_2) \mathbb{E}_{\xi_1, ..., \xi_j}[\| \lambda^j\|^2 ]\right).
             \end{aligned}      
          \end{equation}
       \end{corollary}
       {\it Proof}
          From~\eqref{eq: derivative of tilde L}, we have
          \begin{equation}
             \label{eq: proof of inner algorithm decrease corollary 1}
             \begin{aligned}
             \| \nabla_{\lambda} \ell_{\gamma}(x, z, w^j, \lambda^j) \|^2
             &\leq 
             \| \max(\gamma_2 z - (\gamma_1 + \gamma_2)\lambda^j, \gamma_2 z - \gamma_2 \lambda^j + H(x, w^j)) \|^2 \\ 
             &\leq 4 \gamma_2 \|z \|^2 + 4 M_{H, 0}^2 + 2 (\gamma_1 + \gamma_2) \| \lambda^j\|^2.
             \end{aligned}
          \end{equation}
          Then combining~\eqref{eq: bound of first order oracle of Lk 2},~\eqref{eq: inner algorithm primal decrease},~\eqref{eq: inner algorithm dual decrease} and~\eqref{eq: proof of inner algorithm decrease corollary 1} gives~\eqref{eq: inner algorithm decrease corollary}.
       \qed
       
       By taking decreasing step sizes, we can prove that the dual variable is bounded during the iterations of the inner loop in  expectation.
    
       \begin{lemma}
          \label{lem: bound of dual variable in inner loop}
          Choose the dual step size as
          $
          \rho_j = \frac{\rho}{j+1},
          $
          where the constants $\rho > \frac{1}{\gamma_2}$, then the following bounds hold:
          \begin{equation}
          \label{eq: bound of dual variable in inner loop}
          \| \lambda^j \|^2 \leq M_{\lambda}:=2 \rho^2 \gamma_2^2 \|z\|^2 + 2 p \rho^2 M_{H, 0}^2
          , \quad j=1,..., s .
          \end{equation}
       \end{lemma}
       
       {\it Proof}
          From~\eqref{eq: inner algorithm dual step} and~\eqref{eq: tilde L} we have
          \begin{equation}
             \begin{aligned}
             \lambda^{j+1} &= \mathrm{Proj}_{\mathbb{R}_+^p} \left(\lambda^j + \frac{\rho}{j+1} ( \nabla_z \mathcal{L}_{\gamma_1}(x, w^j, \lambda^j) - \gamma_2 (\lambda^j - z) ) \right)  \\
             &= \mathrm{Proj}_{\mathbb{R}_+^p} \left(\lambda^j + \frac{\rho}{j+1} ( \max(- \gamma_1 \lambda^j, H(x, w^j)) - \gamma_2 \lambda^j + \gamma_2 z ) \right) 
             \end{aligned}
          \end{equation}
          Consider the $i$-th component of $\lambda^{j+1}$, we have
          $$ 
          \begin{aligned}
          \lambda_i^{j+1} 
          &= \max\left(0, \lambda_i^j + \frac{\rho}{j+1} ( - \gamma_1 \lambda_i^j + \gamma_2 z_i + \max(0, H_i(x, w^j)) )\right) \\
          &
          \leq \left| \lambda^j + \frac{\rho}{j+1} ( - \gamma_2 \lambda^j) \right| + \frac{\rho}{j+1}(\gamma_2 |z_i| + | \max(- \gamma_1 \lambda^j_i, H_i(x, w^j)) |
          \end{aligned}
          $$
          Since $\lambda^j_i \geq 0$, it holds that $\max(- \gamma_1 \lambda_i^j, H_i(x, w^j)) \leq \max(0, H_i(x, w^j)) \leq M_{H, 0}$. Then we obtain a recursive relation as
          $$
          \lambda^{j+1}_i \leq ( 1 - \frac{\rho \gamma_2}{j+1}) \lambda^j_i + \frac{\rho}{j+1} \gamma_2 |z_i| + \frac{\rho}{j+1} M_{H, 0}.
          $$
          Multiplying both sides by ${j+1}$ and using $\rho > \frac{1}{\gamma_2}$, we have
          $$
          (j+1) \lambda^{j+1}_i \leq j \lambda^j_i + \rho \gamma_2 |z_i| + \rho M_{H, 0}.
          $$
          Note that the initial point $\lambda^0_i = 0$. Hence $\lambda^{j}_i \leq  \rho \gamma_2 |z_i| + \rho M_{H, 0}), \forall j, i$. 
           This gives the bound of $\lambda^j$ as~\eqref{eq: bound of dual variable in inner loop}.
       \qed
       
       Now we are ready to establish the convergence of the inner loop. In the following Theorem, we show that the output pair $( w^s, \lambda^s)$ is an $\widetilde{O}(\frac{1}{s_k})$-optimal point of the lower-level problem measured by the square norm.
       \begin{theorem}
          \label{thm: point convergence of inner loop}
          By taking the step size as 
          \begin{equation}
             \label{eq: step size of inner loop}
          \eta_j = \frac{\eta}{j+1}, \quad \rho_j = \frac{\rho}{j+1},
          \end{equation}
          where the constants $\eta \geq \frac{1}{\mu_G}, \rho \geq \frac{1}{\gamma_2}$, there exist constants $\phi_1, \phi_2 > 0$ such that
          \begin{equation}
             \label{eq: point convergence of inner loop}
             \begin{aligned}
              \mathbb{E}_{\mathbf{\xi}} [ \|w^s - w^*\|^2] &\leq {\phi}_1 \frac{1+\log(s)}{s} ,\\
              \mathbb{E}_{\mathbf{\xi}}[ \|\lambda^s - \lambda^*\|^2 ]&\leq {\phi}_2 \frac{1+\log(s)}{s} ,
             \end{aligned}
          \end{equation}
       where $\phi_1 = \eta\left( \eta M_{\mathcal{L}, 1} + \rho (4 \gamma_2 \|z \|^2 + 4 M_{H, 0}^2 ) + (\eta M_{\mathcal{L}, 2} + 2 \rho (\gamma_1 + \gamma_2)) M_{\lambda} \right)$, ${\phi}_2 = \frac{\rho}{\eta} \phi_1$.
       \end{theorem}
       {\it Proof}
          Taking $\eta_j = \frac{\eta}{j+1}$, $\rho_j = \frac{\rho}{j+1}$ in~\eqref{eq: inner algorithm decrease corollary} gives
          \begin{equation}
             \label{eq: proof of bound of dual variable in inner loop 1}
             \begin{aligned}
             &\quad   \mathbb{E}_{\xi_1, ..., \xi_j} [\ell_{\gamma}(x, z, w^j, \lambda) - \ell_{\gamma}(x, z, w, \lambda^j)]
              + \frac{j+1}{2\eta} \mathbb{E}_{\xi_1, ..., \xi_j} [\|w^{j+1} - w^*\|^2] \\
              &\quad + \frac{j+1}{2 \rho} \mathbb{E}_{\xi_1, ..., \xi_j} [\| \lambda^{j+1} - \lambda\|^2] \\
             &\leq \frac{1}{2} (\frac{j+1}{\eta} - \mu_G)\mathbb{E}_{\xi_1, ..., \xi_j} [\|w^j - w^*\|^2] + \frac{1}{2}(\frac{j+1}{\rho} - \gamma_2) \mathbb{E}_{\xi_1, ..., \xi_j} [\| \lambda^j - \lambda\|^2 ] 
              \\ &\quad + \frac{\eta}{2(j+1)} (M_{\mathcal{L}, 1} + M_{\mathcal{L}, 2} \mathbb{E}_{\xi_1, ..., \xi_j}[\|\lambda^j\|^2]) 
               \\
               &\quad + \frac{\rho}{2(j+1)} \left(4 \gamma_2 \|z \|^2 + 4 M_{H, 0}^2 +  2(\gamma_1 + \gamma_2) \mathbb{E}_{\xi_1, ..., \xi_j}[\| \lambda^j\|^2 ]\right)  \\
              &\leq \frac{j}{2\eta}\mathbb{E}_{\xi_1, ..., \xi_j} [\|w^j - w^*\|^2] + \frac{j}{2\eta} \mathbb{E}_{\xi_1, ..., \xi_j} [\| \lambda^j - \lambda\|^2 ] 
              \\
              &\quad + \frac{\eta}{2(j+1)} (M_{\mathcal{L}, 1} + M_{\mathcal{L}, 2} \mathbb{E}_{\xi_1, ..., \xi_j}[\|\lambda^j\|^2]) \\
              &\quad + \frac{\rho}{2(j+1)} \left(4 \gamma_2 \|z \|^2 + 4 M_{H, 0}^2 +  2(\gamma_1 + \gamma_2) \mathbb{E}_{\xi_1, ..., \xi_j}[\| \lambda^j\|^2 ]\right).  
             \end{aligned}
          \end{equation}
          Then taking $w = w^*, \lambda = \lambda^*$ and utilizing $\ell_{\gamma}(x, z, w^j, \lambda^*) \geq \ell_{\gamma}(x, z, w^*, \lambda^*) \geq \ell_{\gamma}(x, z, w^*, \lambda^j)$ gives
          \begin{equation}
             \begin{aligned}
             &\quad \frac{j+1}{2\eta} \mathbb{E}_{\xi_1, ..., \xi_j} [\|w^{j+1} - w^*\|^2] + \frac{j+1}{2 \rho} \mathbb{E}_{\xi_1, ..., \xi_j} [\| \lambda^{j+1} - \lambda^*\|^2] \\
             &\leq \frac{j}{2\eta} \mathbb{E}_{\xi_1, ..., \xi_j} [\|w^{j} - w^*\|^2] + \frac{j}{2 \rho} \mathbb{E}_{\xi_1, ..., \xi_j} [\| \lambda^{j} - \lambda^*\|^2] \\
             &\quad + \frac{\eta}{2(j+1)} (M_{\mathcal{L}, 1} + M_{\mathcal{L}, 2} \mathbb{E}_{\xi_1, ..., \xi_j}[\|\lambda^j\|^2])\\
             &\quad  + \frac{\rho}{2(j+1)} \left(4 \gamma_2 \|z \|^2 + 4 M_{H, 0}^2 +  2(\gamma_1 + \gamma_2) \mathbb{E}_{\xi_1, ..., \xi_j}[\| \lambda^j\|^2 ]\right) .
             \end{aligned}
          \end{equation}
          This recursive relationship gives 
          \begin{equation}
             \label{eq: proof of bound of dual variable in inner loop 2}
             \begin{aligned}
             &\quad \frac{j}{2\eta} \mathbb{E}_{\xi_1, ..., \xi_j} [\|w^{j} - w^*\|^2] + \frac{j}{2 \rho} \mathbb{E}_{\xi_1, ..., \xi_j} [\| \lambda^{j} - \lambda^*\|^2] \\
             &\leq
             \sum_{l=0}^{j-1} \frac{\eta}{2(l+1)} (M_{\mathcal{L}, 1} + M_{\mathcal{L}, 2} \mathbb{E}_{\xi_1, ..., \xi_j}[\|\lambda^{l}\|^2]) + \sum_{l=0}^{j-1} \frac{\rho}{2(l+1)} \left(4 \gamma_2 \|z \|^2 + 4 M_{H, 0}^2 \right) \\
             &\quad +  2(\gamma_1 + \gamma_2) \mathbb{E}_{\xi_1, ..., \xi_j}[\| \lambda^l\|^2 ] \\
             &\leq \frac{1}{2}(1 + \log(j)) \left( \eta M_{\mathcal{L}, 1} + \rho (4 \gamma_2 \|z \|^2 + 4 M_{H, 0}^2 ) \right) \\
             &\quad + \frac{1}{2}(\eta M_{\mathcal{L}, 2} +  \rho (\gamma_1 + \gamma_2)) \sum_{l=0}^{j-1} \frac{1}{l+1}\mathbb{E}_{\xi_1, ..., \xi_j}[\|\lambda^{l}\|^2]  .
             \end{aligned}
          \end{equation}

       Substituting~\eqref{eq: bound of dual variable in inner loop} into 
        ~\eqref{eq: proof of bound of dual variable in inner loop 2} 
         and taking $j=s$ we have
          \begin{equation}
            \label{eq: proof of point convergence of inner loop 1}
            \begin{aligned}
            &\quad \frac{s}{2\eta} \mathbb{E}_{\mathbf{\xi}} [\|w^{s} - w^*\|^2] + \frac{s}{2 \rho} \mathbb{E}_{\mathbf{\xi}} [\| \lambda^{j} - \lambda^*\|^2]  \\
             &\leq \frac{1}{2}(1+\log(s)) \left( \eta M_{\mathcal{L}, 1} + \rho (4 \gamma_2 \|z \|^2 + 4 M_{H, 0}^2 ) + (\eta M_{\mathcal{L}, 2} + 2 \rho (\gamma_1 + \gamma_2)) M_{\lambda} \right) .
            \end{aligned}
         \end{equation}
          This gives the boundness of $\mathbb{E}_{\mathbf{\xi}} [\|w^s - w^*\|^2 ]$ and $\mathbb{E}_{\mathbf{\xi}}[ \|\lambda^s - \lambda^*\|^2]$ in~\eqref{eq: point convergence of inner loop}.     
       \qed
    
       By utilizing the measurement as objective function, a similar convergence rate can be obtained as follows.
       \begin{corollary}
          Under the same conditions as in Theorem~\ref{thm: point convergence of inner loop}, it holds that
          \begin{equation}
             \label{eq: objective convergence of inner loop}
             \left| \mathbb{E}_{\mathbf{\xi}} [\ell_{\gamma}(x, z, w^s, \lambda^s)] - E(x, z)\right| \leq \widetilde{O}\left(\frac{1}{s}\right).
          \end{equation}
       \end{corollary}
       {\it Proof}
       From~\eqref{eq: second order derivative of L},~\eqref{eq: derivative of tilde L}, we 
        see $\nabla_w \ell{\gamma}(x, z, w, \lambda) $ is Lipschitz continuous in $w$ with module $L_{\ell, w} = L_{G} + p(\sqrt{M_{\lambda}} + \gamma_1 ^{-1}M_{H, 0})L_{H} + M_{H, 1}^2$, $\nabla_{\lambda} \mathcal{L}(x, z, w, \lambda) $ is Lipschitz continuous in $\lambda$ with module $L_{\ell, \lambda} = 
        \gamma_1 + \gamma_2
        $ given the bound $\|\lambda\|^2 \leq M_{\lambda}$.
       Then under the strongly-convex-strongly-concave condition shown in Lemma~\ref{lem: convex-concav of tilde L}
       , the objective gap
    $$
        \mathbb{E}_{\mathbf{\xi}} [\ell_{\gamma}(x, z, w, \lambda^*)] - \mathbb{E}_{\mathbf{\xi}} [\ell_{\gamma}(x, z, w^*, \lambda)]
        =
        \mathbb{E}_{\mathbf{\xi}} [\ell_{\gamma}(x, z, w, \lambda^*) - \ell_{\gamma}(x, z, w^*, \lambda^*)] + \mathbb{E}_{\mathbf{\xi}} [\ell_{\gamma}(x, z, w, \lambda^*) - \ell_{\gamma}(x, z, w^*, \lambda)] 
        $$
        has a relationship to $\mathbb{E}_{\mathbf{\xi}} [\|w - w^*\|^2] +  \mathbb{E}_{\mathbf{\xi}} [\| \lambda - \lambda^*\|^2]$ as follows:
       $$
       \begin{aligned}
     \frac{\mu_G}{2} \mathbb{E}_{\mathbf{\xi}} [\|w - w^*\|^2] + \frac{\gamma_2}{2}\mathbb{E}_{\mathbf{\xi}} [\| \lambda - \lambda^*\|^2]
    &\leq
    \left| \mathbb{E}_{\mathbf{\xi}} [\ell_{\gamma}(x, z, w, \lambda^*)] - \mathbb{E}_{\mathbf{\xi}} [\ell_{\gamma}(x, z, w^*, \lambda)] \right|\\
    &\leq \frac{L_{\ell, w}}{2} \mathbb{E}_{\mathbf{\xi}} [\|w - w^*\|^2] + \frac{L_{\ell, \lambda} }{2}\mathbb{E}_{\mathbf{\xi}} [\| \lambda^{j} - \lambda^*\|^2].
    \end{aligned}
       $$
       Therefore, the convergence rate of the objective function is also $\widetilde{O}(\frac{1}{s})$, that is
       $$
       \left| \mathbb{E}_{\mathbf{\xi}} \left[\ell_{\gamma}(x, z, w^s, \lambda^*)\right] - \mathbb{E}_{\mathbf{\xi}} \left[\ell_{\gamma}(x, z, w^*, \lambda^s)\right] \right| \leq \widetilde{O}\left(\frac{1}{s}\right).
       $$
       Note that $\ell_{\gamma}(x, z, w^s, \lambda^*) \geq \ell_{\gamma}(x, z, w^s, \lambda^s)$ and $\ell_{\gamma}(x, z, w^*, \lambda^s) \leq \ell_{\gamma}(x, z, w^*, \lambda^*) = E(x, z)$, we have
       $$
       \mathbb{E}_{\mathbf{\xi}} [\ell_{\gamma}(x, z, w^s, \lambda^s)] - E(x, z) \leq \mathbb{E}_{\mathbf{\xi}} [\ell_{\gamma}(x, z, w^s, \lambda^*)] - \mathbb{E}_{\mathbf{\xi}} [\ell_{\gamma}(x, z, w^*, \lambda^s)]\leq \widetilde{O}\left(\frac{1}{s}\right).
       $$
       Similarly, it holds that 
       $$
       E(x, z) - \mathbb{E}_{\mathbf{\xi}} [\ell_{\gamma}(x, z, w^s, \lambda^s)] \leq \widetilde{O}\left(\frac{1}{s}\right).
       $$
       This completes the proof.
    \qed

    \section{Analysis on the outer loop}
    In the outer loop, we apply the stochastic gradient descent (SGD) method to solve the saddle point problem~\eqref{eq: maintext stochastic penalized single level}. Unlike the standard analysis of SGD, this setting presents two main challenges: first, constructing a suitable stochastic gradient oracle for~\eqref{eq: stochastic gradient of Psi}; second, addressing the bias introduced by the inexact solution of the lower-level problem.
    
    In the following proposition, we show that the gradient of $\mathcal{G}$ is continuously differentiable and its gradient is computable given the optimal solution of the subproblem~\eqref{eq: minimax subproblem}. 
    \begin{lemma}
       \label{lem: gradient of mathcal v}
    Assume $(w^*, \lambda^*) = (w^*(x, z), \lambda^*(x, z))$ is defined in~\eqref{eq: minimax subproblem}. 
    Then $\mathcal{G}(x, y, z)$ defined in~\eqref{eq: single level} is continuously differentiable and its gradient is given by
       \begin{subequations}
          \begin{align}
             \nabla_x \mathcal{G}(x, y, z) &= \nabla_x g(x, y) - \nabla_x \ell_{\gamma}(x, z, w^*, \lambda^*), \label{eq: gradient of mathcal G 1} \\
             \nabla_y \mathcal{G}(x, y, z) &= \nabla_y g(x, y), 
             \label{eq: gradient of mathcal G 2} \\
             \nabla_z \mathcal{G}(x, y, z) &= \gamma_2 (\lambda^* - z) \label{eq: gradient of mathcal G 3}.
          \end{align}
       \end{subequations}
    Furthermore substituting the fact $\|z\| \leq B$ we obtain that $\mathcal{G}(x, y, z)$ is Lipschitz continuous with module $L_{\mathcal{G}} = 3M_{G, 1} + \sqrt{p}  M_{H, 1} (B + \frac{\sqrt{p}}{\gamma_2} M_{H, 0} )+  p M_{H, 0}M_{H, 1} + \sqrt{p} M_{H, 0} = L_{\mathcal{G}}$.
    \end{lemma}
    {\it Proof}
       By Theorem 4.24 in~\cite{bonnans2013perturbation}, $E(x, z)$ is continuously differentiable with respect to $z$ and its gradient is given by 
       $$
       \begin{aligned}
       \nabla_x E(x, z) &= \nabla_x \left( \max_{\lambda \in \mathbb{R}_+^p} \left\{ D(x, \lambda) + \frac{\gamma_2}{2} \| \lambda - z\|^2 \right\} \right) \\
       &= \nabla_x \left( D(x, \lambda^*) + \gamma_2 (\lambda^*(x, z) - z) \right). \quad \text{(due to the uniqueness of $\lambda^*$)} \\
       &= \nabla_x \ell_{\gamma}(x, z, w^*, \lambda^*) .
       \end{aligned}
       $$ 
       Hence $\mathcal{G}(x, y, z)$ is continuously differentiable and~\eqref{eq: gradient of mathcal G 1} holds.
       Since $E(x, z)$ is independent of $y$,~\eqref{eq: gradient of mathcal G 2} is straightforward.
       Note that  $- E(x, z)$ is the Moreau envelope of $- D(x, z)$ for any $x$, Proposition~\ref{prop: properties of Moreau envelope}.4 gives $- \nabla_z E(x, z) = \gamma_2 (z - \mathrm{prox}_{- \frac{1}{\gamma_2} D(x, \cdot)}(z))$. 
       On the other hand, the optimality condition of the subproblem~\eqref{eq: minimax subproblem} gives 
       $$
       \mathrm{prox}_{- \frac{1}{\gamma_2} D(x, \cdot)}(z) = \arg \max_{\lambda \in \mathbb{R}_+^p} \left\{ D(x, \lambda) + \frac{\gamma_2}{2} \| \lambda - z\|^2 \right\} = \arg \max_{\lambda \in \mathbb{R}_+^p} \left\{ \min_{w \in Y} \ell_{\gamma}(x, z, w, \lambda) \right\} = \lambda^* .
       $$ 
       This gives~\eqref{eq: gradient of mathcal G 3}. 
       
       Now we show the Lipschitz continuity of $\mathcal{G}(x, y, z)$, that is $\| \nabla \mathcal{G}(x, y, z) \| \leq L_{\mathcal{G}}$.
       Similar to the proof of Lemma~\ref{lem: bound of first order oracle of L with respect to y}, it holds that $\| \nabla_x \ell_{\gamma}(x, z, w^*, \lambda^*) \| \leq M_{G, 1} + \frac{1}{{\gamma_1}} \sum_{i=1}^{p}  ( {\gamma_1} | \lambda^*_i | +  M_{H, 0})M_{H, 1}$.
       Again by Theorem 4.24 in~\cite{bonnans2013perturbation}, $D(x, \lambda) = \max_{y \in Y} \{ \mathcal{L}(x, y, \lambda) \}$ is continuously differentiable with respect to $\lambda$ and its gradient is given by $\nabla_\lambda D(x, \lambda) = \nabla_z \mathcal{L}(x, y^*(x, \lambda), \lambda)$, where $y^*(x, \lambda)$ is the optimal solution of the subproblem $\min_{y \in Y} \mathcal{L}(x, y, \lambda)$. 
       From~\eqref{eq: derivative of L} we know $\| \nabla_z \mathcal{L}(x, y, \lambda) \| \leq \sqrt{p} M_{H, 0}$. Hence $D(x, \lambda)$ is $\sqrt{p} M_{H, 0}-$Lipschitz continuous with respect to $\lambda$. Proposition~\ref{prop: properties of Moreau envelope}.4 implies the Moreau envelope of $E(x, z)$ is also Lipschitz continuous with the same Lipschitz constant. This gives 
       \begin{equation}
       \label{eq: proof of gradient of mathcal G 1}
       \| \nabla_z \mathcal{G}(x, y, z) \| = \gamma_2 \|\lambda^* - z\| \leq \sqrt{p} M_{H, 0}.
       \end{equation}
    Therefore,
       $$
       \begin{aligned}
       \| \nabla \mathcal{G}(x, y, z) \| 
       &\leq \| \nabla_x G(x, y) \| + \| \nabla_x \ell_{\gamma}(x, z, w^*, \lambda^*) \| + \| \nabla_y G(x, y) \| + \| \nabla_z \mathcal{G}(x, y, z) \|\\
       & \leq M_{G,1} + M_{G, 1} + \frac{1}{{\gamma_1}} \sum_{i=1}^{p}  ( {\gamma_1} | \lambda^*_i | +  M_{H, 0})M_{H, 1} +M_{G,1} + \sqrt{p} M_{H, 0}  \\
       &\leq 3M_{G, 1} + \sqrt{p}  M_{H, 1} \| \lambda^* \| +  p M_{H, 0}M_{H, 1} + \sqrt{p} M_{H, 0} \\
       &\leq 3M_{G, 1} + \sqrt{p}  M_{H, 1} (\|z \| + \frac{\sqrt{p}}{\gamma_2} M_{H, 0} )+  p M_{H, 0}M_{H, 1} + \sqrt{p} M_{H, 0}  \quad \text{(by~\eqref{eq: proof of gradient of mathcal G 1})} \\
       &\leq 3M_{G, 1} + \sqrt{p}  M_{H, 1} (B + \frac{\sqrt{p}}{\gamma_2} M_{H, 0} )+  p M_{H, 0}M_{H, 1} + \sqrt{p} M_{H, 0} = L_{\mathcal{G}} .
       \end{aligned}
       $$
       The last inequality follows from that $\|z\| \leq B$).
       This completes the proof.
    \qed
    
    \begin{corollary}
    \label{cor: Lipschitz of Phi}
    $\nabla \Psi(\mathbf{u})$ is Lipschitz continuous with modulus
    \begin{equation}
       \label{eq: Lipschitz continuity of Psi}
       L_{\Psi}  = L_F + c_1 L_{\mathcal{G}} + \frac{c_2}{2} L_{H} M_{H, 0}.
    \end{equation}
    \end{corollary}
    {\it Proof}
        Combining Assumption~\ref{ass: Lipschitz continuity} and Lemma~\ref{lem: gradient of mathcal v} yields the desired result.
    \qed
    
    Assume we have already obtained the approximate optimal solution $(w^k, \lambda^k)$ of the subproblem~\eqref{eq: minimax subproblem} at the $k$-th iteration. Note that $(w^k, \lambda^k)$ is random variables dependent on the inner loop sample $\mathbf{\xi}^k$. In the subsequent analysis, we will take expectation conditioned on $\mathcal{F}^k$ and hence $(w^k, \lambda^k)$ is treated as constants.
    Given sample $\mathbf{\tilde{\xi}}^k =( \tilde \xi^k_1, ..., \tilde \xi^k_{q_k}) \sim \mathcal{D}_{\xi}^{q_k}$, a natural  first order oracle for $\hat{\mathcal{G}}$ in~\eqref{eq: maintext stochastic penalized single level} is considered as
       \begin{equation}
          \label{eq: gradient oracle of mathcal Gk}
          \begin{aligned}
             \nabla_x \mathcal{G}^k(x, y, z; \mathbf{\tilde{\xi}}^k) &= \frac{1}{q_k}\sum_{l=1}^{q_k} \nabla_x g(x, y; \xi^k_l) - \frac{1}{q_k} \sum_{l=1}^{q_k} \nabla_x \ell{\gamma}(x, z, w^k, \lambda^k; \tilde \xi^k_l), \\
             \nabla_y \mathcal{G}^k(x, y, z; \mathbf{\tilde{\xi}}^k) &=\frac{1}{q_k}\sum_{l=1}^{q_k}  \nabla_y g(x, y; \tilde \xi^k_l), \\
             \nabla_z \mathcal{G}^k(x, y, z; \mathbf{\tilde{\xi}}^k) &={\gamma_2} (\lambda^k - z).
          \end{aligned}
       \end{equation}
    Conditioned on $\mathbf{\xi}^k$, the bias of  $\nabla \mathcal{G}^k(x, y, z; \mathbf{\tilde{\xi}}^k)$ arise from 
    the term $\frac{1}{q_k} \sum_{l=1}^{q_k} \nabla_x \ell{\gamma}(x, z, w^k, \lambda^k; \tilde \xi^k_l)$. However when $\mathbf{\xi}^k$ is also treated as random variable, the estimation of bias is much more complicated. 
    We establish two lemmas to control this bias.
    Lemma~\ref{lem: bound of bias of mathcal G} is conditioned on $\mathcal{F}^k$, while Lemma~\ref{lem: bound of  first order oracle of Gk} is conditioned on $\tilde{\mathcal{F}}_{k-1}$. In the following analysis, $\mathcal{G}(\mathbf{u})$ is the abbreviation of $\mathcal{G}(x, y, z)$ for simplicity.
    \begin{lemma}
    \label{lem: bound of bias of mathcal G}
       $\nabla \mathcal{G}^k(\mathbf{u}; \mathbf{\tilde{\xi}}^k)$ has a controllable bias conditioned on $\mathcal{F}^k$ as
       \begin{equation}
          \label{eq: bound of bias of mathcal G}
          \begin{aligned}
             &\quad \| \mathbb{E}_{\mathbf{\tilde{\xi}}^k} [\nabla \mathcal{G}^k(\mathbf{u}; \mathbf{\tilde{\xi}}^k) | \mathcal{F}^k] - \nabla \mathcal{G}(\mathbf{u}) \|\\
              &\leq 
             \frac{\gamma_2}{2}  \|\lambda^k - \lambda^*(x^{k-1}, z^{k-1})\|^2  + (M_{H, 1} + \sqrt{p})\| \lambda^k - \lambda^*(x^{k-1}, z^{k-1})\| \\
             + &(L_G + \frac{1}{\gamma_1}(M_{H, 2}^2 + \gamma_1 L_H \|z\| + \frac{\gamma_1 L_H}{\gamma_2} M_{H, 0} + M_{H, 0}L_H)) \| w^k - w^*(x^{k-1}, z^{k-1})\|,
          \end{aligned}
       \end{equation}
       and a controllable conditional variance as
       \begin{equation}
          \label{eq: bound of variance of mathcal G}
          \begin{aligned}
             \mathbb{V}_{\mathbf{\tilde{\xi}}^k} [\nabla \mathcal{G}^k(\mathbf{u}; \mathbf{\tilde{\xi}}^k) | \mathcal{F}_k] \leq  \frac{  4 \sigma_g^2 }{q_k}.
          \end{aligned}
       \end{equation}
       Here $\mathbb{E}_{\mathbf{\tilde{\xi}}^k}[\cdot], \mathbb{V}_{\mathbf{\tilde{\xi}}^k}[\cdot]$ are the abbreviation of $\mathbb{E}_{\mathbf{\tilde{\xi}}^k \sim \mathcal{D}_{\xi}^{q_k}}[\cdot], \mathbb{V}_{\mathbf{\tilde{\xi}}^k \sim \mathcal{D}_{\xi}^{q_k}}[\cdot]$, respectively.
    \end{lemma}
    {\it Proof}
       When condition on $\mathcal{F}^k$, $(w^k, \lambda^k), (w^*(x^{k-1}, z^{k-1}), \lambda^*(x^{k-1}, z^{k-1}))$ are constants.
       Utilizing~\eqref{eq: gradient of mathcal G 1} and taking expectation over $\mathbf{\tilde{\xi}}^k$ gives
       \begin{equation}
          \begin{aligned}
             &\quad \mathbb{E}_{\mathbf{\tilde{\xi}}^k} [\nabla_x \mathcal{G}^k(\mathbf{u}; \mathbf{\tilde{\xi}}^k)| {\mathcal{F}}_k] - \nabla_x \mathcal{G}(\mathbf{u}) \\
             &=\mathbb{E}_{\mathbf{\tilde{\xi}}^k} \left[ \frac{1}{q_k} \sum_{l=1}^{q_k} \nabla_x g(x, y; \tilde{\xi}_l) -  \frac{1}{q_k} \sum_{l=1}^{q_k} \nabla_x \ell_{\gamma}(x, z, w^k, \lambda^k; \tilde{\xi}_l) | {\mathcal{F}}_k\right] \\
             &\quad - \nabla_x g(x, y) + \nabla_x \ell_{\gamma}(x, z, w^*, \lambda^*)  \\
             &= - \mathbb{E}_{\tilde \xi^k_1} \left[\nabla_x \ell_{\gamma}(x, z, w^k, \lambda^k; \tilde{\xi}_1)| {\mathcal{F}}_k\right] + \nabla_x \ell_{\gamma}(x, z, w^*, \lambda^*)  .
          \end{aligned}
       \end{equation}
       From the expression of $\tilde{\mathcal{L}}$ and \eqref{eq: derivative of L},~\eqref{eq: gradient oracle of L}, we can further compute
       \begin{equation}
          \label{eq: proof of bound of bias of mathcal G 1}
          \begin{aligned}
             &\quad \mathbb{E}_{\mathbf{\tilde{\xi}}^k} [\nabla_x \mathcal{G}^k(\mathbf{u}; \mathbf{\tilde{\xi}}^k)| {\mathcal{F}}_k] - \nabla_x \mathcal{G}(\mathbf{u}) \\
             &= - \mathbb{E}_{\tilde \xi^k_1}  \left[\nabla_x g(x, w^k; \xi_1) + \frac{1}{\gamma_1} \sum_{i=1}^{p} [\gamma_1 \lambda_i^k +  H_i(x, w^k)]_+ \nabla_x H_i(x, w^k) - \frac{\gamma_2}{2} \| \lambda^k - z\|^2 | {\mathcal{F}}_k \right] \\
             &\quad + \nabla_x g(x, w^*) + \frac{1}{\gamma_1} \sum_{i=1}^{p} [\gamma_1 \lambda_i^* +  H_i(x, w^*)]_+ \nabla_x H_i(x, w^*)  - \frac{\gamma_2}{2} \| \lambda^* - z\|^2  \\
             &=- \nabla_x g(x, w^k) - \frac{1}{\gamma_1} \sum_{i=1}^{p} [\gamma_1 \lambda_i^k +  H_i(x, w^k)]_+ \nabla_x H_i(x, w^k) + \frac{\gamma_2}{2} \| \lambda^k - z\|^2  \\
             &\quad + \nabla_x g(x, w^*) + \frac{1}{\gamma_1} \sum_{i=1}^{p} [\gamma_1 \lambda_i^* +  H_i(x, w^*)]_+ \nabla_x H_i(x, w^*)  - \frac{\gamma_2}{2} \| \lambda^* - z\|^2  .
          \end{aligned}
       \end{equation}
       By Assumption~\ref{ass: Lipschitz continuity}, it holds that
       \begin{equation}
          \label{eq: proof of bound of bias of mathcal G 2}
           \begin{aligned}
               &\quad \left\| \frac{1}{\gamma_1} \sum_{i=1}^{p} [\gamma_1 \lambda_i^k +  H_i(x, w^k)]_+ \nabla_x H_i(x, w^k) - \frac{1}{\gamma_1}\sum_{i=1}^{p} [\gamma_1 \lambda_i^* +  H_i(x, w^*)]_+ \nabla_x H_i(x, w^*) \right\|\\
                &\leq \frac{1}{\gamma_1} \left\|\sum_{i=1}^{p} \left([\gamma_1 \lambda_i^k +  H_i(x, w^k)]_+ \nabla_x H_i(x, w^k) -  [\gamma_1 \lambda_i^* +  H_i(x, w^*)]_+ \nabla_x H_i(x, w^k)
                \right)\right\|\\
                &+ \frac{1}{\gamma_1}\left\| \sum_{i=1}^{p} \left([\gamma_1 \lambda_i^* +  H_i(x, w^*)]_+ \nabla_x H_i(x, w^k) - [\gamma_1 \lambda_i^* +  H_i(x, w^*)]_+ \nabla_x H_i(x, w^*) \right)\right\| \\
                &\leq\frac{1}{\gamma_1} \sum_{i=1}^{p} \left( \left| [\gamma_1 \lambda_i^k +  H_i(x, w^k)]_+  -  [\gamma_1 \lambda_i^* +  H_i(x, w^*)]_+ \right| \cdot\left\| \nabla_x H_i(x, w^k) \right\|
                \right)\\
                &+ \frac{1}{\gamma_1} \sum_{i=1}^{p} \left( [\gamma_1 \lambda_i^* +  H_i(x, w^*)]_+  \cdot \left\| \nabla_x H_i(x, w^k) -  \nabla_x H_i(x, w^*) \right\| \right) \\
                &\leq \frac{1}{\gamma_1}M_{H, 1} (\gamma_1 \| \lambda^k - \lambda^* \| +M_{H, 1} \|w^k - w^* \| )   + \frac{1}{\gamma_1} (\gamma_1 \|\lambda^*\| + M_{H, 0})L_{H} \| w^k  -w^*\| \\
                &=M_{H, 1} \| \lambda^k - \lambda^*\| + \frac{1}{\gamma_1}(M_{H, 2}^2 + \gamma_1 \|\lambda^*\|L_H + M_{H, 0}L_H) \| w^k - w^*\|.
             \end{aligned}
       \end{equation}
    The last inequality uses the fact that $| [a]_+ - [b]_+| \leq |a - b|$ and the assumptions that
    $\|H_i(x, w^k) - H_i(x, w^*)\| \leq M_{H, 1} \|w^k - w^*\|$, $\|H_i(x, w^*)\| \leq M_{H, 0}$ and $\| \nabla H_i(x, w^k) - \nabla H_i(x, w^*)\| \leq L_H \|w^k - w^*\|$. 
    By~\eqref{eq: proof of gradient of mathcal G 1}, it holds that $\|\lambda^k - z\|^2 - \|\lambda^* - z\|^2 = \|\lambda^k - \lambda^*\| + 2 (\lambda^k - \lambda^*)^T(\lambda^* - z) \leq \|\lambda^k - \lambda^*\|^2 + \frac{2\sqrt{p}}{\gamma_2} M_{H, 0} \|\lambda^k - \lambda^*\|$. 
    Then by combining~\eqref{eq: proof of bound of bias of mathcal G 1} and~\eqref{eq: proof of bound of bias of mathcal G 2},
    we bound the bias of $\nabla_x \mathcal{G}^k(\mathbf{u}; \mathbf{\tilde{\xi}}^k)$ as
    \begin{equation*}
       \begin{aligned}
          &\quad \| \mathbb{E}_{\mathbf{\tilde{\xi}}^k} [\nabla_x \mathcal{G}^k(\mathbf{u}; \tilde \xi) | \mathcal{F}^k] - \nabla_x \mathcal{G}(\mathbf{u}) \| \\
          &\leq L_G \| w^k - w^* \| +M_{H, 1} \| \lambda^k - \lambda^*\| + \frac{1}{\gamma_1}(M_{H, 2}^2 + \gamma_1 \|\lambda^*\|L_H + M_{H, 0}L_H) \| w^k - w^*\| \\
          &\quad + \frac{\gamma_2}{2}  \|\lambda^k - \lambda^*\|^2 
           + \sqrt{p} M_{H, 0} \|\lambda^k - \lambda^*\| \\
          &= \frac{\gamma_2}{2}  \|\lambda^k - \lambda^*\|^2  + (M_{H, 1} + \sqrt{p})\| \lambda^k - \lambda^*\| \\
          &\quad + (L_G + \frac{1}{\gamma_1}(M_{H, 2}^2 + \gamma_1 \|\lambda^*\|L_H+ M_{H, 0}L_H)) \| w^k - w^*\| \\
          &\leq\frac{\gamma_2}{2}  \|\lambda^k - \lambda^*\|^2  + (M_{H, 1} + \sqrt{p})\| \lambda^k - \lambda^*\| \\
          &\quad + (L_G + \frac{1}{\gamma_1}(M_{H, 2}^2 + \gamma_1 L_H \|z\| + \frac{\gamma_1 L_H}{\gamma_2} M_{H, 0} + M_{H, 0}L_H)) \| w^k - w^*\| .
       \end{aligned}
    \end{equation*}
    The last inequality follows from~\eqref{eq: proof of bound of bias of mathcal G 1} again. Besides
    $
      \mathbb{E} [\nabla_y \mathcal{G}^k(x, y, z; \mathbf{\tilde{\xi}}^k)| \tilde{\mathcal{F}}_k]  - \nabla_y \mathcal{G}(x, y, z) = 0
    $ and $
       \mathbb{E} [\nabla_z \mathcal{G}^k(x, y, z; \mathbf{\tilde{\xi}}^k)| \tilde{\mathcal{F}}_k]  - \nabla_y \mathcal{G}(x, y, z) = \gamma_2 \| \lambda^k - \lambda^* \|
    $ are straightforward.
    Therefore~\eqref{eq: bound of bias of mathcal G} holds. 
    
    From the expression of $\nabla \mathcal{G}^k(\mathbf{u}; \mathbf{\tilde{\xi}}^k)$, we can compute the conditional variance of $\nabla \mathcal{G}^k(x, y, z; \mathbf{\tilde{\xi}}^k)$ as
       \begin{equation*}
          \begin{aligned}
             &\quad \mathbb{V}_{\mathbf{\tilde{\xi}}^k} \left[\nabla_x \mathcal{G}^k(\mathbf{u}; \mathbf{\tilde{\xi}}^k) | {\mathcal{F}}_k\right]  \\
             &= \mathbb{V}_{\mathbf{\tilde{\xi}}^k} \left[
             \frac{1}{q_k} \sum_{l=1}^{q_k} \nabla_x g(x, y; \xi^k_l) - \frac{1}{q_k}\sum_{l=1}^{q_k} \nabla_x \ell{\gamma}(x, z, w^k, \lambda^k; \tilde \xi^k_l) | {\mathcal{F}}_k\right] \\
             &= \frac{1}{q_k} \mathbb{V}_{\tilde \xi^k_1} \left[\nabla_x g(x, y; \tilde \xi^k_1)-\nabla_x  g(x, w^k; \tilde \xi^k_1) - \frac{1}{\gamma_1} \sum_{i=1}^{p} [\gamma_1 \lambda_i^k +  H_i(x, w^k)]_+ \nabla_x H_i(x, w^k) \right.\\
             &\quad \left. + \frac{\gamma_2}{2} \| \lambda^k - z\|^2 | \mathcal{F}_k \right] \\
             &= \frac{1}{q_k} \mathbb{V}_{\tilde \xi^k_1} \left[\nabla_x g(x, y; \tilde \xi^k_1)-\nabla_x  g(x, w^k; \tilde \xi^k_1) | \mathcal{F}_k \right] \quad \text{ ($\lambda^k, w^k$ are constants conditioned on $\mathcal{F}^k$)} \\
             &\leq \frac{2}{q_k} \mathbb{V}_{\tilde \xi^k_1} \left[\nabla_x g(x, y; \tilde \xi^k_1) | \mathcal{F}_k \right] + \frac{2}{q_k} \mathbb{V}_{\tilde \xi^k_1} \left[\nabla_x  g(x, w^k; \tilde \xi^k_1) | \mathcal{F}_k\right] .
             \end{aligned}
             \end{equation*}
            and    
            \begin{equation*}
            \begin{aligned}
             &\quad \mathbb{V}_{\mathbf{\tilde \xi}^k}  \left[\nabla_y \mathcal{G}^k(\mathbf{u}; \tilde \xi^k) | {\mathcal{F}}_k\right] 
             = \mathbb{V}_{\mathbf{\tilde \xi}^k} \left[\frac{1}{q_k} \sum_{l=1}^{q_k} \nabla_y g(x, y; \xi^k_l)\right] 
             = \frac{1}{q_k}\mathbb{V}_{{\tilde \xi}^k_1} \left[\nabla_y g(x, y; \xi^k_1) | \mathcal{F}_k \right] \\
             &\quad \mathbb{V}_{\mathbf{\tilde{\xi}}^k}  \left[\nabla_z \mathcal{G}^k(\mathbf{u}; \mathbf{\tilde{\xi}}^k) | {\mathcal{F}}_k\right] = 0. 
          \end{aligned}
       \end{equation*}
    
       Hence 
       \[
          \begin{aligned}
       &\quad \mathbb{V}_{\mathbf{\tilde{\xi}}^k} [\| \nabla \mathcal{G}^k(\mathbf{u}; \mathbf{\tilde{\xi}}^k) | \mathcal{F}_k\|] \\
       &\leq\frac{2}{q_k} \mathbb{V}_{\tilde \xi^k_1} \left[\nabla_x g(x, y; \tilde \xi^k_1) | \mathcal{F}_k \right] + \frac{2}{q_k} \mathbb{V}_{\tilde \xi^k_1} \left[\nabla_x  g(x, w^k; \tilde \xi^k_1) | \mathcal{F}_k\right] + \frac{1}{q_k}\mathbb{V}_{{\tilde \xi}^k_1} \left[\nabla_y g(x, y; \xi^k_1) | \mathcal{F}_k \right] \\
       &\leq  \frac{2}{q_k} \mathbb{V}_{\tilde \xi^k_1} \left[\nabla g(x, y; \xi^k_1)\right] + \frac{2}{q_k} \mathbb{V}_{\tilde \xi^k_1} \left[\nabla g(x, w^k; \xi^k_1)\right] \leq \frac{ 4 \sigma_g^2 }{q_k}.
       \end{aligned}
       \]
       This completes the proof.
    \qed

    In the following lemma, $w^k, \lambda^k$ are treated as random variables and we try to control the bias and variance of the first order oracle of $\mathcal{G}^k$ conditioned on $\tilde{\mathcal{F}}_{k-1}$.
    
    \begin{lemma}
    \label{lem: bound of  first order oracle of Gk}
        The first order oracle of $\mathcal{G}^k$ has a bounded conditional bias and variance as
        \begin{subequations}
          \begin{align}
          \left\| \mathbb{E}_{\mathbf{\tilde{\xi}}^k, \mathbf{\xi}^k} [\nabla \mathcal{G}^k(\mathbf{u}; \mathbf{\tilde{\xi}}^k) | \tilde{\mathcal{F}}_{k-1}] -  \nabla \mathcal{G}(\mathbf{u}) \right\| &\leq \epsilon_{\mathcal{G}}^k, \label{eq: bound of bias of mathcal Gk} \\
            \mathbb{V}_{\mathbf{\tilde{\xi}}^k, \mathbf{\xi}^k}[\nabla \mathcal{G}^k(\mathbf{u}; \mathbf{\tilde{\xi}}^k) | \tilde{\mathcal{F}}_{k-1}] \leq (\sigma_{\mathcal{G}}^k)^2 , \label{eq: bound of variance of mathcal Gk}
          \end{align}
        \end{subequations}
        where 
        \begin{equation}
          \label{eq: epsilon of mathcal G}
          \begin{aligned}
          \epsilon_{\mathcal{G}}^k &=
          \frac{\gamma_2}{2}  \mathbb{E} [\|\lambda^k - \lambda^*(x^{k-1}, z^{k-1})\|^2  |\tilde{ \mathcal{F}}_{k-1} ] + (M_{H, 1} + \sqrt{p}) \mathbb{E} [\| \lambda^k - \lambda^*(x^{k-1}, z^{k-1})\| | \mathcal{F}_{k-1}] \\
             + &(L_G + \frac{1}{\gamma_1}(M_{H, 2}^2 + \gamma_1 L_H \|z\| + \frac{\gamma_1 L_H}{\gamma_2} M_{H, 0} + M_{H, 0}L_H)) \mathbb{E} [\| w^k - w^*(x^{k-1}, z^{k-1})\| | \tilde{\mathcal{F}}_{k-1} ] \\
          \end{aligned}
        \end{equation}
        and 
        \begin{equation}
          \label{eq: sigma of mathcal G}
          \begin{aligned}
        (\sigma_{\mathcal{G}}^k)^2 &= \frac{4\sigma_g^2}{q_k} + 2 (L_G^2 + \frac{(\gamma_1 M_{\lambda} + M_{H, 0})^2}{\gamma_1^2}L_{H}^2) \mathbb{E}_{\mathbf{\xi}^k} \left[\| w^k - w^*(x^{k-1}, z^{k-1}) \|^2 | \tilde{\mathcal{F}}_{k-1}\right] \\
        & \quad + \gamma_2^2 \mathbb{E}_{\mathbf{\xi}^k} \left[\| \lambda^k - \lambda^*(x^{k-1}, z^{k-1}) \|^2 | \tilde{\mathcal{F}}_{k-1}\right].
          \end{aligned}
       \end{equation}
    
    \end{lemma}
    {\it Proof}
       Taking expectation over $\mathbf{\xi}^k$ in~\eqref{eq: bound of bias of mathcal G} gives~\eqref{eq: bound of bias of mathcal Gk}.
       Utilizing the property of conditional variance and Lemma~\ref{lem: bound of  first order oracle of Gk}, we have
        \begin{equation}
          \label{eq: proof of bound of  first order oracle of Gk 1}
            \begin{aligned} 
            &\quad  \mathbb{V}_{\mathbf{\tilde{\xi}}^k, \mathbf{\xi}^k} [\nabla \mathcal{G}^k(\mathbf{u}; \mathbf{\tilde{\xi}}^k) | \tilde{\mathcal{F}}_{k-1}] \\
                 &= \mathbb{E}_{\mathbf{\xi}_k}[ \mathbb{V}_{\mathbf{\tilde{\xi}}^k} [\nabla \mathcal{G}^k(\mathbf{u}; \mathbf{\tilde{\xi}}^k) | {\mathcal{F}}_k] | \tilde{\mathcal{F}}_{k-1}]
                +   \mathbb{V}_{\mathbf{\xi}^k} [ \mathbb{E}_{\mathbf{\tilde{\xi}}^k}[\nabla \mathcal{G}^k(\mathbf{u}; \mathbf{\zeta}) | {\mathcal{F}}_k] | \tilde{\mathcal{F}}_{k-1}]\\
                 &\leq \frac{4\sigma_g^2}{q_k} + \mathbb{V}_{\mathbf{\xi}^k} [ \nabla_x g(x, y) - \nabla_x \ell_{\gamma}(x, z, w^k, \lambda^k) | \tilde{\mathcal{F}}_{k-1}] + \mathbb{V}_{\mathbf{\xi}^k} [ \nabla_y g(x, y) | \tilde{\mathcal{F}}_{k-1}] \\
                 &\quad + \mathbb{V}_{\mathbf{\xi}^k} [  \gamma_2 (\lambda^k -z) | \tilde{\mathcal{F}}_{k-1}]  \\
                 &= \frac{4\sigma_g^2}{q_k} + \mathbb{V}_{\mathbf{\xi}^k} [ \nabla_x \ell_{\gamma}(x, z, w^k, \lambda^k) | \tilde{\mathcal{F}}_{k-1}] + 0 + \gamma_2^2 \mathbb{V}_{\mathbf{\xi}^k} [\lambda^k | \tilde{\mathcal{F}}_{k-1}]. \\
            \end{aligned}
       \end{equation}
       The inequality is due to~\eqref{eq: gradient oracle of mathcal Gk} and Lemma~\ref{lem: bound of bias of mathcal G}.
        The second term in the right-hand side of~\eqref{eq: proof of bound of  first order oracle of Gk 1} is bounded by
        \begin{equation}
          \label{eq: proof of bound of  first order oracle of Gk 2}
          \begin{aligned}
          &\quad \mathbb{V}_{\mathbf{\xi}^k} [ \nabla_x \ell_{\gamma}(x, z, w^k, \lambda^k) | \tilde{\mathcal{F}}_{k-1}] \\
          &= \mathbb{V}_{\mathbf{\xi}^k} [\nabla_x g(x, w^k) + \frac{1}{\gamma_1} \sum_{i=1}^{p} [\gamma_1 \lambda_i^k +  H_i(x, w^k)]_+ \nabla_x H_i(x, w^k)| \tilde{\mathcal{F}}_{k-1}] \\
          &\leq 2 \mathbb{V}_{\mathbf{\xi}^k} [\nabla_x g(x, w^k) | \tilde{\mathcal{F}}_{k-1}] + \frac{2(\gamma_1 M_{\lambda} + M_{H, 0})^2}{\gamma_1^2}\mathbb{V}_{\mathbf{\xi}^k} [\nabla_x H(x, w^k)| \tilde{\mathcal{F}}_{k-1}] 
          \end{aligned}
          \end{equation}
    
          Since $w^*(x^{k-1}, z^{k-1})$ is a constant conditioned on $\tilde{\mathcal{F}}_{k-1}$, we have 
          $$
          \begin{aligned}
          \mathbb{V}_{\mathbf{\xi}^k} [\nabla_x g(x, w^k) | \tilde{\mathcal{F}}_{k-1}] &= \mathbb{V}_{\mathbf{\xi}^k} [\nabla_x g(x, w^k) - \nabla_x g(x, w^*(x^{k-1}, z^{k-1})) | \tilde{\mathcal{F}}_{k-1}] \\
          &\leq \mathbb{E}_{\mathbf{\xi}^k} [\| \nabla_x g(x, w^k) - \nabla_x g(x, w^*(x^{k-1}, z^{k-1})) \|^2 | \tilde{\mathcal{F}}_{k-1}] \\
          &\leq L_{G}^2 \mathbb{E}_{\mathbf{\xi}^k} \left[ \| w^k - w^*(x^{k-1})\|^2 | \tilde{\mathcal{F}}_{k-1}\right].
          \end{aligned}
          $$
          The first inequality uses the fact that variance is smaller than second order moment, and the second inequality is due to the Lipschitz continuity of $G$. Similarly, we have 
          $$\mathbb{V}_{\mathbf{\xi}^k} [\nabla_x H(x, w^k) | \tilde{\mathcal{F}}_{k-1}] \leq L_{H}^2 \mathbb{E}_{\mathbf{\xi}^k} \left[ \| w^k - w^*(x^{k-1})\|^2 | \tilde{\mathcal{F}}_{k-1}\right]. 
          $$ Substituting these two inequalities into~\eqref{eq: proof of bound of  first order oracle of Gk 2} gives
          \begin{equation}
          \begin{aligned}    
             \label{eq: proof of bound of  first order oracle of Gk 3}
           &\quad \mathbb{V}_{\mathbf{\xi}^k} [ \nabla_x \ell_{\gamma}(x, z, w^k, \lambda^k) | \tilde{\mathcal{F}}_{k-1}] \\
          &\leq 2 \left(L_{G}^2 + \frac{(\gamma_1 M_{\lambda} + M_{H, 0})^2}{\gamma_1^2} L_{H}^2 \right)\mathbb{E} [\| w^k - w^*(x^{k-1}, z^{k-1}) \|^2 | \tilde{\mathcal{F}}_{k-1}].
          \end{aligned}
        \end{equation}
        The last term in the right-hand side of~\eqref{eq: proof of bound of  first order oracle of Gk 1} is bounded by
        \begin{equation}
          \label{eq: proof of bound of  first order oracle of Gk 4}
           \begin{aligned}
              \gamma_2^2\mathbb{V}_{\mathbf{\xi}^k} [\lambda^k | \tilde{\mathcal{F}}_{k-1}] 
              &= \gamma_2^2 \mathbb{V} [ \lambda^k - \lambda^*(x^{k-1}, z^{k-1}) | \tilde{\mathcal{F}}_{k-1} ]
              \leq \gamma_2^2\mathbb{E} [\| \lambda^k - \lambda^*(x^{k-1}, z^{k-1}) \|^2 | \tilde{\mathcal{F}}_{k-1}] .
           \end{aligned}
        \end{equation}
          Combining~\eqref{eq: proof of bound of  first order oracle of Gk 1},~\eqref{eq: proof of bound of  first order oracle of Gk 3} and~\eqref{eq: proof of bound of  first order oracle of Gk 4} gives the desired result. 
    \qed

    \begin{lemma}
       \label{lem: bound of bias of mathcal G and variance of mathcal G}
       Under 
       the conditions of Theorem~\ref{thm: point convergence of inner loop}, $\epsilon_{\mathcal{G}}^k$ and $\sigma_{\mathcal{G}}^k $ are bounded as
       \begin{subequations}
          \begin{align}
             \| \epsilon_{\mathcal{G}}^k \| &\leq \frac{\gamma_2}{2}  \phi_2\frac{1+\log(s_k)}{s_k} + (M_{H, 1} + \sqrt{p}) \left(\phi_2 \frac{1+\log(s_k)}{s_k}\right)^{\frac{1}{2}} \\
             + &\left(L_G + \frac{1}{\gamma_1}(M_{H, 2}^2 + \gamma_1 L_H B + \frac{\gamma_1 L_H}{\gamma_2} M_{H, 0} + M_{H, 0}L_H)\right) \left(\phi_1\frac{1+\log(s_k)}{s_k}\right)^{\frac{1}{2}},\\
             (\sigma_{\mathcal{G}}^k)^2  &\leq \frac{2\sigma_g^2}{q_k} + ((L_G^2 + \frac{(\gamma_1 \bar M_{\lambda} + M_{H, 0})^2}{\gamma_1^2}L_{H}^2) \bar{\phi}_1 + \gamma_2^2\bar{\phi}_2 )\frac{1+\log(s_k)}{s_k}.
          \end{align}
       \end{subequations}
       where $\bar{\phi}_1 = \eta\left( \eta M_{\mathcal{L}, 1} + \rho (4 \gamma_2 B + 4 M_{H, 0}^2 ) + (\eta M_{\mathcal{L}, 2} + 2 \rho (\gamma_1 + \gamma_2)) M_{\lambda} \right)$, $\bar{\phi}_2 = \frac{\rho}{\eta} \bar{\phi}_1$, $\bar M_{\lambda}=2 \rho^2 \gamma_2^2 B^2 + 2 p \rho^2 M_{H, 0}^2$ are constants obtained by replacing  $\|z\|$ with  $B$.
    \end{lemma}
    
    {\it Proof}
    It follows from  Cauchy-Schwarz inequality that $\mathbb{E} [\| \lambda^k - \lambda^*(x^{k-1}, z^{k-1})\| | \mathcal{F}_k] \leq (\mathbb{E} [\| \lambda^k - \lambda^*(x^{k-1}, z^{k-1})\|^2 | \mathcal{F}_k])^{\frac{1}{2}} $ and $\mathbb{E} [\| w^k - w^*(x^{k-1}, z^{k-1})\| |  \mathcal{F}_k] \leq (\mathbb{E} [\| w^k - w^*(x^{k-1}, z^{k-1})\|^2 | \mathcal{F}_k])^{\frac{1}{2}} $
       Substituting the results of Theorem~\ref{thm: point convergence of inner loop} into~\eqref{eq: epsilon of mathcal G} and~\eqref{eq: sigma of mathcal G} gives the desired results. Note that the notation $ w^s - w^*, \lambda^s - \lambda^*$ in~\eqref{eq: point convergence of inner loop} are replaced by $w^k - w^*(x^{k-1}, z^{k-1}), \lambda^k - \lambda^*(x^{k-1}, z^{k-1})$ in the current context. 
    \qed
    
    Denote $\nabla f(x, y; \mathbf{\zeta}^k) = \frac{1}{r_k} \sum_{l=1}^{r_k} \nabla f(x, y; \zeta^k_l)$.
    From the definition of $\Psi^k(\mathbf{u}; \mathbf{\zeta}^k, \mathbf{\tilde{\xi}}^k)$ in~\eqref{eq: stochastic gradient of Psi} and~\eqref{eq: Psi k}, we derive the relationship between the bias of $\nabla \Psi^k(\mathbf{u}; \mathbf{\zeta}^k, \mathbf{\tilde{\xi}}^k)$ and $\nabla \mathcal{G}^k(\mathbf{u}; \mathbf{\tilde{\xi}}^k)$ as follows.
    \begin{equation}
       \label{eq: relationship between bias of Psi and G}
       \nabla \Psi^k(\mathbf{u}; \mathbf{\zeta}^k, \mathbf{\tilde{\xi}}^k) - \nabla \Psi(\mathbf{u}) = \nabla f(x, y; \mathbf{\zeta}^k) - \nabla F(x, y) + c_1(\mathcal{G}^k(\mathbf{u}; \mathbf{\tilde{\xi}}^k) - \nabla \mathcal{G}(\mathbf{u})).
    \end{equation}
    Since $\nabla f(x, y; \mathbf{\zeta}^k)$ is unbiased, the bias of $\nabla \Psi^k(\mathbf{u}; \mathbf{\zeta}^k, \mathbf{\tilde{\xi}}^k)$ is fully determined by the bias of $\nabla \mathcal{G}^k(\mathbf{u}; \mathbf{\tilde{\xi}}^k)$. 
    
    From the relationship between variance and moment, the bias $b^k$ defined in~\eqref{eq: bk} can be bounded as follow.
    \begin{lemma}
       \label{lem: bound of bias}
        The bias $b^k$ has a controllable momentum as 
        \begin{equation}
           \label{eq: bound of bias}
           \begin{aligned}
            \mathbb{E} [ \| b^k \|^2 | \tilde{\mathcal{F}}_{k-1}] 
            &\leq   
             2 \left(\frac{\sigma_{f}^2}{r_k} + c_1^2 \frac{(\sigma_{\mathcal{G}}^k)^2}{q_k} + c_1^2(\epsilon_{\mathcal{G}}^k)^2  \right).
       \end{aligned}
        \end{equation}
        Here $\mathbb{E}[\cdot] $ is the abbreviation of $\mathbb{E}_{\mathbf{\zeta}^k, \mathbf{\tilde{\xi}}^k, \mathbf{\xi}^k}[\cdot]$.
    \end{lemma}
    {\it Proof}
    By Lemma~\ref{lem: bound of  first order oracle of Gk}, it holds that 
    \begin{equation*}
       \begin{aligned}
       &\quad \mathbb{E} [ \| b^k \|^2 | \widetilde{\mathcal{F}}_{k-1}] \\
       &\leq 2 \mathbb{V} [ \nabla \Psi^k(\mathbf{u}^k; \mathbf{\zeta}^k, \mathbf{\tilde{\xi}}^k) | \widetilde{\mathcal{F}}_{k-1} ] + 2 \mathbb{E} [\| \mathbb{E} [\nabla \Psi^k(\mathbf{u}^k; \mathbf{\zeta}^k, \mathbf{\tilde{\xi}}^k) | \widetilde{\mathcal{F}}_{k-1} ] - \nabla \Psi^k(\mathbf{u}^k) \|^2 | \widetilde{\mathcal{F}}_{k-1}] \\
       &=2\left( \mathbb{V} \left[ \nabla f(\mathbf{u}^k; \mathbf{\zeta}^k) | \widetilde{\mathcal{F}}_{k-1} \right] + \mathbb{V} \left[c_1  \nabla \mathcal{G} (\mathbf{u}^k; \mathbf{\tilde{\xi}}^k) | \widetilde{\mathcal{F}}_{k-1} \right] \right. \\
       &\quad \left. + \mathbb{E} \left[c_1^2 \|\mathbb{E} [
     \nabla \mathcal{G}(\mathbf{u}^k; \mathbf{\zeta}^k) | \widetilde{\mathcal{F}}_{k-1} ] - \nabla \mathcal{G}(\mathbf{u}^k) \|^2 | \widetilde{\mathcal{F}}_{k-1}\right]  \right)\\
       &\leq2 \left(\frac{\sigma_{f}^2}{r_k} + c_1^2 \frac{\sigma_{\mathcal{G}}^2 }{q_k} + c_1^2(\epsilon_{\mathcal{G}}^k)^2 \right).
       \end{aligned}
    \end{equation*}
    The last inequality follow from Lemma~\ref{lem: bound of bias of mathcal G} and~\ref{lem: bound of  first order oracle of Gk}. 
    This completes the proof. 
    \qed

    The convergence of Algorithm~\ref{alg: main} is established in the following theorem.
    \begin{theorem}
       \label{thm: convergence of outer loop}
       Assume  the step sizes satisfy $\alpha_k < \frac{1}{2 L_{\Psi}}$,
       then the sequence $\{\mathbf{u}^k\}$ satisfies
    \begin{equation}
    \label{eq: convergence of outer loop}
    \begin{aligned}
        &\quad \mathbb{E} \left[\frac{1}{\sum_{k=0}^{K-1} \alpha_k}\sum_{k=0}^{K-1} \frac{1}{\alpha_k} \|\mathbf{u}^{k+1} - \mathbf{u}^k \|^2 \right] \\
        &
          \leq \frac{4}{\sum_{k=0}^{K-1} \alpha_k} \mathbb{E} [\Psi(\mathbf{u}^0) - \Psi(\mathbf{u}^K)] +  \frac{2}{\sum_{k=0}^{K-1} \alpha_k}\sum_{k=0}^{K-1} \alpha_k \mathbb{E} [\|b^k\|^2].
    \end{aligned}
    \end{equation}
    Here the expectation is taken over all $\mathbf{\xi^k} \sim \mathcal{D}_{\xi}^{s_k}, \mathbf{\tilde{\xi}}^k \sim \mathcal{D}_{{\xi}}^{q_k}, \mathbf{\zeta}^k \sim \mathcal{D}_{\zeta}^{r_k}$, $k = 0, ..., K-1$.
    \end{theorem}
    {\it Proof}
       The projection gradient step $\mathbf{u}^{k+1} = \mathrm{Proj}_{\mathcal{U}}(\mathbf{u}^k - \alpha_k \nabla \Psi^k(\mathbf{u}^k; \mathbf{\zeta}^k))$ gives
       \[
          \langle \mathbf{u}^{k+1} - \mathbf{u}^k, \mathbf{u}^{k+1} - (\mathbf{u}^k - \alpha_k \nabla \Psi^k(\mathbf{u}^k; \mathbf{\zeta}^k), \mathbf{\tilde{\xi}}^k)  \rangle \leq 0 .
       \]
       This is equivalent to
       \[
       \langle \nabla \Psi^k(\mathbf{u}^k; \mathbf{\zeta}^k), \mathbf{u}^{k+1} - \mathbf{u}^k \rangle \leq - \frac{1}{\alpha_k} \|\mathbf{u}^{k+1} - \mathbf{u}^k \|^2   .
       \]
       The Lipschitz property of $\nabla \Psi(\mathbf{u})$ gives
       \begin{equation}
       \label{eq: proof of convergence of outer loop 1}
       \begin{aligned}
           &\quad \Psi(\mathbf{u}^{k+1}) - \Psi(\mathbf{u}^k)  \\
          &\leq \Psi(x^{k+1}, y^{k+1}, z^{k+1}, \nu^k, \rho^k) - \Psi(x^k, y^k, z^k, \nu^{k+1}, \rho^{k+1}) \\
          &\leq  \langle \nabla \Psi(\mathbf{u}^k), \mathbf{u}^{k+1} - \mathbf{u}^k \rangle + \frac{L_{\Psi}}{2} \|\mathbf{u}^{k+1} - \mathbf{u}^k \|^2 \\
          &= \langle \nabla \Psi^k(\mathbf{u}^k; \mathbf{\zeta}^k, \mathbf{\tilde{\xi}}^k), \mathbf{u}^{k+1} - \mathbf{u}^k \rangle + \langle b^k, \mathbf{u}^{k+1} - \mathbf{u}^k \rangle + \frac{L_{\Psi}}{2} \|\mathbf{u}^{k+1} - \mathbf{u}^k \|^2 \\
          &\leq - \frac{1}{\alpha_k} \|\mathbf{u}^{k+1} - \mathbf{u}^k \|^2 + \frac{\alpha_k}{2}\|b^k\|^2 + \frac{1}{2\alpha_k}\|\mathbf{u}^{k+1} - \mathbf{u}^k \|^2 +  \frac{L_{\Psi}}{2} \|\mathbf{u}^{k+1} - \mathbf{u}^k \|^2\\
          &= -\frac{1}{2} (\frac{1}{\alpha_k} - L_{\Psi}) \|\mathbf{u}^{k+1} - \mathbf{u}^k \|^2 + \frac{\alpha_k}{2} \|b^k\|^2 \\
          &\leq -\frac{1}{4 \alpha_k} \|\mathbf{u}^{k+1} - \mathbf{u}^k \|^2 + \frac{\alpha_k}{2} \|b^k\|^2.
       \end{aligned}   
    \end{equation}
    The last inequality uses the assumption $\alpha_k \leq \frac{1}{2L_{\Psi}}$.
    Then summing up~\eqref{eq: proof of convergence of outer loop 1} over $k=0,..., K-1$ gives
    \begin{equation}
        \label{eq: proof of convergence of outer loop 2}
    \begin{aligned}
       \sum_{k=0}^{K-1} \frac{1}{\alpha_k}\|\mathbf{u}^{k+1} - \mathbf{u}^k \|^2 
       \leq& 4 (\Psi(\mathbf{u}^0) - \Psi(\mathbf{u}^K)) 
       + 2 \sum_{k=0}^{K-1} \alpha_k \|b^k\|^2.
    \end{aligned}   
    \end{equation}
    
    Taking expectation on both sides and multiplying $\frac{1}{\sum_{k=0}^{K-1} \alpha_k}$
    yields~\eqref{eq: convergence of outer loop}.
    This completes the proof.
    \qed
    
    We consider measuring the convergence by the deviation from the optimal condition as
    \begin{equation}
        \mathrm{dist}(0, \nabla \Psi(\mathbf{u})  + \mathcal{N}_{\mathcal{U}}(\mathbf{u}) ).
    \end{equation}
    Let $\delta^k = \mathbf{u}^{k+1} - \mathbf{u}^k + \alpha^k \nabla \Psi^k(\mathbf{u}^k; \mathbf{\zeta}^k, \mathbf{\tilde{\xi}}^k)$. The projection step gives $- \delta^k \in \mathcal{N}_{\mathcal{U}}(\mathbf{u}^{k+1})$.
    We can derive the following bound on this measure in Theorem~\ref{thm: gradient bound of Psi}.
    
    \begin{theorem}
       \label{thm: gradient bound of Psi}
       Assume $\alpha_k \leq \frac{1}{2 L_{\Psi}}$.
       Then it holds that 
    \begin{equation}
    \label{eq: gradient bound of Psi}
    \begin{aligned}
     &\quad \frac{1}{\sum_{k=0}^{K-1} \alpha_k}\sum_{k=0}^{K-1} \alpha_k \mathbb{E} [\mathrm{dist}(0, \nabla \Psi(\mathbf{u}^{k+1}) + \mathcal{N}_{\mathcal{U}}(\mathbf{u}^{k+1}))^2] 
         \\
         &\leq \frac{18}{\sum_{k=0}^{K-1} \alpha_k} (\Psi(\mathbf{u}^0) - \Psi(\mathbf{u}^K)) + \frac{11}{\sum_{k=0}^{K-1} \alpha_k} \sum_{k=0}^{K-1} \alpha_k \mathbb{E}[\|b^k\|^2] .
    \end{aligned}
    \end{equation}
    
    \end{theorem}
    \label{thm: convergence of outer loop 1}
    {\it Proof}
       Since $-\delta^k \in \mathcal{N}_{\mathcal{U}}(\mathbf{u}^{k+1})$, we have
       \begin{equation}
          \label{eq: convergence of outer loop 1}
           \mathrm{dist}(0, \nabla \Psi^k(\mathbf{u}^{k}) + \mathcal{N}_{\mathcal{U}}(\mathbf{u}^{k+1}))
           \leq \frac{1}{\alpha_k} \| \mathbf{u}^{k+1} - \mathbf{u}^k\| .
       \end{equation}
       The Lipschitz property of $\nabla \Psi(\mathbf{u})$ gives
       \begin{equation}
       \begin{aligned}
          \mathrm{dist}(0, \nabla \Psi(\mathbf{u}^{k+1}) + \mathcal{N}_{\mathcal{U}}(\mathbf{u}^{k+1})) 
           &\leq 
           \mathrm{dist}(0, \nabla \Psi(\mathbf{u}^{k}) + \mathcal{N}_{\mathcal{U}}(\mathbf{u}^{k+1}))
           + L_{\Psi} \| \mathbf{u}^{k+1} - \mathbf{u}^k\|\\
           &\leq 
           \mathrm{dist}(0, \nabla \Psi^k(\mathbf{u}^{k}) + \mathcal{N}_{\mathcal{U}}(\mathbf{u}^{k+1})) + \| b^k\|
           + L_{\Psi} \| \mathbf{u}^{k+1} - \mathbf{u}^k\|\\
           &\leq \|b^k\| + (\frac{1}{\alpha_k} + L_{\Psi}) \| \mathbf{u}^{k+1} - \mathbf{u}^{k}\| \quad \text{( by~\eqref{eq: convergence of outer loop 1})}\\
           &  \leq \|b^k\| + \frac{3}{2\alpha_k} \| \mathbf{u}^{k+1} - \mathbf{u}^{k}\|.
           \end{aligned}
       \end{equation}
       Taking square  we obtain
       \begin{equation}
       \begin{aligned}
          \mathrm{dist}(0, \nabla \Psi(\mathbf{u}^{k+1}) + \mathcal{N}_{\mathcal{U}}(\mathbf{u}^{k+1}))^2 
          & \leq (\|b^k\| + \frac{3}{2\alpha_k} \| \mathbf{u}^{k+1} - \mathbf{u}^{k}\|)^2 \\
          &\leq 2 \|b^k\|^2 + \frac{9}{2\alpha_k^2} \| \mathbf{u}^{k+1} - \mathbf{u}^{k}\|^2.
       \end{aligned}
       \end{equation}
       Substituting the above inequality into~\eqref{eq: convergence of outer loop} gives~\eqref{eq: gradient bound of Psi}. This completes the proof.
    \qed
    
    \begin{corollary}
    \label{cor: rate of outer loop}
       Take step size as~\eqref{eq: step size of inner loop} in the inner loop.     
       Take constant step size $\alpha_k = {\alpha} < \frac{1}{2L_{\Psi}}$ in the outer loop and take constant sample sizes as $r_k = r$, $q_k = q$ and $s_k = s$ in Algorithm~\ref{alg: main}. 
        Randomly choosing an index $R$ from $\{1, ..., K\}$ with probability $\mathrm{Prob}(R = k) = \frac{\alpha_{k-1}}{\sum_{k=1}^{K} \alpha_{k-1}}$,
        then we have 
        \begin{equation*}
    \begin{aligned}
    &\quad \mathbb{E} \left[\frac{1}{\sum_{k=0}^{K-1} \alpha_k} \sum_{k=0}^{K-1} \frac{1}{\alpha_k} \|\mathbf{u}^{k+1} - \mathbf{u}^{k} \|^2\right] \leq \widetilde{\mathcal{O}}\left(\frac{1}{\alpha K} + \frac{1}{r} + \frac{c_1^2}{q} + \frac{c_1^2}{s} \right)     , \\
    &\quad \mathbb{E} [\mathrm{dist}(0, \nabla \Psi(\mathbf{u}^R) + \mathcal{N}_{\mathcal{U}}(\mathbf{u}^{R}))^2] \leq \widetilde{\mathcal{O}}\left(\frac{1}{\alpha K} + \frac{1}{r} + \frac{c_1^2}{q} + \frac{c_1^2}{s} \right).
    \end{aligned}
    \end{equation*}
    \end{corollary}
    {\it Proof}
    From Lemma~\ref{lem: bound of bias of mathcal G and variance of mathcal G},~\ref{lem: bound of bias} we know that
    $\mathbb{E}[\|b^k\|^2] \leq \widetilde{\mathcal{O}}(\frac{1}{r_k} + \frac{c_1^2}{q_k} + \frac{c_1^2}{s_k})
    $. Substituting this into Theorem~\ref{thm: convergence of outer loop} and~\ref{thm: gradient bound of Psi} gives
       \begin{equation*}
    \begin{aligned}
    \mathbb{E} \left[\frac{1}{\sum_{k=0}^{K}\alpha_k}\sum_{k=0}^{K-1} \frac{1}{\alpha_k} \|\mathbf{u}^{k+1} - \mathbf{u}^k \|^2 \right]
          &\leq \frac{4}{\alpha K} \mathbb{E} [\Psi(\mathbf{u}^0) - \Psi(\mathbf{u}^K)] +  \frac{2}{K}\sum_{k=0}^{K-1}  \mathbb{E} [\|b^k\|^2] \\
          &\leq \widetilde{\mathcal{O}}\left(\frac{1}{\alpha K} + \frac{1}{r} + \frac{c_1^2}{q} + \frac{c_1^2}{s} \right),  \\
    \mathbb{E} [\mathrm{dist}(0, \nabla \Psi(\mathbf{u}^R) + \mathcal{N}_{\mathcal{U}}(\mathbf{u}^{R}))^2]  
          &\leq \frac{18}{\alpha K} \mathbb{E} [\Psi(\mathbf{u}^0) - \Psi(\mathbf{u}^K)] +  \frac{11}{K}\sum_{k=0}^{K-1}  \mathbb{E} [\|b^k\|^2] \\
          &\leq \widetilde{\mathcal{O}}\left(\frac{1}{\alpha K} +\frac{1}{r} + \frac{c_1^2}{q} + \frac{c_1^2}{s}\right) .
    \end{aligned}
    \end{equation*}
    This completes the proof.
    \qed
    
    \begin{remark}
       Let $c = \max(c_1, c_2)$.
       The step size condition $\alpha_k < \frac{1}{2L_{\Psi}}$ and~\eqref{eq: Lipschitz continuity of Psi} implies $\alpha_k$ is a most $\widetilde{\mathcal{O}}(c^{-1})$. 
       With $\alpha \sim \mathcal{O}(c^{-1})$, 
       $r \sim \mathcal{O}(\epsilon^{-1})$,
       $q \sim \mathcal{O}(c_1^2 \epsilon^{-1})$,  $s \sim {\mathcal{O}}(c_1^2 \epsilon^{-1})$, $K \sim {\mathcal{O}}(c \epsilon^{-1})$, the right side of the above inequality is $\widetilde{\mathcal{O}}(\epsilon)$.
      Then the sample complexity on $\xi$ is $\sum_{k=0}^{K-1} {s_k} +\sum_{k=1}^{k} q_k  = sK + qK = \widetilde{\mathcal{O}}(c c_1^2 \epsilon^{-2})$ and the sample complexity on $\zeta$ is $\sum_{k=0}^{K-1} {r_k} = r K = \widetilde{\mathcal{O}}(c\epsilon^{-2})$. Theorem~\ref{thm: convergence of outer loop} shows the algorithm converges to the $\epsilon$-stationary point of the problem~\eqref{eq: maintext stochastic penalized single level} for any fixed $c_1 > 0$ with this sample complexity.
    \end{remark}
    
    \begin{remark}
       By Theorem~\ref{thm: equivalence of stochastic single level}, if we take $c_1 \sim \mathcal{O}(\epsilon^{-1})$, $c_2 \sim \mathcal{O}(\epsilon^{-3}), \delta \sim \mathcal{O}(\epsilon^{-2})$, then~\eqref{eq: maintext stochastic penalized single level} is equivalent to the original BLO~\eqref{eq: BLO} in the sense of $\epsilon$-accuracy. Under this condition, the sample complexity on $\xi$ is $\widetilde{\mathcal{O}}(\epsilon^{-7})$ and the sample complexity on $\zeta$ is $\widetilde{\mathcal{O}}( \epsilon^{-5})$. 
    \end{remark}
    

    \subsection{Analysis on variance reduction}
    In this section, 
    we introduce a stronger assumption on the Lipschitz continuity of the gradient of $f(x, y)$ and $g(x, y)$ as Assupmtion~\ref{ass: averaged Lipschitz}.
    From~\eqref{lem: bound of  first order oracle of Gk} and \eqref{eq: relationship between bias of Psi and G}, we know $\nabla \Psi^k(\mathbf{u}; \mathbf{\zeta}^k, \mathbf{\tilde{\xi}}^k)$ is $L_{\Psi}'$-averaged Lipschitz continuous
    conditioned on $\mathcal{F}_{k-1}$
    with module
    \begin{equation}
       \label{eq: averaged Lipschitz}
    L_{\Psi}' = L_f + c_1 \epsilon_{\mathcal{G}}^k \leq \mathcal{O}(L_f + \frac{c_1}{s_k}).
    \end{equation}
    Define the error of the direction as
    $$
    e^k = d^k - \mathbb{E} [\nabla \Psi^k (\mathbf{u}^k; \mathbf{\zeta}^k, \mathbf{\tilde{\xi}}^k) | \tilde{\mathcal{F}}_{k-1}].
    $$
    First we derive the decrease in variance reduction by a single iteration.
    
    \begin{lemma}
    \label{lem: decrease in variance reduction}
    The sequence $\{\mathbf{u}^k\}$ satisfies
         \begin{equation}
    \label{eq: decrease in variance reduction}
          \quad \Psi(\mathbf{u}^{k+1}) - \Psi(\mathbf{u}^k) \leq  -\frac{1}{4\alpha_k} \|  \mathbf{u}^{k+1} - \mathbf{u}^{k} \|^2  + \alpha_k \|e^k\|^2  + \alpha_k c_1^2 (\epsilon_{\mathcal{G}}^k)^2 ,
    \end{equation}
    where $\epsilon_{\mathcal{G}}^k$ is defined in Lemma~\ref{lem: bound of  first order oracle of Gk}.
    \end{lemma}
    {\it Proof}
    Similar to~\eqref{eq: proof of convergence of outer loop 1}, it holds that 
     \begin{equation}
       \begin{aligned}
          &\quad \Psi(\mathbf{u}^{k+1}) - \Psi(\mathbf{u}^k)  \\
          &\leq  \langle \nabla \Psi(\mathbf{u}^k), \mathbf{u}^{k+1} - \mathbf{u}^k \rangle + \frac{L_{\Psi}}{2} \|\mathbf{u}^{k+1} - \mathbf{u}^k \|^2 \\
          &=  \langle d^k, \mathbf{u}^{k+1} - \mathbf{u}^k \rangle - \langle  e^k + \mathbb{E} [\nabla \Psi^k (\mathbf{u}^k; \mathbf{\zeta}^k, \mathbf{\tilde{\xi}}^k) | \widetilde{\mathcal{F}}_{k-1}] - \nabla \Psi(\mathbf{u}^k) , \mathbf{u}^{k+1} - \mathbf{u}^k \rangle \\
          &\quad +\frac{L_{\Psi}}{2} \|\mathbf{u}^{k+1} - \mathbf{u}^k \|^2 \\
          &\leq -\frac{1}{ \alpha_k} \|  \mathbf{u}^{k+1} - \mathbf{u}^{k} \|^2 -  \langle  e^k+ \mathbb{E} [\nabla \Psi^k (\mathbf{u}^k; \mathbf{\zeta}^k, \mathbf{\tilde{\xi}}^k) | \widetilde{\mathcal{F}}_{k-1}] - \nabla \Psi(\mathbf{u}^k), \mathbf{u}^{k+1} - \mathbf{u}^k \rangle \\
          &\quad + \frac{L_{\Psi}}{2} \|\mathbf{u}^{k+1} - \mathbf{u}^k \|^2 \\
            &\leq -\frac{1}{ \alpha_k} \|  \mathbf{u}^{k+1} - \mathbf{u}^{k} \|^2 + \frac{\alpha_k}{2} \|  e^k + \mathbb{E} [\nabla \Psi^k (\mathbf{u}^k; \mathbf{\zeta}^k, \mathbf{\tilde{\xi}}^k)| \widetilde{\mathcal{F}}_{k-1}]- \nabla \Psi(\mathbf{u}^k)  \|^2 \\
            &\quad + (\frac{1}{2\alpha_k} + L_{\Psi})\| \mathbf{u}^{k+1} - \mathbf{u}^k \|^2 \\
          &\leq -\frac{1}{2} (\frac{1}{\alpha_k} - L_{\Psi})\|  \mathbf{u}^{k+1} - \mathbf{u}^{k} \|^2  + \alpha_k \|e^k\|^2  + \alpha_k \| \nabla \Psi(\mathbf{u}^k) -  \mathbb{E} [\nabla \Psi^k (\mathbf{u}^k; \mathbf{\zeta}^k, \mathbf{\tilde{\xi}}^k) | \widetilde{\mathcal{F}}_{k-1}] \|^2 \\ 
          &\leq  -\frac{1}{2} (\frac{1}{\alpha_k} - L_{\Psi}) \|  \mathbf{u}^{k+1} - \mathbf{u}^{k} \|^2  + \alpha_k \|e^k\|^2  + \alpha_k c_1^2 (\epsilon_{\mathcal{G}}^k)^2 .
       \end{aligned}   
    \end{equation}
    The first inequality follows from the Lipschitz property, the second inequality is due to the projection gradient step and the third inequality uses Young's inequality.
    This completes the proof. 
    \qed
    
    Next we show that the sequence of error $\{e^k\}$ has recursive relationship as follows, hence the error is decreasing.
    
    \begin{lemma}
    \label{lem: bound of error in variance reduction}
        The conditional expectation of the error $e^{k+1}$ is bounded as
        \begin{equation}
            \label{eq: bound of error in variance reduction}
    \begin{aligned}
            \mathbb{E} [\| e^{k+1}\|^2| \widetilde{\mathcal{F}}_k]
            &\leq {2}\beta_{k+1}^2 \left( \frac{\sigma_{f}^2}{r_{k+1}} + c_1^2 (\sigma_{\mathcal{G}}^k)^2 \right)
            + 8 (1 - \beta_{k+1})^2  \|e^{k} \|^2 \\
            &\quad +\left. 8 \frac{(L_{\Psi}')^2 + L_{\Psi}^2}{r_{k+1}} \mathbb{E}[\|\mathbf{u}^{k+1} - \mathbf{u}^k \|^2 | \tilde{\mathcal{F}}_{k}]  + \frac{8c^2}{r_{k+1}} \big((\epsilon_{\mathcal{G}}^k)^2 +(\epsilon_{\mathcal{G}}^{k+1})^2 \big)\right)
    \end{aligned}
        \end{equation}
        Here the expectation is taken over $\mathbf{\zeta}^{k+1}, \mathbf{\tilde{\xi}}^{k+1}, \mathbf{\xi}^{k+1}$.
    \end{lemma}
    {\it Proof}
    Let $\Delta^{k+1}$ be the bias of the gradient at $\mathbf{u}^{k+1}$ and $\tilde{\Delta}^{k}$ be the error of the gradient at $\mathbf{u}^{k+1}$ comparing the expectation of the gradient of the last iteration, that is,
    $$
    \begin{aligned}
       \Delta^{k+1} &=   \nabla \Psi^{k+1}(\mathbf{u}^{k+1}; \mathbf{\zeta}^{k+1}, \mathbf{\tilde{\xi}}^{k+1}) -  \mathbb{E} [\nabla \Psi^{k+1} (\mathbf{u}^{k+1}; \mathbf{\zeta}^{k+1}, \mathbf{\tilde{\xi}}^{k+1}) | \tilde{\mathcal{F}}_k], \\
       \tilde \Delta^{k} &=  \nabla \Psi^{k+1}(\mathbf{u}^{k}; \mathbf{\zeta}^{k+1}, \mathbf{\tilde{\xi}}^{k+1}) -  \mathbb{E} [\nabla \Psi^{k}(\mathbf{u}^{k}; \mathbf{\zeta}^{k}, \mathbf{\tilde{\xi}}^{k}) | \tilde{\mathcal{F}}_{k-1}].
    \end{aligned}
    $$ 
    
       From the definition of $e^k$ and~\eqref{eq: direction update in variance reduction} we have
        \begin{equation}
        \label{eq: proof of bound of error in variance reduction 1}
            \begin{aligned}
                e^{k+1} 
                &=d^{k+1} - \mathbb{E} [\nabla \Psi^{k+1}(\mathbf{u}^{k+1}; \mathbf{\zeta}^{k+1}, \mathbf{\tilde{\xi}}^{k+1}) | \tilde{\mathcal{F}}_k ] \\
                &= \nabla \Psi^{k+1}(\mathbf{u}^{k+1}; \mathbf{\zeta}^{k+1}, \mathbf{\tilde{\xi}}^{k+1}) - \mathbb{E} [\nabla \Psi^{k+1}(\mathbf{u}^{k+1}; \mathbf{\zeta}^{k+1}, \mathbf{\tilde{\xi}}^{k+1}) | \tilde{\mathcal{F}}_k ] \\
                &\quad  + (1 - \beta_{k+1})(d^{k} - \nabla \Psi^{k+1}(\mathbf{u}^{k}; \mathbf{\zeta}^{k+1}, \mathbf{\tilde{\xi}}^{k+1})) \\
                &= \nabla \Psi^{k+1}(\mathbf{u}^{k+1}; \mathbf{\zeta}^{k+1}, \mathbf{\tilde{\xi}}^{k+1}) - \mathbb{E} [\nabla \Psi^{k+1}(\mathbf{u}^{k+1}; \mathbf{\zeta}^{k+1}, \mathbf{\tilde{\xi}}^{k+1}) | \tilde{\mathcal{F}}_k ] \\
                &\quad + (1 - \beta_{k+1})(e^{k} + \mathbb{E} [\nabla \Psi^{k}(\mathbf{u}^{k}; \mathbf{\zeta}^{k}, \mathbf{\tilde{\xi}}^{k}) | \mathcal{F}_{k-1}]- \nabla \Psi^{k+1}(\mathbf{u}^{k}; \mathbf{\zeta}^{k+1}, \mathbf{\tilde{\xi}}^{k+1})) \\
                &= \beta_{k+1} \left(\nabla \Psi^{k+1}(\mathbf{u}^{k+1}; \mathbf{\zeta}^{k+1}, \mathbf{\tilde{\xi}}^{k+1} -\mathbb{E}[ \nabla \Psi^{k+1}(\mathbf{u}^{k+1}; \mathbf{\zeta}^{k+1}, \mathbf{\tilde{\xi}}^{k+1}) | \tilde{\mathcal{F}}_k ] \right) \\
                &\quad + (1 - \beta_{k+1})\left(\nabla \Psi^{k+1}(\mathbf{u}^{k+1}; \mathbf{\zeta}^{k+1}, \mathbf{\tilde{\xi}}^{k+1})
                - \mathbb{E}[ \nabla \Psi^{k+1} (\mathbf{u}^{k+1}; \mathbf{\zeta}^{k+1}, \mathbf{\tilde{\xi}}^{k+1}) | \tilde{\mathcal{F}}_k ] \right.\\
                &\quad \left.+ e^{k} + \mathbb{E} [\nabla \Psi^{k}(\mathbf{u}^{k}; \mathbf{\zeta}^{k}, \mathbf{\tilde{\xi}}^{k}) | \mathcal{F}_{k-1}] - \nabla \Psi^{k+1}(\mathbf{u}^{k}; \mathbf{\zeta}^{k+1}, \mathbf{\tilde{\xi}}^{k+1}) \right) \\
                &= \beta_{k+1} \Delta^{k+1} + (1- \beta_{k+1}) e^{k} + (1 - \beta_{k+1})(\Delta^{k+1} - \tilde \Delta^{k}) .
            \end{aligned}
        \end{equation}
        It follows from~\eqref{eq: bound of variance of mathcal G} and Lemma~\ref{lem: bound of  first order oracle of Gk} that 
        $\mathbb{E}[\|\Delta^{k+1} \|^2| \widetilde{\mathcal{F}}_k] = \mathbb{V}[\nabla^k(\mathbf{u}^{k+1}; \mathbf{\zeta}^{k+1}, \mathbf{\tilde{\xi}}^k] = \frac{\sigma_{f}^2}{ r_{k+1}} +  c_1^2 (\sigma_{\mathcal{G}}^k)^2$, where $\sigma_{\mathcal{G}}^k$ is defined in~\eqref{eq: sigma of mathcal G}. 
        Since $e^{k}$ is a constant term conditioned on $\widetilde{\mathcal{F}}_{k}$
        and $\mathbb{E} [ \| \Delta^{k+1} \| | \widetilde{\mathcal{F}}_k] = 0$,
       ~\eqref{eq: proof of bound of error in variance reduction 1} implies that 
        \begin{equation}
        \begin{aligned}
            \label{eq: proof of bound of error in variance reduction 2}
            &\quad \mathbb{E} [\| e^{k+1}\|^2| \widetilde{\mathcal{F}}_{k}] \\ 
            &= \mathbb{E} [ \| \beta_{k+1} \Delta^{k+1} + (1 - \beta_{k+1})(\Delta^{k+1} - \tilde \Delta^{k}) \|^ 2| \widetilde{\mathcal{F}}_{k} ] + (1 - \beta_{k+1})^2 \| e^{k} \|^2  \\
            &\leq 2 \beta_{k+1}^2 \mathbb{E}[\| \Delta^{k+1} \|^2 | \widetilde{\mathcal{F}}_{k}]
             + 2 (1 - \beta_{k+1})^2 \mathbb{E}[\|\Delta^{k+1} - \tilde \Delta^{k} \|^2 | \widetilde{\mathcal{F}}_{k}]  + (1 - \beta_{k+1})^2 \| e^{k} \|^2  \\
              &\leq {2} \beta_{k+1}^2 \left( \frac{\sigma_{f}^2}{r_{k+1}} + c_1^2 (\sigma_{\mathcal{G}}^k)^2\right)
             + 2 \mathbb{E}[\|\Delta^{k+1} - \tilde \Delta^{k} \|^2 | \widetilde{\mathcal{F}}_{k}] + (1 - \beta_{k+1})^2 \| e^{k} \|^2  \\
             &\leq {2} \beta_{k+1}^2 \left( \frac{\sigma_{f}^2}{r_{k+1}} + c_1^2 (\sigma_{\mathcal{G}}^k)^2\right)
             + 2 (1 - \beta_{k+1})^2 \mathbb{E}[\|\Delta^{k+1} - \tilde \Delta^{k} \|^2 | \widetilde{\mathcal{F}}_{k}] + (1 - \beta_{k+1})^2 \| e^{k} \|^2  .
        \end{aligned}  
        \end{equation}
        Now we come to handle the term $\mathbb{E} [ \| \Delta^{k+1} - \tilde \Delta^{k} \|^2 | \tilde{\mathcal{F}}_{k}]$ .
        Let 
        $$
        \begin{aligned}
        \mathrm{R}^{k+1}(\mathbf{u}^{k+1}, \mathbf{u}^{k}; \mathbf{\zeta}^{k+1},\mathbf{ \tilde{\xi}}^k) &= \nabla 
        \Psi^{k+1}(\mathbf{u}^{k+1}; \mathbf{\zeta}^{k+1}, \mathbf{\tilde{\xi}}^k) -\nabla \Psi^{k+1}(\mathbf{u}^{k}; \mathbf{\zeta}^{k+1}, \mathbf{\tilde{\xi}}^k), \\
        \mathrm{R}^{k+1}(\mathbf{u}^{k+1}, \mathbf{u}^{k}) &= \nabla 
        \Psi(\mathbf{u}^{k+1}) - \Psi(\mathbf{u}^{k}).
          \end{aligned}
          $$
          be the difference of the gradients between the points $\mathbf{u}^{k+1}$ and $\mathbf{u}^{k}$, respectively.
    The averaged Lipschitz property of $\nabla \Psi^{k+1}(\mathbf{u}; \mathbf{\zeta}, \mathbf{\tilde{\xi}})$
    and the Lipschitz property of $\nabla \Psi(\mathbf{u})$ give
        $\mathbb{E} [\| \mathrm{R}^{k+1}(\mathbf{u}^{k+1}, \mathbf{u}^{k}; \mathbf{\zeta}^{k+1}, \mathbf{\tilde{\xi}}^{k+1}) \|^2 | \tilde{\mathcal{F}}_k ]\leq 
        (L_{\Psi}')^2 \mathbb{E} [\|\mathbf{u}^{k+1} - \mathbf{u}^{k}\|^2 | \tilde{\mathcal{F}}_k] $ and 
        $\| \mathrm{R}(\mathbf{u}^{k+1}, \mathbf{u}^{k}) \| \leq 
        L_{\Psi} \|\mathbf{u}^{k+1} - \mathbf{u}^{k}\|$.
        Then it holds that 
        $$
        \begin{aligned}
        \Delta^{k+1} - \tilde \Delta^{k} &=
        \mathrm{R}^{k+1}(\mathbf{u}^{k+1}, \mathbf{u}^{k}; \mathbf{\zeta}^{k+1}, \mathbf{\tilde{\xi}}^{k+1}) - \mathrm{R}^{k}(\mathbf{u}^{k+1}, \mathbf{u}^{k})\\
        &\quad +  \mathbb{E} [\nabla \Psi^{k+1} (\mathbf{u}^{k+1}; \mathbf{\zeta}^{k+1}, \mathbf{\tilde{\xi}}^{k+1}) | \widetilde{\mathcal{F}}_k] 
         - \nabla \Psi^{k+1}(\mathbf{u}^{k+1}) \\
         &\quad - \mathbb{E} [\nabla \Psi^{k}(\mathbf{u}^{k}; \mathbf{\zeta}^{k}, \mathbf{\tilde{\xi}}^{k}) | \widetilde{\mathcal{F}}_{k-1}] + \Psi^{k}(\mathbf{u}^k),
        \end{aligned}
        $$
        and 
        \begin{equation}
            \label{eq: proof of bound of error in variance reduction 3}
        \begin{aligned}
           \mathbb{E}[\|\Delta^{k+1} - \tilde \Delta^{k} \|^2 | \widetilde{\mathcal{F}}_{k}] 
            &\leq 2 \mathbb{E}[\|\mathrm{R}^{k+1}(\mathbf{u}^{k+1}, \mathbf{u}^{k}; \mathbf{\zeta}^{k+1}, \mathbf{\tilde{\xi}}^k) - \mathrm{R}^{k+1}(\mathbf{u}^{k+1}, \mathbf{u}^{k}) \|^2 |\widetilde{\mathcal{F}}_{k} ] \\
            &\quad + 4 \| \mathbb{E} [ \nabla \Psi^{k+1}(\mathbf{u}^{k+1}; \mathbf{\zeta}^{k+1}, \mathbf{\tilde{\xi}}^k) - \nabla \Psi(\mathbf{u}^{k}) \| | \widetilde{\mathcal{F}}_k] \\
            &\quad + 4 \| \mathbb{E} [ \nabla \Psi^{k} (\mathbf{u}^{k}; \mathbf{\zeta}^{k}, \mathbf{\tilde{\xi}}^{k}) |  \widetilde{\mathcal{F}}_{k}] -  \nabla \Psi (\mathbf{u}^{k}) \|^2   \\
            &\leq \frac{4\left((L_{\Psi}')^2 +  (L_{\Psi}^2)^2\right)}{r_{k+1}} \mathbb{E}[\| \mathbf{u}^{k+1} - \mathbf{u}^{k} \|^2 | \widetilde{\mathcal{F}}_{k}] + \frac{4 c_1^2}{r_{k+1}}  \big((\epsilon_{\mathcal{G}}^k)^2 + (\epsilon_{\mathcal{G}}^{k+1})^2 \big).
        \end{aligned}
        \end{equation}
        The last inequality follows from Lemma~\ref{lem: bound of  first order oracle of Gk} and
       ~\eqref{eq: relationship between bias of Psi and G}.
        Combining~\eqref{eq: proof of bound of error in variance reduction 2} and~\eqref{eq: proof of bound of error in variance reduction 3} gives 
       ~\eqref{eq: bound of error in variance reduction}. This completes the proof.
    \qed
    Using the above two lemmas, we can derive the convergence of the variance reduction.
    \begin{theorem}
    \label{thm: decrease of variance reduction}
    The sequence $\{\mathbf{u}^{k+1}\}$ generated by Algorithm~\ref{alg: variance reduce} satisfies
    \begin{equation}
    \label{eq: decrease of variance reduction}
          \begin{aligned}
             &\quad 
             \frac{1}{\sum_{k=0}^{K-1} \alpha_k}
             \sum_{k=0}^{K-1} \mathbb{E}\left[ \frac{1}{\alpha_k}\|\mathbf{u}^{k+1} - \mathbf{u}^k\|^2 \right] 
             \leq  \frac{1}{\sum_{k=0}^{K-1} \alpha_k} \mathbb{E} [\Psi(\mathbf{u}^0) + \theta_0 \|e^0\|^2 - \Psi(\mathbf{u}^K) - \theta_K \|e^K\|^2] \\
             &\quad \quad  +  \frac{1}{\sum_{k=0}^{K-1} \alpha_k} \sum_{k=0}^K \left\{\alpha_kc^2 \mathbb{E} [ (\epsilon_{\mathcal{G}}^k)^2]  +
             2 \frac{\theta}{\alpha_{k}} \beta_{k+1}^2 \big( \frac{\sigma_{f}^2}{r_{k+1}} + c_1^2 \mathbb{E} [ (\sigma_{\mathcal{G}}^k)^2 ]\big) \right.\\
             &\quad \quad \left.
             + \frac{8c_1^2\theta}{\alpha_k r_{k+1}}\big(\mathbb{E}[(\epsilon_{\mathcal{G}}^k)^2] +\mathbb{E}[ (\epsilon_{\mathcal{G}}^{k+1})^2 ] \big)
             \right\}.
          \end{aligned}
    \end{equation}
    with  $\theta = \frac{1}{
    64 ((L_{\Psi}')^2 + L_{\Psi}^2)
    }$, $\alpha_k \leq \frac{1}{8L_{\Psi}}$, 
    $
       \beta_{k} \geq 1 - \sqrt{\frac{\frac{\theta}{\alpha_k} - \alpha_k}{8 }}
    $.
    \end{theorem}
    
    {\it Proof}
        Consider a merit function as $\Psi(\mathbf{u}^{k}) + \theta_k \|e^k\|^2$, where $\theta_k$ satisfies that 
       \begin{subequations}
        \begin{align}
            \alpha_k + 8 (1-\beta_{k+1})^2\theta_{k+1} - \theta_k \leq 0 ,    \label{eq: condition on theta 1} \\
            - \frac{1}{2 \alpha_k} + \frac{L_{\Psi}}{2} + 8\theta_{k+1}\frac{(L_{\Psi}')^2 + L_{\Psi}^2}{r_{k+1}} &\leq - \frac{1}{4\alpha_k}\label{eq: condition on theta 2} .
        \end{align}
        \end{subequations}
        Considering the reduction of the merit function, we have
        \begin{equation}
        \label{eq: proof of  decrease of variance reduction 1}
        \begin{aligned}
           &\quad \mathbb{E}[ \Psi(\mathbf{u}^{k+1}) + \theta_{k+1} \|e^{k+1}\|^2 - \Psi(\mathbf{u}^k) - \theta_{k} \|e^{k}\|^2| \widetilde{\mathcal{F}}_{k} ] \\
           &\leq \mathbb{E}\left[  -\frac{1}{2} (\frac{1}{\alpha_k} - L_{\Psi}) \|  \mathbf{u}^{k+1} - \mathbf{u}^{k} \|^2  + \alpha_k \|e^k\|^2  + \alpha_k c_1^2 (\epsilon_{\mathcal{G}}^k)^2 + \theta_{k+1} \|e^{k+1}\|^2 \right.\\
           &\quad \left. - \theta_k \|e^{k}\|^2 | \widetilde{\mathcal{F}}_{k} \right] \text{(by~\eqref{eq: decrease in variance reduction})}\\
          &\leq\mathbb{E}\left[ -\frac{1}{2} (\frac{1}{\alpha_k} - L_{\Psi}) \|  \mathbf{u}^{k+1} - \mathbf{u}^{k} \|^2  + \alpha_k \|e^k\|^2  + \alpha_k c_1^2 (\epsilon_{\mathcal{G}}^k)^2 \right. \\
          &\quad + \theta_{k+1} \left(
          2 \beta_{k+1}^2 \big( \frac{\sigma_{f}^2}{r_{k+1}} + c_1^2 (\sigma_{\mathcal{G}}^k)^2 \big) 
            + 8 (1 - \beta_{k+1})^2 \|e^{k} \|^2 \right.  \\
            &\quad +\left.\left. 8 \frac{(L_{\Psi}')^2 + L_{\Psi}^2}{r_{k+1}} \|\mathbf{u}^{k+1} - \mathbf{u}^k \|^2   + \frac{8c^2}{r_{k+1}} \big((\epsilon_{\mathcal{G}}^k)^2 + (\epsilon_{\mathcal{G}}^{k+1})^2 \big) \right) - \theta_k \|e^{k}\|^2 | \tilde{\mathcal{F}}_{k} \right]
            \quad \text{(by~\eqref{eq: bound of error in variance reduction})}
            \\
            &\leq\mathbb{E}\left[  \left(- \frac{1}{2 \alpha_k} + \frac{L_{\Psi}}{2} + 8\theta_{k+1}\frac{(L_{\Psi}')^2 + L_{\Psi}^2}{r_{k+1}} \right) \|  \mathbf{u}^{k+1} - \mathbf{u}^{k} \|^2  + \alpha_k c_1^2 (\epsilon_{\mathcal{G}}^k)^2  \right. \\
            &\quad + \left. \theta_{k+1} \left(
          2 \beta_{k+1}^2 \big( \frac{\sigma_{f}^2}{r_{k+1}} + c_1^2 (\sigma_{\mathcal{G}}^k)^2 \big)  + \frac{8c^2}{r_{k+1}} \big((\epsilon_{\mathcal{G}}^k)^2 + (\epsilon_{\mathcal{G}}^{k+1})^2 \big) \right)  | \widetilde{\mathcal{F}}_{k} \right]
            \quad \text{(by~\eqref{eq: condition on theta 1})}
            \\
            &\leq\mathbb{E}\left[ - \frac{1}{4\alpha_k} \|  \mathbf{u}^{k+1} - \mathbf{u}^{k} \|^2  + \alpha_k c_1^2 (\epsilon_{\mathcal{G}}^k)^2 \right. \\
            &\quad \left.+ \theta_{k+1} \left(
          2 \beta_{k+1}^2 \big( \frac{\sigma_{f}^2}{r_{k+1}} + c_1^2 (\sigma_{\mathcal{G}}^k)^2 \big)  + \frac{8c^2}{r_{k+1}} \big((\epsilon_{\mathcal{G}}^k)^2 + (\epsilon_{\mathcal{G}}^{k+1})^2 \big) \right)  | \widetilde{\mathcal{F}}_{k} \right] .
            \quad \text{(by~\eqref{eq: condition on theta 2})}
        \end{aligned}
    \end{equation}
    To ensure the conditions~\eqref{eq: condition on theta 1} and~\eqref{eq: condition on theta 2}, we take 
    \begin{equation}
        \label{eq: theta}
       \theta_{k+1}  = \frac{\theta}{\alpha_k} \quad \text{with}~\theta = \frac{1}{
    64 ((L_{\Psi}')^2 + L_{\Psi}^2)
       }.
    \end{equation}
    and let $\alpha_k \leq \frac{1}{8L_{\Psi}}$, 
    $
       \beta_{k} \geq 1 - \sqrt{\frac{\theta_k - \alpha_k}{8 }}
    $.
     Then~\eqref{eq: proof of  decrease of variance reduction 1} is simplified to 
    \begin{equation}
        \label{eq: proof of  decrease of variance reduction 2}
        \begin{aligned}
           &\quad \mathbb{E}[ \Psi(\mathbf{u}^{k+1}) + \theta_{k+1} \|e^{k+1}\|^2 - \Psi(\mathbf{u}^k) - \theta_{k} \|e^{k}\|^2| \widetilde{\mathcal{F}}_{k} ] \\
            &\leq\mathbb{E}\left[ - \frac{1}{4\alpha_k} \|  \mathbf{u}^{k+1} - \mathbf{u}^{k} \|^2  + \alpha_k c_1^2 (\epsilon_{\mathcal{G}}^k)^2  \right.\\
            &\quad \left.+ \frac{\theta}{\alpha_k} \big(
          2 \beta_{k+1}^2 \big( \frac{\sigma_{f}^2}{r_{k+1}} + c_1^2 (\sigma_{\mathcal{G}}^k)^2 \big)  + \frac{8c^2}{r_{k+1}} \big((\epsilon_{\mathcal{G}}^k)^2 + (\epsilon_{\mathcal{G}}^{k+1})^2 \big) \big)  \widetilde{\mathcal{F}}_{k} \right] .
        \end{aligned}
    \end{equation}
    Summing up the above inequality from $k=0$ to $K-1$, taking expectation over all the random variables 
    and multiplying $\frac{1}{\sum_{k=0}^{K-1} \alpha_k}$ on both sides, we have~\eqref{eq: decrease of variance reduction}.
    This completes the proof.
    \qed

    \begin{corollary}
          \label{cor: rate of variance reduction}
          Take constant sample sizes as $r_k = r, q_k = q, s_k = s$, $\beta_{k} = \beta \alpha_k^2$ with $\beta = \mathcal{O}(\alpha^{-2})$, and $\alpha_k = \alpha (k+1)^{-\frac{1}{3}}$, where $\alpha = \mathcal{O}(L_\Psi^{-1}) = \mathcal{O}(c^{-1})$ satisfying $\alpha_k \leq \frac{1}{2L_{\Psi}}$ as required in Theorem~\ref{thm: decrease of variance reduction}. Then the sequence $\{\mathbf{u}^{k}\}$ generated by Algorithm~\ref{alg: variance reduce} satisfies
          \begin{equation}
                \begin{aligned}
                   \frac{1}{\sum_{k=0}^{K-1} \alpha_k}\sum_{k=0}^{K-1} \mathbb{E}\left[ \frac{1}{\alpha_k}\|\mathbf{u}^{k+1} - \mathbf{u}^k\|^2 \right] 
                   &\leq\widetilde{\mathcal{O}} \left(\frac{c}{K^{\frac{2}{3}}} + \frac{c_1^2}{q} + \frac{c_1^2}{s}
                  + \frac{1}{K^{\frac{2}{3}}r} + \frac{K^{\frac{2}{3}}}{r} (\frac{1}{q} + \frac{1}{s})
                   \right) .
                 \end{aligned}
          \end{equation}
       \end{corollary}
       {\it Proof}
       From Corollary~\ref{cor: Lipschitz of Phi} and Lemma~\ref{lem: bound of bias of mathcal G and variance of mathcal G},~\ref{lem: bound of bias} we know 
       $L_{\Psi} = \mathcal{O}(c)$, $L_{\Psi}' = \mathcal{O}(c + \frac{c_1}{s})$,
       $\mathbb{E} [(\epsilon_{\mathcal{G}}^k)^2] \leq \widetilde{\mathcal{O}}(\frac{1}{q_k} + \frac{1}{s_k})$ ,
       $\mathbb{E} [(\sigma_{\mathcal{G}}^k)^2] \leq \widetilde{\mathcal{O}}(\frac{1}{q_k} + \frac{1}{s_k})$ and
       $\mathbb{E}[\|b^k\|^2] \leq \widetilde{\mathcal{O}}(\frac{1}{r_k} + \frac{c_1^2}{q_k} + \frac{c_1^2}{s_k})
       $. 
       Besides, it follows from~\eqref{eq: theta} that  $\theta = \mathcal{O}(L_\Psi^{-2}) = \mathcal{O}(c^{-2})$. By substituting these settings into~\eqref{eq: decrease of variance reduction}, we have
       \begin{equation}
             \begin{aligned}
                &\quad \frac{1}{\sum_{k=0}^{K-1} \alpha_k}\sum_{k=0}^{K-1} \mathbb{E}\left[ \frac{1}{\alpha_k}\|\mathbf{u}^{k+1} - \mathbf{u}^k\|^2 \right]  \\
                &\leq \widetilde{\mathcal{O}} \left( c K^{-\frac{2}{3}} \left( 
                1 + \alpha c_1^2 K^{\frac{2}{3}} \widetilde{\mathcal{O}}(\epsilon_{\mathcal{G}}^2) + \theta \beta^2 \alpha^3  \widetilde{\mathcal{O}}(r^{-1} + c_1^2(\sigma_{\mathcal{G}}^k )^2) + \frac{\theta c_1^2}{\alpha r} K^{\frac{4}{3}}
                \widetilde{\mathcal{O}}(\epsilon_{\mathcal{G}}^2)
                \right)
                \right) \\
                &= \widetilde{\mathcal{O}} \left( c K^{-\frac{2}{3}} \left( 
                1 + \frac{c_1^2}{c} K^{\frac{2}{3}} \widetilde{\mathcal{O}}(\epsilon_{\mathcal{G}}^2) + c^{-1}  \widetilde{\mathcal{O}}(r^{-1} + c_1^2(\sigma_{\mathcal{G}}^k )^2) + c^{-1} r^{-1} K^{\frac{4}{3}}
                \widetilde{\mathcal{O}}(\epsilon_{\mathcal{G}}^2)
                \right)
                \right) \\
                &\leq \widetilde{\mathcal{O}} \left( c K^{-\frac{2}{3}} \left( 
                1 + \frac{c_1^2}{c} K^{\frac{2}{3}}( \frac{1}{q} + \frac{1}{s}) + c^{-1} (\frac{1}{r} + \frac{c_1^2}{q} + \frac{c_1^2}{s}) + c^{-1} r^{-1}K^{\frac{4}{3}} (\frac{1}{q} + \frac{1}{s})
                 \right)
                \right) \\
                &=\widetilde{\mathcal{O}} \left(\frac{c}{K^{\frac{2}{3}}} + \frac{c_1^2}{q} + \frac{c_1^2}{s}
               + \frac{1}{K^{\frac{2}{3}}}(\frac{1}{ r} + \frac{c_1^2}{q} + \frac{c_1^2}{s}) + \frac{K^{\frac{2}{3}}}{r} (\frac{1}{q} + \frac{1}{s})
                \right) 
                \\
                &=\widetilde{\mathcal{O}} \left(\frac{c}{K^{\frac{2}{3}}} + \frac{c_1^2}{q} + \frac{c_1^2}{s}
               + \frac{1}{K^{\frac{2}{3}}r} + \frac{K^{\frac{2}{3}}}{r} (\frac{1}{q} + \frac{1}{s})
                \right) .
              \end{aligned}
       \end{equation}
       This completes the proof.
    \qed
       \begin{remark}
       Further take  $K \sim \mathcal{O}(c^{1.5}\epsilon^{-1.5}  )$, $r \sim \mathcal{O}(1)$, $q \sim \mathcal{O}(c_1^2 \epsilon^{-1})$, $s \sim \mathcal{O}(c_1^2 \epsilon^{-1})$, then the right side is
       $\widetilde{\mathcal{O}}(\epsilon)$.  The sample complexity on $\xi$ is $\sum_{k=0}^{K-1} {s_k} + \sum_{k=0}^{K-1} {q_k} = (s+ q)K = \widetilde{\mathcal{O}}(c^{1.5} c_1^2\epsilon^{-2.5})$ and the sample complexity on $\zeta$ is $\sum_{k=0}^{K-1} {r_k} = rK = \widetilde{\mathcal{O}}(c^{1.5} \epsilon^{-1.5})$ .
       \end{remark}

\bibliographystyle{plainnat} 
\bibliography{ref}

\begin{thebibliography}{31}
\providecommand{\natexlab}[1]{#1}
\providecommand{\url}[1]{\texttt{#1}}
\expandafter\ifx\csname urlstyle\endcsname\relax
  \providecommand{\doi}[1]{doi: #1}\else
  \providecommand{\doi}{doi: \begingroup \urlstyle{rm}\Url}\fi

\bibitem[Alpaydin and Alimoglu(1996)]{alpaydin1996pen}
Ethem Alpaydin and Fevzi Alimoglu.
\newblock Pen-based recognition of handwritten digits.
\newblock \emph{(No Title)}, 1996.

\bibitem[Bennett et~al.(2006)Bennett, Hu, Ji, Kunapuli, and
  Pang]{bennett2006model}
Kristin~P Bennett, Jing Hu, Xiaoyun Ji, Gautam Kunapuli, and Jong-Shi Pang.
\newblock Model selection via bilevel optimization.
\newblock In \emph{The 2006 IEEE International Joint Conference on Neural
  Network Proceedings}, pages 1922--1929. IEEE, 2006.

\bibitem[Bonnans and Shapiro(2013)]{bonnans2013perturbation}
J~Fr{\'e}d{\'e}ric Bonnans and Alexander Shapiro.
\newblock \emph{Perturbation analysis of optimization problems}.
\newblock Springer Science \& Business Media, 2013.

\bibitem[Chen et~al.(2021)Chen, Sun, and Yin]{chen2021tighter}
Tianyi Chen, Yuejiao Sun, and Wotao Yin.
\newblock Tighter analysis of alternating stochastic gradient method for
  stochastic nested problems.
\newblock \emph{arXiv preprint arXiv:2106.13781}, 2021.

\bibitem[Cutkosky and Orabona(2019)]{cutkosky2019momentum}
Ashok Cutkosky and Francesco Orabona.
\newblock Momentum-based variance reduction in non-convex sgd.
\newblock \emph{Advances in neural information processing systems}, 32, 2019.

\bibitem[Domahidi et~al.(2013)Domahidi, Chu, and Boyd]{domahidi2013ecos}
Alexander Domahidi, Eric Chu, and Stephen Boyd.
\newblock Ecos: An socp solver for embedded systems.
\newblock In \emph{2013 European control conference (ECC)}, pages 3071--3076.
  IEEE, 2013.

\bibitem[Franceschi et~al.(2018)Franceschi, Frasconi, Salzo, Grazzi, and
  Pontil]{franceschi2018bilevel}
Luca Franceschi, Paolo Frasconi, Saverio Salzo, Riccardo Grazzi, and
  Massimiliano Pontil.
\newblock Bilevel programming for hyperparameter optimization and
  meta-learning.
\newblock In \emph{International conference on machine learning}, pages
  1568--1577. PMLR, 2018.

\bibitem[Gao et~al.(2023)Gao, Ye, Yin, Zeng, and Zhang]{gao2023moreau}
Lucy~L Gao, Jane~J Ye, Haian Yin, Shangzhi Zeng, and Jin Zhang.
\newblock {Moreau envelope based difference-of-weakly-convex reformulation and
  algorithm for bilevel programs}.
\newblock \emph{arXiv preprint arXiv:2306.16761}, 2023.

\bibitem[Grimmer et~al.(2023)Grimmer, Lu, Worah, and
  Mirrokni]{grimmer2023landscape}
Benjamin Grimmer, Haihao Lu, Pratik Worah, and Vahab Mirrokni.
\newblock The landscape of the proximal point method for nonconvex--nonconcave
  minimax optimization.
\newblock \emph{Mathematical Programming}, 201\penalty0 (1):\penalty0 373--407,
  2023.

\bibitem[Hong et~al.(2023)Hong, Wai, Wang, and Yang]{hong2023two}
Mingyi Hong, Hoi-To Wai, Zhaoran Wang, and Zhuoran Yang.
\newblock A two-timescale stochastic algorithm framework for bilevel
  optimization: Complexity analysis and application to actor-critic.
\newblock \emph{SIAM Journal on Optimization}, 33\penalty0 (1):\penalty0
  147--180, 2023.

\bibitem[Ji et~al.(2020)Ji, Yang, and Liang]{ji2020provably}
Kaiyi Ji, Junjie Yang, and Yingbin Liang.
\newblock Provably faster algorithms for bilevel optimization and applications
  to meta-learning.
\newblock \emph{Neural Information Processing Systems}, 2020.

\bibitem[Jiang et~al.(2024{\natexlab{a}})Jiang, Xiao, Tenorio, Real-Rojas,
  Marques, and Chen]{jiang2024primal}
Liuyuan Jiang, Quan Xiao, Victor~M Tenorio, Fernando Real-Rojas, Antonio
  Marques, and Tianyi Chen.
\newblock A primal-dual-assisted penalty approach to bilevel optimization with
  coupled constraints.
\newblock \emph{arXiv preprint arXiv:2406.10148}, 2024{\natexlab{a}}.

\bibitem[Jiang et~al.(2024{\natexlab{b}})Jiang, Li, Hong, and
  Zhang]{jiang2024barrier}
Xiaotian Jiang, Jiaxiang Li, Mingyi Hong, and Shuzhong Zhang.
\newblock A barrier function approach for bilevel optimization with coupled
  lower-level constraints: Formulation, approximation and algorithms.
\newblock \emph{arXiv preprint arXiv:2410.10670}, 2024{\natexlab{b}}.

\bibitem[Kang et~al.(2023)Kang, Jung, Jeoung, Hong, and Hong]{kang2023bi}
Hyuna Kang, Seunghoon Jung, Jaewon Jeoung, Juwon Hong, and Taehoon Hong.
\newblock A bi-level reinforcement learning model for optimal scheduling and
  planning of battery energy storage considering uncertainty in the
  energy-sharing community.
\newblock \emph{Sustainable Cities and Society}, 94:\penalty0 104538, 2023.

\bibitem[Khanduri et~al.(2023)Khanduri, Tsaknakis, Zhang, Liu, Liu, Zhang, and
  Hong]{khanduri2023linearly}
Prashant Khanduri, Ioannis Tsaknakis, Yihua Zhang, Jia Liu, Sijia Liu, Jiawei
  Zhang, and Mingyi Hong.
\newblock Linearly constrained bilevel optimization: A smoothed implicit
  gradient approach.
\newblock In \emph{International Conference on Machine Learning}, pages
  16291--16325. PMLR, 2023.

\bibitem[Lu and Mei(2024)]{lu2024first}
Zhaosong Lu and Sanyou Mei.
\newblock First-order penalty methods for bilevel optimization.
\newblock \emph{SIAM Journal on Optimization}, 34\penalty0 (2):\penalty0
  1937--1969, 2024.

\bibitem[MacKay et~al.(2019)MacKay, Vicol, Lorraine, Duvenaud, and
  Grosse]{mackay2019self}
Matthew MacKay, Paul Vicol, Jon Lorraine, David Duvenaud, and Roger Grosse.
\newblock Self-tuning networks: Bilevel optimization of hyperparameters using
  structured best-response functions.
\newblock \emph{arXiv preprint arXiv:1903.03088}, 2019.

\bibitem[Nocedal and Wright(1999)]{nocedal1999numerical}
Jorge Nocedal and Stephen~J Wright.
\newblock \emph{Numerical optimization}.
\newblock Springer, 1999.

\bibitem[Qin et~al.(2023)Qin, Song, and Jiang]{qin2023bi}
Xiaorong Qin, Xinhang Song, and Shuqiang Jiang.
\newblock Bi-level meta-learning for few-shot domain generalization.
\newblock In \emph{Proceedings of the IEEE/CVF Conference on Computer Vision
  and Pattern Recognition}, pages 15900--15910, 2023.

\bibitem[Rockafellar and Wets(2009)]{rockafellar2009variational}
R~Tyrrell Rockafellar and Roger J-B Wets.
\newblock \emph{Variational analysis}, volume 317.
\newblock Springer Science \& Business Media, 2009.

\bibitem[Shen and Chen(2023)]{shen2023penalty}
Han Shen and Tianyi Chen.
\newblock On penalty-based bilevel gradient descent method.
\newblock In \emph{International Conference on Machine Learning}, pages
  30992--31015. PMLR, 2023.

\bibitem[Shen et~al.(2024)Shen, Yang, and Chen]{shen2024principled}
Han Shen, Zhuoran Yang, and Tianyi Chen.
\newblock Principled penalty-based methods for bilevel reinforcement learning
  and rlhf.
\newblock \emph{arXiv preprint arXiv:2402.06886}, 2024.

\bibitem[Sinha et~al.(2020)Sinha, Khandait, and Mohanty]{sinha2020gradient}
Ankur Sinha, Tanmay Khandait, and Raja Mohanty.
\newblock A gradient-based bilevel optimization approach for tuning
  hyperparameters in machine learning.
\newblock \emph{arXiv preprint arXiv:2007.11022}, 2020.

\bibitem[Stadie et~al.(2020)Stadie, Zhang, and Ba]{stadie2020learning}
Bradly Stadie, Lunjun Zhang, and Jimmy Ba.
\newblock Learning intrinsic rewards as a bi-level optimization problem.
\newblock In \emph{Conference on Uncertainty in Artificial Intelligence}, pages
  111--120. PMLR, 2020.

\bibitem[Tsaknakis et~al.(2022)Tsaknakis, Khanduri, and
  Hong]{tsaknakis2022implicit}
Ioannis Tsaknakis, Prashant Khanduri, and Mingyi Hong.
\newblock An implicit gradient-type method for linearly constrained bilevel
  problems.
\newblock In \emph{ICASSP 2022-2022 IEEE International Conference on Acoustics,
  Speech and Signal Processing (ICASSP)}, pages 5438--5442. IEEE, 2022.

\bibitem[Tsaknakis et~al.(2023)Tsaknakis, Khanduri, and
  Hong]{tsaknakis2023implicit}
Ioannis Tsaknakis, Prashant Khanduri, and Mingyi Hong.
\newblock An implicit gradient method for constrained bilevel problems using
  barrier approximation.
\newblock In \emph{ICASSP 2023-2023 IEEE International Conference on Acoustics,
  Speech and Signal Processing (ICASSP)}, pages 1--5. IEEE, 2023.

\bibitem[Xu and Zhu(2023)]{xu2023efficient}
Siyuan Xu and Minghui Zhu.
\newblock Efficient gradient approximation method for constrained bilevel
  optimization.
\newblock In \emph{Proceedings of the AAAI Conference on Artificial
  Intelligence}, volume~37, pages 12509--12517, 2023.

\bibitem[Yang et~al.(2024)Yang, Gao, and Yuan]{yang2024bilevel}
Yan Yang, Bin Gao, and Ya-xiang Yuan.
\newblock Bilevel reinforcement learning via the development of hyper-gradient
  without lower-level convexity.
\newblock \emph{arXiv preprint arXiv:2405.19697}, 2024.

\bibitem[Yao et~al.(2024{\natexlab{a}})Yao, Yin, Zeng, and
  Zhang]{yao2024overcoming}
Wei Yao, Haian Yin, Shangzhi Zeng, and Jin Zhang.
\newblock Overcoming lower-level constraints in bilevel optimization: A novel
  approach with regularized gap functions.
\newblock \emph{arXiv preprint arXiv:2406.01992}, 2024{\natexlab{a}}.

\bibitem[Yao et~al.(2024{\natexlab{b}})Yao, Yu, Zeng, and
  Zhang]{yao2024constrained}
Wei Yao, Chengming Yu, Shangzhi Zeng, and Jin Zhang.
\newblock {Constrained bi-level optimization: Proximal Lagrangian value
  function approach and Hessian-free algorithm}.
\newblock \emph{arXiv preprint arXiv:2401.16164}, 2024{\natexlab{b}}.

\bibitem[Zhu et~al.(2020)Zhu, Li, Wu, Zhao, Ding, and Shi]{zhu2020personalized}
Hancheng Zhu, Leida Li, Jinjian Wu, Sicheng Zhao, Guiguang Ding, and Guangming
  Shi.
\newblock Personalized image aesthetics assessment via meta-learning with
  bilevel gradient optimization.
\newblock \emph{IEEE Transactions on Cybernetics}, 52\penalty0 (3):\penalty0
  1798--1811, 2020.

\end{thebibliography}

\end{document}